%% file: toricmaps.tex
\documentclass[a4paper,10pt]{article}

\usepackage{color}
\newenvironment{red}{\color{red}}{}
\newcommand{\bred}{\begin{red}}
\newcommand{\ered}{\end{red}}
\usepackage{preamble}
\usepackage[all,line]{xy}

 \usepackage{vmargin}
 \setmargrb{1in}{1in}{1in}{1in} 

\newdir{smallhead}{
:(15,1)\dir^{>}
*:(15,-1)\dir_{>}}

\newdir{multihead}{
\dir{smallhead}
*!(0,1.1)\dir{smallhead}
*!(0,-1.1)\dir{smallhead}
*{\POS(-6,0); \POS(0, 1.1) **\crv{~**\dir{-} (-4.5,1.1)}}
*{\POS(-6,0); \POS(0,-1.1) **\crv{~**\dir{-} (-4.5,-1.1)}}
}

\newcommand{\multito}{
\
\begin{xy}
  =(0,0) "SOURCE";
  =(9,0) "TARGET";
  {\ar@{multihead} "SOURCE"; "TARGET"}
\end{xy}
\
}

\newcommand{\multimapsto}{
\
\begin{xy}
  =(0,0) "SOURCE";
  =(9,0) "TARGET";
  {\ar@{|-multihead} "SOURCE"; "TARGET"}
\end{xy}
\
}

\newcommand{\xymultito}[1]{\ar@{multihead}[#1]}

\newcommand{\latticemarker}{\bullet}

\newcommand{\xylattice}{
\POS( 0, 0)*{\latticemarker};
\POS( 1, 0)*{\latticemarker};
\POS( 0, 1)*{\latticemarker};
\POS( 1, 1)*{\latticemarker};
\POS(-1, 0)*{\latticemarker};
\POS( 0,-1)*{\latticemarker};
\POS(-1,-1)*{\latticemarker};
\POS(-1, 1)*{\latticemarker};
\POS( 1,-1)*{\latticemarker};
\POS( 2, 2)*{\latticemarker};
\POS( 2, 1)*{\latticemarker};
\POS( 2, 0)*{\latticemarker};
\POS( 2,-1)*{\latticemarker};
\POS( 2,-2)*{\latticemarker};
\POS( 1,-2)*{\latticemarker};
\POS( 0,-2)*{\latticemarker};
\POS(-1,-2)*{\latticemarker};
\POS(-2,-2)*{\latticemarker};
\POS(-2,-1)*{\latticemarker};
\POS(-2, 0)*{\latticemarker};
\POS(-2, 1)*{\latticemarker};
\POS(-2, 2)*{\latticemarker};
\POS(-1, 2)*{\latticemarker};
\POS( 0, 2)*{\latticemarker};
\POS( 1, 2)*{\latticemarker};
}

\newcommand{\xybiglattice}{
\xylattice
\POS( -2, 3)*{\latticemarker};
\POS( -1, 3)*{\latticemarker};
\POS(  0, 3)*{\latticemarker};
\POS(  1, 3)*{\latticemarker};
\POS(  2, 3)*{\latticemarker};
}

\newcommand{\halfsingularity}{
      \xylattice
      {\ar@{-}(0,0);(1, 1)};
      {\ar@{-}(0,0);(2, 2)};
      {\ar@{-}(0,0);(1,-1)};
      {\ar@{-}(0,0);(2,-2)};
}

\newcommand{\affinepieceofresolutionofhalfsingularity}{
      \xylattice
      {\ar@{-}(0,0);(1, 1)};
      {\ar@{-}(0,0);(2, 2)};
      {\ar@{-}(0,0);(2, 0)};
}

\newcommand{\fanoffakewps}{
      \xylattice
      {\ar@{-}(0,0);(1,2)};
      {\ar@{-}(0,0);(-2,-1)};
      {\ar@{-}(0,0);(2,-2)};
}

\newcommand{\gavina}{
      \xylattice
      {\ar@{-}(0,0);(1,2)};
      {\ar@{-}(0,0);(-1,-2)};
      {\ar@{-}(0,0);(2,0)};
      {\ar@{-}(0,0);(0,2)};
      \POS(0, -3)*{\text{fan of $Y$}};
}

\newcommand{\gavinb}{
      \xybiglattice
      {\ar@{-}(0,0);(1,3)};
      {\ar@{-}(0,0);(-1,-2)};
      {\ar@{-}(0,0);(2,2)};
      \POS(0, -3)*{\text{fan of $X$}};
}

\newcommand{\fanoffakewpswithextraray}{
\fanoffakewps
      {\ar@{-}(0,0);(2,0)};
}

\newcommand{\makefan}[1]{
\[
  \begin{xy}
    ( 0, 0);( 6, 0):;
     #1
   \end{xy}
\]
}

\newcommand{\functionfield}[1]{\CC(#1)}

\newcommand{\phibar}{{\overline{\varphi}}}
\newcommand{\phireg}{{\varphi_\mathrm{reg}}}
\newcommand{\phireginv}{\varphi_\mathrm{reg}^{-1}}
\def\Reg{\operatorname{Reg}}
\def\homog#1{{#1^{\mathrm{hgs}}}}
\def\dualcone#1{{#1}^{\vee}}
\def\agree#1#2{{\operatorname{Agr}(#1, #2)}}
\def\agr{\mathfrak{agr}}
\def\irrel#1{\operatorname{Irrel}(#1)}
\def\coxmonolatt#1{TM(#1)}
\def\monoposcone#1{TM[#1]}
\newcommand{\raylattice}[1]{R_{#1}}
\newcommand{\crl}{S}         
\newcommand{\mrl}{\Gamma}    
\newcommand{\cox}[1]{\crl[#1]}
\newcommand{\coxfield}[1]{\crl(#1)}

\newcommand{\mapring}[1]{\mrl(#1)}

\newcommand{\mvsa}{\gamma}   
\newcommand{\mvsb}{\delta}   
\DeclareMathOperator{\Bl}{Bl}    
\newcommand\Iff{{\ \Longleftrightarrow\ }}
\newcommand{\zerolocus}[1]{\set{#1 = 0}}
\newcommand{\starfan}[1]{\operatorname{Star}(#1)}
\newcommand{\quofan}{\Sigma_{Y(\sigma)}}

\newcommand{\DiagramForImage}[7]{
\ensuremath{
   \xymatrix @C 0pc{
      \text{$#1$}
    &
    \ccO_X(X_g) \ar@{^(->}[d]^{\deg 0 \text{ part}}  && \text{$ \qquad \quad  \qquad $}&&
                \ccO_Y(V)  \ar@{^(->}[d]_{\deg 0 \text{ part}} \ar[llll]_{\varphi^*} & \text{$#5$}
   \\
   \text{$#2$}
    &
    \cox X [\inv{g}] &&&&
                \cox Y [ \inv{\upsilon}]  &
   \text{$#6$}
   \\
   \text{$#3$}
    &
    \CC[\Reg \Phi]\ar@{^(->}[u] \ar@{^(->}[d] &&&&  \cox Y  \ar@{^(->}[u] \ar[dllll]_{\Phi^*} &
   \text{$#7$}
   \\
   \text{$#4$}
    &
    \mapring \Phi
     }
}
}

\title{Maps of toric varieties in Cox coordinates}
\author{Gavin Brown \and \JaBu}
\date{}

\begin{document}

\renewcommand{\theenumi}{\textnormal{(\roman{enumi})}}
\renewcommand{\labelenumi}{\theenumi}
\newcommand{\Dind}{D_{\mathrm{ind}}}
\newcommand{\Dirrel}{D_{\mathrm{irrel}}}

\maketitle


\begin{abstract}
The Cox ring provides a coordinate system on a toric variety analogous
to the homogeneous coordinate ring of projective space. Rational maps
between projective spaces are described using polynomials in the
coordinate ring, and we generalise this to toric varieties, providing a
unified description of arbitrary rational maps between toric varieties
in terms of their Cox coordinates. Introducing formal roots of polynomials
is necessary even in the simplest examples.
\end{abstract}


\medskip
{\footnotesize
\noindent\textbf{AMS Mathematical Subject Classification 2010:}
Primary: 14M25; Secondary:  14E05, 14Q99;}

\setcounter{tocdepth}{2}
\tableofcontents

\input{toricmapsintro}

\input{toricmapsbackground}

\input{toricmapsmultivalued}

\input{toricmapsdescriptions}


\bibliography{toricmaps}

\bibliographystyle{alpha}
\vfill
\noindent
\begin{tabular}{p{0.45\textwidth}p{0.48\textwidth}}
Gavin Brown                  & \JaBu                               \\
School of Mathematics        & Institute of Mathematics of the     \\
Loughborough University      & Polish Academy of Sciences          \\
LE11 3TU,                    & ul.~\'Sniadeckich 8,  P.O.~Box 21,  \\
United Kingdom               & 00-956 Warszawa, Poland             \\
\verb|G.D.Brown@lboro.ac.uk| & \verb|jabu@mimuw.edu.pl|            \\
\end{tabular}
\vfill

\end{document}

%% file: toricmapsintro.tex

\section{Introduction}
\label{sect:intro}

This paper describes maps between toric varieties in terms of
Cox coordinates, that is, using the usual generators of the Cox rings
of the source and target.
The results are not confined to maps that preserve the toric structures,
but to arbitrary rational maps of such varieties.

Any rational map between two projective spaces can
be lifted to a morphism between their {\em Cox covers},
their affine GIT covering spaces:
it is described by a sequence of homogeneous polynomials of the same degree.
To generalise this to maps between any toric varieties,
we need descriptions which also use roots of polynomials, and so
we cannot hope to lift the maps to morphisms, or even to rational maps,
between covering spaces: instead, we consider multi-valued maps
like $x\mapsto\pm\sqrt{x}$, which we denote by $x\multimapsto\sqrt{x}$
to emphasise that it is not a map in the usual sense.

The use of radical expressions to define maps is well established in
some toric and orbifold contexts; writing weighted blow ups of cyclic quotient
singularities, for example. The radicals define a map on the orbifold cover,
and this paper generalises such calculations to all rational maps
of toric varieties.

The Cox ring literature has several treatments of the
functors both of toric varieties, initiated by Cox \cite{cox_functor}
and generalised by Kajiwara \cite{kajiwara}, who also uses radical
expressions explicitly, and of more general varieties satisfying
certain finiteness conditions by Berchtold and Hausen \cite{berchtold_hausen}.
There is also an approach by Berchtold and Hausen \cite{berchtold_hausen04},
Theorem~9.2, using {\em bunches} of cones.
These are mainly concerned with morphisms,
whereas the treatment here considers all rational
maps of all toric varieties, and uses all Weil divisors rather than
(sufficiently many) Cartier divisors.

Our main result, stated more precisely as Theorem~\ref{thm:main} and
in final form as Theorem~\ref{thm:existence_of_complete_descriptions}, is this.
Let $\varphi\colon X\dashrightarrow Y$ be a rational map between
toric varieties (not necessarily respecting their toric structures).
Then there is a `multi-valued map' $\Phi\colon\CC^m\multito\CC^n$
between the Cox covers of $X$ and $Y$ which is defined
using radical expressions in the Cox coordinates of $X$ and has
the following properties:
\begin{description}
\item[Evaluation at points:]
if $\varphi$ is defined at $x\in X$ and $\xi\in\CC^m$ is an expression
for $x=[\xi]$ in Cox coordinates, then $\varphi(x) = [\Phi(\xi)]\in Y$.
\item[Pullback of divisors:]
If $D=(f)$ is a Cartier divisor on $Y$, where $f$ lies in the Cox ring
$\cox Y$ of $Y$, then the divisors $\varphi^*(D)$ and $(\Phi^*f)$ on $X$
agree on the open subset where $\varphi$ is regular.
\end{description}
These are the two essential properties of the {\em complete description} $\Phi$,
refined as properties~A--\ref{item:complete_description_coordinates}
in Sections~\ref{sect:descriptions}--\ref{sect:DEF},
but it has other good features: for example, it allows easy computation
of the image and preimage of subschemes under~$\varphi$ (\S\ref{sect:image}--\ref{sect:preimage}).

In the rest of this introduction we present some examples
and briefly survey enough of the Cox ring approach to toric geometry
to be able to state the main result more precisely.
Section~\ref{sect:simple_ring_extensions} explains a class
of radical extensions of rings which we apply in
Section~\ref{sect:multival} to make a basic theory of multi-valued maps.
These two sections are the technical heart of the paper. The practical
theory for describing maps that we build on this is natural, but it
succeeds because we work in carefully controlled extensions
of the Cox rings when writing the coordinates of maps.
In Section~\ref{sect:descriptions}, we say what it means for
a multi-valued map to describe a rational map between toric varieties,
and we prove the main Theorem~\ref{thm:existence_of_complete_descriptions}
on the existence of a complete description $\Phi$.
Section~\ref{sect:properties_of_descriptions} explains the composition
of descriptions and the computation of images and preimages.

We work over the complex numbers $\CC$. The foundational
aspects of toric geometry \cite{kempf_knutsen_mumford_saintdonat} work
over any field,
but our presentation relies on Cox's construction \cite{cox_homogeneous},
and that is given over $\CC$.

\subsection{Motivating examples} \label{sect:Motivating_example}

\subsubsection{A line on a quadric}
\label{example:coordinate_axis_on_P112}

A weighted projective space $\PP(\fromto{a_1}{a_n})$ 
  is a (usually) singular algebraic variety obtained  as the quotient 
  $(\CC^n \setminus \set{0})/\CC^*$, where the action of $\CC^*$ has weights $(\fromto{a_1}{a_n})$,
  that is: 
  \[
    t \cdot (\fromto{y_1}{y_n}) = (\fromto{{t}^{a_1} y_1}{{t}^{a_n} y_n}).
  \]
This is completely analogous to the case of an ordinary projective space $\PP^{n-1} = \PP(\fromto{1}{1})$
  and just as in the case of $\PP^{n-1}$ we can consider $\fromto{y_1}{y_n}$
  to be the \emph{homogeneous coordinates} on $\PP(\fromto{a_1}{a_n})$.

Consider the weighted projective space $\PP(1,1,2)$ with homogeneous
coordinates $y_1$, $y_2$, $y_3$. The coordinate axis $\Gamma=(y_2) \subset \PP(1,1,2)$
is a smooth rational curve $\Gamma\cong\PP^1$. In coordinates $x_1,x_2$
on $\PP^1$, we can describe the embedding
$\PP^1 \to \Gamma\subset\PP(1,1,2)$ by
 \[
   [x_1,x_2] \multimapsto [\sqrt{x_1},0, x_2]
 \]
(see \S\ref{sect:multi_valued_sections} for our formal definition
of $\sqrt{x_1}$). Multiplying through by $\sqrt{x_1}$ with the
given weights $(1,1,2)$ gives an alternative:
 \[
   [x_1,x_2] \mapsto [x_1,0, x_1 x_2].
 \]
 We discuss two benefits of the first. 

The first issue is to calculate images of points. For instance, to
see the image of the point $[0,1] \in \PP^1$ using the first
description, we immediately compute $[0,0,1]$.
With the second description, we are in trouble, because the description
of the map evaluates to $[0,0,0]$ and so does not help.
The square root
is not too bad. The image of the point $[1,0] \in \PP^1$
computed by the first description is either $[1,0,0]$ or $[-1,0,0]$
depending on which root we take; but these are the same point in
$\PP(1,1,2)$, so either expression is fine.

The second issue is to pull back divisors. For instance, to pull
back a Cartier divisor from the linear system of $\ccO_{\PP(1,1,2)}(2)$,
we would like simply to substitute the definining equations of the map.
For example, suppose we pull back $y_3=0$. Clearly, this coordinate
axis meets $\Gamma$ transversely in one point $[1,0,0]$. Using the
first description, we pullback the function $y_3$ to get the function $x_2$,
whose vanishing locus on $\PP^1$ is exactly $[1,0]$ as we would like.
The second description, however, is not good enough in this respect either:
the naive pull back is $x_1 x_2$.

\subsubsection{Weighted blow ups: the affine $\frac{1}{2}(1,1)$ singularitity}
\label{sect:example_weighted_blow_up}

In the first example, the square root merely simplified some
calculations. Now we give an example where it is unavoidable.
Consider the simplest singular toric variety $Y$: the affine
$\half(1,1)$ singularity, that is, the quotient of $\CC^2$ by
$\ZZ/2$ acting by
\[
  (y_1,y_2) \mapsto (-y_1, -y_2).
\]
Let $X$ be an affine piece of its resolution, $X= \CC^2 \subset \Bl_{[0,0]} Y$.
In fan terminology this corresponds to the following embedding of cones:
\makefan{
      \affinepieceofresolutionofhalfsingularity
      {\ar (3,0); (5,0)};
      (8,0);(9,0):;
      \halfsingularity}
The map $\varphi\colon X \to Y$ as a map of affine varieties,
\[
  \varphi\colon \Spec \CC[x_1, x_2] \to \Spec \CC[{y_1}^2, y_1 y_2, {y_2}^2],
\]
corresponds, via the dual map of cones, to the affine coordinate
ring homomorphism:
\begin{align*}
  \varphi^*\colon\CC[{y_1}^2, y_1 y_2, {y_2}^2] & \to \CC[x_1,x_2] \\
\text{sending}\quad  {y_1}^2 \mapsto x_1, \quad
  y_1 y_2 & \mapsto x_1 x_2  \quad\text{and}\quad
  {y_2}^2 \mapsto x_1 {x_2}^2.
\end{align*}
Therefore if we hope to extend $\varphi^*$ to the full
Cox ring $\cox{Y}=\CC[y_1,y_2]$
\[
  \xymatrix{
    \varphi^*\colon\CC[{y_1}^2, y_1 y_2, {y_2}^2] \ar[r] \ar@{^{(}->}[d]
    	& \CC[x_1,x_2]\ar@{^{(}->}[d]  \\
    \Phi^*\colon \CC[y_1, y_2] \ar[r] & \text{some new ring}
  }
\]
we need a map $\Phi^*$ doing either
\[
\begin{array}{rl}
  y_1 &\multimapsto \sqrt{x_1}\\
  y_2 &\multimapsto x_2\sqrt{x_1}
\end{array}
\quad\mathrm{or}\quad
\begin{array}{rl}
  y_1 &\multimapsto -\sqrt{x_1}\\
  y_2 &\multimapsto -x_2\sqrt{x_1}.
\end{array}
\]
Introducing the square roots is necessary for such a description.
We are allowed to choose either square root of $x_1$, but we must
make the choice only once: having picked the root of $x_1$ for the
first coordinate, the root of $x_1$ used in the second coordinate must
be the same.

\subsubsection{Fake weighted projective space}

Descriptions of maps that require roots also arise for maps
between projective toric varieties. Let $\Sigma_Y$ be the fan
\makefan{\fanoffakewps}
and $Y$ the associated toric variety; this is the simplest example
of a fake projective space (see \cite{wero_fps}, \cite{kasprzyk_fps})
and is the quotient of $\PP^2$ by $\ZZ/3$ acting with weights $(2,1,0)$.

Let $X$ be a weighted blow up of any of the $3$ singular points of $Y$,
for example given by the fan
\makefan{\fanoffakewpswithextraray}
Then every description of the blow up map $X \to Y$ will involve at
least $3$rd roots of polynomials.
For instance, if we encode the actions defining $X$ and $Y$ as
\[
\begin{pmatrix}
1 & 1 & 1 & 0\\
2 & 1 & 0 & -3
\end{pmatrix}
\quad\text{and}\quad
\begin{pmatrix}
1 & 1 & 1\\
2/3 & 1/3 & 0
\end{pmatrix}
\]
in coordinates $x_1,\dots,x_4$ on $X$, $y_1,\dots,y_3$ on $Y$---treating
the second row of the weights of $Y$ as the homogeneity imposed by the finite
$\ZZ/3$ action---then the map is defined by
\[
[x_1,\dots,x_4] \multimapsto [ x_1\sqrt[3]{x_4}^2, x_2\sqrt[3]{x_4}, x_3].
\]
(The second row of the grading matrix of $Y$ only permits scaling by
cube roots of unity, so it cannot be used to eliminate the radical here;
the notation is slightly clumsy.)

\subsubsection{Ideals of subvarieties of toric varieties}
\label{sect:ideals_of_subvarieties_of_toric_varieties}

The use of Cox rings to describe subschemes of toric varieties
includes a small, well-known catch \cite[Thm 3.7]{cox_homogeneous}:
significantly different ideals can determine the same subscheme.
This problem arises when considering maps too.
Consider $X = \PP^2$ and an action of $\ZZ/2$ on $X$ with weights
$(0,0,1)$. The quotient of $X$ by $\ZZ/2$ is $Y = \PP(1,1,2)$, and,
in coordinates, the quotient map $\varphi\colon X\rightarrow Y$ is
\[
  [x_1,x_2,x_3] \mapsto [x_1, x_2, {x_3}^2].
\]
This description of $\varphi$ has the two properties mentioned at
the outset (it is well defined on every point of $X$ and Cartier
divisors can be pulled back by simple substitution), but there
is still a difficulty when calculating the preimage of subschemes
in $Y$. For instance, in coordinates $y_i$ on $Y$, the subschemes
$B_1$ and $B_2$ of $Y$ defined, respectively, by the ideals
\[
  \langle y_1 \rangle
\text{ and }
  \langle {y_1}^2, y_1 y_2 \rangle
\]
are equal and both reduced, but the ideal defining $B_2$ is not
radical even though both ideals are saturated at the irrelevant
maximal ideal. Local calculations
show that the preimage subscheme $A=\varphi^{-1}(B_2)$ is non-reduced
and equal to scheme defined by $\langle {x_1}^2, x_1 x_2 \rangle$. On the other hand,
if we pullback the defining equations of $B_1$ we get the reduced
scheme $A' = (x_1)$. Although $A$ and $A'$ are certainly not
equal as schemes, their scheme structures are equal on the
preimage of smooth locus of $Y$. This is the best we can
hope for and is explained generally in
Theorem~\ref{thm:image_under_description_agrees}.

\subsubsection{Reading toric birational maps from complete descriptions}
\label{sect:example_toric_birational}

Let $X=\PP(1,1,2)$ with Cox coordinates $x_1, x_2, x_3$ and $Y$ be the toric
variety with Cox coordinates $y_1, y_2, y_3, y_4$, bi-grading given by the matrix
of weights
\[
\begin{pmatrix}
1&2&0&-1\\0&0&1&1
\end{pmatrix}
\quad
\text{and irrelevant ideal}
\quad
B_Y = (y_1,y_2) \cap (y_3,y_4).
\]
Suppose $Y$ and $X$ are described by fans in a common lattice $N=\ZZ^2$
as follows.
\makefan{
      \gavina
      {\ar (3,0); (5,0)};
      (8,0);(9,0):;
      \gavinb}
The implicit birational map between $X$ and $Y$ is
\begin{align*}
\varphi\colon X &\dashrightarrow Y
&\text{and}\quad
\psi\colon Y &\dashrightarrow X\\
[x_1, x_2, x_3] &\longmapsto [x_1,x_1 x_2, x_1 x_3, x_1 x_2] &
[y_1, y_2, y_3, y_4] &\longmapsto [{y_1}^2 y_4, y_2 y_4, y_1 y_2 y_3 y_4].
\end{align*}

The geometry of this birational equivalence is evident in the fans but we
hope to read it from equation descriptions.
It is better seen using complete descriptions, which
we can make easily from the original monomial
descriptions: the question is simply how much can we cancel.
For $\varphi$ we
can use the first grading of $Y$ to remove a $\sqrt{x_1}$ factor,
and then the second grading to remove a further $x_1$: thus
\begin{align*}
[x_1,x_1 x_2, x_1 x_3, x_1 x_2] &\text{ becomes }
[\sqrt{x_1}, x_2, x_1 x_3, x_1 x_2 \sqrt{x_1}]\\
 & \text{ which in turn becomes }
[\sqrt{x_1}, x_2, x_3, x_2 \sqrt{x_1}].
\end{align*}
Similarly we can modify the description of $\psi$, so the result is
\begin{align*}
\varphi\colon [x_1, x_2, x_3] &\multimapsto [\sqrt{x_1}, x_2, x_3, x_2 \sqrt{x_1}]\\
\psi\colon [y_1, y_2, y_3, y_4] &\multimapsto [y_1^2 \sqrt{y_4},y_2 \sqrt{y_4}, y_1 y_2 y_3].
\end{align*}
Many features of the birational geometry are now clear.
The map $\varphi$ is not defined at the three $0$-strata of $X$,
while $\psi$ is not defined on the
$0$-strata $(1,0,1,0)$ and $(0,1,1,0)$ in $Y$.
The coordinate loci $(x_1)$ and $(x_2)$ in $X$ are contracted,
and similarly $(y_1)$, $(y_2)$ and $(y_4)$ in $Y$ are contracted.
Furthermore, comparing with the weighted blow ups above, we see that
$(x_1)$ and $(y_4)$ are contracted as $\frac{1}2(1,1)$ exceptional divisors,
while $(y_1)$ is a $(2,1)$ weighted blow up of a smooth point,
and $(y_2)$ and $(x_2)$ are ordinary (smooth) blow ups of smooth points.

\subsubsection{Spaces of maps}
\label{sect:using_descriptions}

For two toric varieties $X$ and $Y$ and a map $\alpha\colon\Pic Y \to \Pic X$
we can use the structure of descriptions to classify all the regular maps
$\varphi\colon X \to Y$ for which $\varphi^*=\alpha$, precisely because
our results on descriptions apply to all maps.
We illustrate by computing all maps from a toric del Pezzo surface
to a certain weighted projective $5$-space; the conclusion is that
the map is unique and toric up to coordinate choice.
For brevity, we will assume that the image of $\varphi$ is not contained
in any toric stratum of $Y$, not even after a change of coordinates on $Y$.

Let $X = \FF_1$, simply $\PP^2$ blown up in a single point,
a del Pezzo surface of degree $8$.
Thus
\[
\cox{X} = \CC[x_1, x_2, x_3, x_4]
\quad\text{graded by $\ZZ^2$ with gradings}\quad
 \begin{pmatrix}
   1&1&1&0\\
   0&0&1&1
 \end{pmatrix}
\]
and irrelevant ideal $(x_1,x_2) \cap (x_3,x_4)$.
Consider also $Y = \PP(1,1,1,2,2,2)$:
\[
\cox{Y} = \CC[y_1, y_2, y_3, y_4, y_5, y_6]
\quad\text{graded by $\ZZ$ with gradings}\quad
 \begin{pmatrix}
   1&1&1&2&2&2
 \end{pmatrix}.
\]

For the demonstration, we assume that the regular map $\varphi\colon X \to Y$
pulls back a divisor in $\ccO_Y(2)$ (the ample generator of $\Pic Y$)
to an anticanonical divisor of $X$ in $\ccO(3,2)$. We claim that
the map $\varphi$ is unique up to changes of coordinates on
$X$ and $Y$. (This is analogous to nondegenerate quadratic maps
$\PP^1\rightarrow\PP^2$ being the usual conic in the right coordinates.)
We use the results of this paper, in particular that complete descriptions
exist (Theorem~\ref{thm:existence_of_complete_descriptions})
and satisfy properties A--D (Definitions~\ref{defin:homogrel}
and \ref{defin:properties_DEFG} and Propositions~\ref{prop:description_satisfy_homog_and_relev}
and \ref{prop:Y_Q_factorial_then_max_agr_not_necessary}).

Let $\Phi\colon \CC^4 \multito \CC^6$ be a complete description of $\varphi$,
a well-defined expression of the map using rational functions and radicals as
above (with, loosely speaking, as much cancellation as possible already done).
We may assume that each component of $\Phi$ is of the form $p\cdot q^{1/r}$
for polynomials $p,q$ and some $r\in\NN$. (This holds in general
for regular maps by Corollary~\ref{cor:regularity_locus_from_complete_description}.)
By condition~\ref{item:complete_description_Cartier} (the pullback of Cartier divisors is given by $\Phi^*$ on the regular locus of $\varphi$)
the expressions
$(\Phi^*y_1)^2$,
$(\Phi^*y_2)^2$,
$(\Phi^*y_3)^2$,
$\Phi^*y_4$,
$\Phi^*y_5$ and
$\Phi^*y_6$ are rational forms.
Applying the homogeneity condition~\ref{item:homogeneity_cond_on degree_0}
(the usual homogeneity condition that rational functions pull back to
rational functions)
we see that $\frac{\Phi^*y_2}{\Phi^*y_1}$ and $\frac{\Phi^*y_3}{\Phi^*y_1}$
are rational functions.
Thus we can write
\begin{align*}
  \Phi\colon \CC^4 & \multito \CC^6\\
                   x & \multimapsto [f_1 \sqrt{g}, f_2 \sqrt{g}, f_3 \sqrt{g}, f_4, f_5, f_6]
\end{align*}
for polynomials $f_i, g \in \cox X$, and apply
condition~\ref{item:homogeneity_cond_on degree_0} once more to see
\begin{gather*}
     \deg f_1 = \deg {f_2} =\deg {f_3},\qquad
     \deg f_4 = \deg {f_5} =\deg {f_6},\\
     \text{and}\quad 2 \deg f_1 + \deg g = \deg f_4 = (3,2).
\end{gather*}
The last condition narrows the possibilities for the multidegree
of $f_1$:
\[
  \deg f_1 \in \set{(0,0), (1,0),(0,1),(1,1)}.
\]
But the linear systems in multidegrees $(0,0), (1,0),(0,1)$ are small,
and allowing the degree of $f_1$ to be any those would force
the $3$ sections $f_1,f_2, f_3$ to be linearly dependent.
A suitable coordinate change on $Y$ would then transform (at least)
one of $f_1, f_2, f_3$ to $0$, presenting the image of $\varphi$
inside some toric stratum, which is exactly what our simplifying
assumption forbids. So $\deg f_1 = (1,1)$.

So $\deg g =(1,0)$, and changing coordinates on $X$ we may assume $g = x_1$.
Also the $\CC$-linear span of $f_1,f_2, f_3$ is equal to the span
$x_1 x_4, x_2 x_4, x_3$, so changing coordinates on $Y$ we may assume
\[
   f_1 = x_1 x_4,\quad f_2 = x_2 x_4 \quad\text{and}\quad f_3 = x_3.
\]
The linear system $(3,2)$ is spanned by the nine monomials:
\[
  {x_1} {x_3}^2, {x_2} {x_3}^2,
  {x_1}^2  {x_3} {x_4}, {x_1} {x_2}  {x_3} {x_4},  {x_2}^2  {x_3} {x_4},
  {x_1}^3 {x_4}^2, {x_1}^2 {x_2} {x_4}^2, {x_1} {x_2}^2 {x_4}^2, {x_2}^3 {x_4}^2.
\]
However if any of $f_4, f_5, f_6$ contains any summand divisible by $x_1$
then we can change the coordinates on $Y$ to get rid of this summand.
For instance if $f_4 = {x_2} {x_3}^2 +  {x_1} {x_2}  {x_3} {x_4}$, then
$f_4 - (f_2 \sqrt{g}) (f_3 \sqrt{g}) = {x_2} {x_3}^2$.
Therefore we may assume $f_4, f_5, f_6$ are contained in span of
$
  {x_2} {x_3}^2,
  {x_2}^2  {x_3} {x_4},
  {x_2}^3 {x_4}^2,
$
and changing coordinates on $Y$ again we may assume
\[
   f_4 = {x_2} {x_3}^2,\quad
   f_5 = {x_2}^2  {x_3} {x_4} \quad\text{and}\quad
   f_6 = {x_2}^3 {x_4}^2.
\]
Thus every map $\varphi$ satisfying the assumptions can be written as
\begin{align*}
  \varphi \colon X & \longrightarrow Y\\
                 x & \multimapsto
        [x_1 x_4 \sqrt{x_1}, x_2 x_4 \sqrt{x_1}, x_3 \sqrt{x_1},
              {x_2} {x_3}^2,  {x_2}^2  {x_3} {x_4},  {x_2}^3 {x_4}^2]
\end{align*}
in some homogeneous coordinates on $X$ and $Y$.

\subsubsection{Multi-valued multi-linear systems}
\label{sect:mvmls}

It is worth noting that the homogeneity conditions of
Definition~\ref{defin:homogrel}
are more precise than simply arranging for the degrees of components
of a map being correct.

Let $X=\PP(1,1,2)$ with coordinates $x_1, x_2, x_3$ and
$Y=\PP(1,2,3)$ with coordinates $y_1, y_2, y_3$.
Let $f={x_1}^3-x_2 x_3$ and $\mvsa=\sqrt{f}$.
Then
\[
\Phi\colon (x_1, x_2, x_3) \multimapsto (\sqrt{x_1},x_2,\mvsa)
\]
has the correct degrees but nevertheless
fails to determine a rational map: indeed
\[
\Phi^*(y_2/{y_1}^2) = x_2/x_1
\quad\textrm{is nice, but}\quad
\Phi^*(y_3/{y_1}^3) = \sqrt{1 - \frac{x_2 x_3}{{x_1}^3}}
\]
is not a rational function on $X$ (which is what the homogeneity condition
requires; or, using the homogeneity condition \ref{item:homogeneity_cond_all}
and the language of Definition~\ref{def:multivalsec} instead,
$\Phi^*({y_1}^3+y_3)$ is not a homogeneous multi-valued section.)
Simply arranging for the correct homogeneous degrees is not the
full content of the homogeneity condition.
It is better thought of
as requiring all defining sections to be elements of a single
vector space of multi-valued sections together with its multiples.
If $\mvsa$ is the third coordinate, then the degree~3
sections defining the map must all have $\mvsa$ as their common
irrational part; formally speaking, this is the conclusion of
Proposition~\ref{prop:image_has_constant_degree_and_irrational_part}.

But they do not: $\Phi^*(y_1 y_2) = \sqrt{x_1}\cdot x_2$ has
irrational part $\sqrt{x_1}$ not equal to that of $\Phi^*(y_3)$.
Forcing $\Phi^*y_3 = \mvsa$ requires $\sqrt[r]{f}$ for $r=6$ and $4$
respectively as a factor
into the first two components; but then we can scale the entire
irrational part away in any case.

However, defining a (different) map as
\[
\Phi\colon (x_1, x_2, x_3) \multimapsto (\mvsa,{x_2}^3,\mvsa^3 + \mvsa x_1 x_3)
\]
is fine, since now
\[
\Phi^*(y_2/{y_1}^2) = {x_2}^3/f
\quad\textrm{and}\quad
\Phi^*(y_3/{y_1}^3) = 1 + (x_1 x_3 / f).
\]
(And, at least as a first test, $\Phi^*(y_1^3+y_3)$ is now $\mvsa \cdot (2f + x_1 x_3)$,
which is a homogeneous multi-valued section.)

If we regard a map to a weighted projective space as being determined
by a basis of a graded ring
$
V = \bigoplus_{d\in\NN} V_d
$
where each $V_d\subset \overline{\coxfield X}$ is a finite-dimensional
vector space consisting only of multi-valued
sections of degree $d/N$, for some fixed denominator $N\in\NN$,
then we must ensure that each $V_d$ has the same irrational part $\mvsa^d$,
for some $\mvsa\in\overline{\cox X}$.
In the corrected example, this holds:
\[
V_1 = \mvsa\cdot\CC,\quad
V_2 = \mvsa^2 \cdot \CC\left< 1, {x_2}^3/f \right>,\quad
V_3 = \mvsa^3 \cdot \CC\left< 1, {x_2}^3/f, x_1 x_3/f \right>
\]
and so on---the irrational parts of these spaces of sections are
visibly the same (up to the power that fixes their degree).

\subsection{Maps of toric varieties in Cox coordinates}

\subsubsection{Cox coordinates on toric varieties}
\label{sect:Cox_coordinates}

We review the standard elements of toric geometry that we use
throughout this paper, closely following three of the standard sources
\cite{cox_homogeneous}, \cite{danilov} and \cite{fulton},
without further comment or citation.
A toric variety $X$ of dimension $d$ is defined by a fan $\Sigma_X$
spanning a  (possibly strict) subspace of a $d$-dimensional lattice $N_X$.
The rays of $\Sigma_X$, which, by minor abuse of notation, we can
take as the primitive vectors $\rho_1,\dots,\rho_m$ on the
$1$-skeleton $\Sigma_X^{(1)}$, play two roles.
First, treating them as independent symbols,
they generate a new lattice $\raylattice{X}\cong\ZZ^m$,
\textbf{the ray lattice} of $X$, with chosen basis the $\rho_i$.
 The natural map $\rho_X\colon \raylattice{X}\rightarrow N_X$ sends
each symbol $\rho_i$ to the primitive vector. 
When considering the Cox quotient construction,
one usually assumes for convenience that $X$ has no torus factors, but
this is not necessary in our approach (see also \cite[\S5.1]{cox_book}).
If $X = X' \times (\CC^*)^k$, where $X'$ is has no torus factors,
then the fan $\Sigma_X$ spans a linear subspace
$\langle \Sigma_X \rangle\subset N_X \otimes \RR$ of codimension $k$.
We choose primitive lattice vectors $\fromto{\rho_{m+1}}{\rho_{m+k}}$ in $N_X$
such that the lattice $\langle \Sigma_X \rangle \cap N_X$ together with
$\fromto{\rho_{m+1}}{\rho_{m+k}}$ generate the lattice $N_X$.
These additional lattice vectors are called \textbf{virtual rays}
and they play the role of place holders for variables corresponding to
coordinates on $(\CC^*)^k$. The ray lattice is then extended to
$\raylattice{X}\cong\ZZ^{m+k}$ with the bigger basis $\fromto{\rho_{1}}{\rho_{m+k}}$,
and the map $\rho_X\colon \raylattice{X}\rightarrow N_X$ is extended
accordingly to take account of these virtual rays.

Second, we denote the elements of the basis dual to the $\rho_i$ in
$\raylattice{X}$ by $x_i$, and interprete them as the indeterminates
of a polynomial ring.
The ring the $x_i$ generate is the famous \textbf{Cox ring} $\cox X$ of
$X$, also known as its homogeneous, or total, coordinate ring.
It is graded by the divisor class group $\Cl(X)$.
The \textbf{irrelevant ideal}
$B_X\subset\cox X$ is defined by standard generators, one for each maximal
cone $\sigma\in\Sigma_X$, defined as $\mu_\sigma = \prod x_i$, where the
product is taken over those rays $\rho_i$ not contained in $\sigma$ (one sets
$B_X=\cox X$ if there is only one cone of maximal dimension).
Note that if $\rho_i$ is a virtual ray then
the monomial $\mu_{\sigma}$ is divisible by $x_i$ for every $\sigma$.

Thus $X=\CC\times\CC^*$, determined by a fan with a single ray
in $N_X = \ZZ^2$ as its unique maximal cone, has Cox ring
$\cox X=\CC[x_1,x_2]$ and irrelevant ideal $B_X = \left<x_2\right>$
(rather than $\cox X = B_X = \CC[x_1,x_2,1/x_2]$, for example), where the
variable $x_1$ corresponds to the 1-skeleton of the fan and $x_2$ 
to a virtual ray chosen arbitrarily to extend the rational span of the fan
to the entire $\ZZ^2$.

We also treat the $x_i$ in their own right, namely as a basis of
the lattice dual to $\raylattice{X}$,
the \textbf{Cox monomials lattice} $\coxmonolatt{X}$.
We write $\monoposcone{X}$ for the positive orthant in $\coxmonolatt{X}$.
The lattice $M_X$ of monomials, the dual of $N_X$,
embeds $M_X\hookrightarrow\coxmonolatt{X}$
as the dual map to $\rho_X$.

The {\bf Cox cover} of $X$ is defined to be $\Spec\cox X$; it is
isomorphic to $\CC^m$ with standard coordinates $x_i$, and we
usually write it as such with its heritage implicit.
The gradings describe the action of a group $G_X = \Hom(\Cl(X), \CC^*) \simeq T\oplus A$,
where $T\cong \Gm^d$ is an algebraic torus and $A$ is a finite
abelian group.
Cox proves, Theorem~2.1 of \cite{cox_homogeneous},
that $X$ is a quotient of $\CC^m$ by $G_X$ in the sense of GIT.
Indeed, there is a rational map $\pi_X\colon \CC^m\dashrightarrow X$
that is a morphism precisely on $\Reg\pi_X$, the complement of
the {\bf irrelevant locus} $\irrel{X}=V(B_X)\subset\CC^m$,
and is a categorical quotient there.
Thus one thinks of elements $\xi\in\CC^m$ as representative
coordinate expressions for their images $x=\pi_X(\xi)\in X$;
we also denote $\pi_X(\xi)$ by $[\xi]$.
These are the Cox coordinates on $X$ that we use systematically.

Denoting the field of fractions of $\cox X$ by $\coxfield{X}$,
the function field $\functionfield{X}$ of $X$ is naturally isomorphic
to the subfield of $\coxfield{X}$ of $G_X$-invariant functions.
We treat these as
being the rational functions on $\CC^m$ of degree~0, just as
for rational functions on projective space.
We refer to elements of $\cox X$ and $\coxfield{X}$ as {\bf polynomial
and rational sections} on $X$ respectively, rather than as functions.
We say that section $f\in \coxfield{X}$ is regular on
$U \subset X$ if $f$ is a regular function on
$\pi_X^{-1}(U) = \{ \xi\in\Reg\pi_X \mid \pi_X(\xi)\in U\}$.

The Cox ring has a more intrinsic definition.
Suppose in the first place that $X$ has no torus factors.
Then
\[
\cox X = \bigoplus H^0(X,D)
\]
where the sum is taken over the Weil class group $\Cl(X)$, with
$D$ being a representative Weil divisor in the particular class (chosen
systematically so that multiplication is defined automatically).
The natural isomorphism between these two descriptions
follows from the association of a Weil divisor $D_\rho$ to each
ray~$\rho$: $D_\rho$ is the irreducible divisor supported
on the image of $\set{x_\rho=0}\subset\CC^m$ in $X$, where $x_\rho$ is
the Cox coordinate corresponding to $\rho$.
In general, when $X = X' \times (\CC^*)^k$ with virtual rays
$\fromto{\rho_{m+1}}{\rho_{m+k}}$,
\[
\cox X [\fromto{x_{m+1}^{-1}}{x_{m+k}^{-1}}] =
 \bigoplus H^0(X,D)
\]
where
\[
 \cox X [\fromto{x_{m+1}^{-1}}{x_{m+k}^{-1}}]
= \CC[\fromto{x_1}{x_m}, \fromto{x_{m+1}, x_{m+1}^{-1}}{x_{m+k}, x_{m+k}^{-1}}].
\]
We take these
isomorphisms as implicit, so for each homogeneous rational section
$f\nolinebreak\in\nolinebreak\coxfield{X}$ there is a Weil divisor, denoted $(f)$.
The converse is also true and follows from the same isomorphism:
if $D$ is a Weil divisor on a toric variety $X$, then $D=(f)$
for some non-zero homogeneous function $f \in \coxfield{X}$.
Moreover, in the case $X$ has no torus factors,
$D$ is effective if and only if $f \in \cox X$;
if $X$ does have torus factors, the criterion is instead that
$f \in \cox X[\fromto{x_{m+1}^{-1}}{x_{m+k}^{-1}}]$.
In any case, if $D$ is effective,
then there exists non-zero homogeneous section $f \in \cox{X}$
such that $D = (f)$.
This association also obeys the natural calculus: $(f g) = (f) + (g)$.

Given $f \in \cox{X}$, as well as considering the divisor $(f)$ on $X$,
we will also consider the zero set of $f$ in the Cox cover $\CC^m$ of $X$.
To avoid confusion we will always denote this affine zero set by
$\zerolocus{f} \subset \CC^m$.

\subsubsection{The main results}
\label{sect:results}

The elementary examples of \S\ref{sect:Motivating_example} are
part of a general theory. The first result is that any rational
map between toric varieties has a description by radicals of
Cox coordinates.

\begin{thm}
\label{thm:main}
Let $X$ and $Y$ be toric varieties over $\CC$ with Cox rings
$\cox{X} = \CC[\fromto{x_1}{x_m}]$ and $\cox{Y} = \CC[\fromto{y_1}{y_n}]$
and corresponding Cox covers $\CC^m$ and $\CC^n$.

If $\varphi\colon X \dashrightarrow Y$ is a rational map, then there
are homogeneous rational sections $q_i \in \coxfield{X}$ and an expression
\[
    \Phi\colon [\fromto{x_1}{x_m}] \multimapsto
    \bigl[\fromto{\sqrt[r_1]{q_1}}{\sqrt[r_n]{q_n}}\bigr],
\]
which satisfies the following properties:
\begin{enumerate}
\item
If $\xi \in \CC^m$ and $\varphi$ is regular at $x=[\xi]$, then
$y=[\Phi(\xi)]$ is a well-defined point of $Y$ and $\varphi(x) = y$.
\item
If $D=(f)$ is a Cartier divisor on $Y$ whose suport does not contain the image of $\varphi$,
where $f\in \coxfield{Y}$,
then $\varphi^*D$ and $(\Phi^*f)$ are equal as divisors on $X$
when restricted to the regular locus of $\varphi$.
\item
If $A \subset X$ is a closed subscheme defined by a saturated ideal
$I_A \ideal \cox X$, then the image $\varphi(A) \subset Y$
is defined by the preimage under $\Phi^*$ of the span of
$I_A$ in some extension of $\cox X$.
\item
If $B \subset Y$ is a closed subscheme defined by an ideal
$I_B \ideal \cox Y$, then the preimage $\varphi^{-1}(B) \subset X$
is defined on $\varphi^{-1}(Y_0)$, by the ideal
$\langle \Phi^*(I_B) \rangle \cap \cox X$ of $\cox X$,
where $Y_0$ is the smooth locus of $Y$.
\end{enumerate}
\end{thm}

This statement needs some explanation.
In \S\ref{sect:agreement}, we explain what it means for an expression
$
    \Phi\colon [\fromto{x_1}{x_m}] \multimapsto
    [\fromto{\sqrt[r_1]{q_1}}{\sqrt[r_n]{q_n}}],
$
to be a description of a rational map $X \dashrightarrow Y$, and
Definition~\ref{defin:complete_agreement} specifies `complete descriptions'.
This theorem gathers some results for complete descriptions proved in
Theorems~\ref{thm:existence_of_complete_descriptions},
\ref{thm:image_under_description_agrees},
\ref{thm:preimage_under_description_agrees},
Proposition~\ref{prop:very_complete_description_Cartier} and their subsequent
comments and corollaries.
Those results are more general and detailed; the statements above are
special cases.
The statement on preimage above does not explain the extension (in fact, it is
simply a map ring $\mapring\Phi$ as discussed next), but the precise
details are in Corollary~\ref{cor:preimage_indep}.

Furthermore, care is needed when defining $\Phi(\xi)$.
Recall from \S\ref{sect:example_weighted_blow_up} that the root
of a polynomial can be chosen arbitrarily but only chosen once.
If the same root of the same polynomial occurs again in the
expression for $\Phi$ (even if not in an explicit form), then
we must use the root chosen before. We make this book-keeping
precise by introducing simple extensions of rings in
\S\ref{sect:ring_extensions} and map rings $\Gamma(\Phi)$ for $\Phi$ in
\S\ref{sect:map_ring}.
The point is that we work in extensions $\mapring{\Phi}$ of $\cox{X}$
containing the image of $\Phi^*$ which cannot be made arbitrarily;
the notion of `simple' extension assembles just enough conditions for
our purposes here. The ideal spans of the form
$\langle J \rangle$ appearing in the statement are taken inside
these $\mapring{\Phi}$.
Theorems~\ref{thm:image_under_description_agrees}
and \ref{thm:preimage_under_description_agrees} explain this precisely,
and the latter also
explains how to achieve the exact preimage over the singular locus.

\begin{rem}
The statement of the theorem might suggest that using the descriptions one is able to pullback the Weil divisors, 
  even those that are not $\QQ$-Cartier.
However, in the situation, when $Y$ is not $\QQ$-factorial,
  the complete descriptions as in Theorem~\ref{thm:main} are not unique,
  see Example~\ref{ex:complete_descriptions_are_not_unique}.
It is implicit in the statement, that the divisor $(\Phi^* f)$ does not depend on the choice of $\Phi$,
  whenever $(f)$ is Cartier on $Y$.
But when $(f)$ only defines a Weil divisor, then $(\Phi^* f)$ depends on $\Phi$.
\end{rem}

The second result gives a criterion for a radical expression like
$\Phi$ above to determine a rational map of toric varieties;
this is spelled out in
Theorem~\ref{thm:Phi_homog_and_relevant_is_a_description}.
\begin{thm}
Let $\Phi\colon \CC^m\multito\CC^n$ be a multi-valued map
between the Cox covers of toric varieties $X$ and $Y$.
If $\Phi$ satisfies the homogeneity and relevance conditions of
Definition~\ref{defin:homogrel}, then there is a unique rational
map $\varphi\colon X\dashrightarrow Y$ that $\Phi$ describes.
\end{thm}
In other words, subject only to natural conditions of
homogeneity with respect to all gradings and relevance
(and the precise specification of what
is allowed as a radical expression to define~$\Phi$), a sequence
of radical expressions in Cox coordinates does indeed determine a
rational map.

It was pointed out by an anonymous referee that the methods here should
apply more generally to the Mori dream spaces of Hu and Keel
\cite{hukeel} in a fairly natural way: the Cox ring of
a Mori dream space is
a quotient of a polynomial ring, so we can work with
coordinates (and so also multivalued coordinates) as usual.
However we have not checked the details of this: the relations in the Cox ring
add more relevance conditions and also relate the radical multivalued
expressions (which we keep independent by use of the simple extensions of
\S\ref{sect:ring_extensions}), so there is something to check.

\subsection*{Acknowledgements}

It is our pleasure to thank Janko B\"ohm for his patient listening
and many comments, Grzegorz Kapustka for suggesting examples and
David Cox and Paulo Lima-Filho for explaining and clarifying several points.
The authors were supported by EPSRC grant EP/E000258/1 at the University
of Kent, UK, during much of this work, and the second author was a
Marie Curie International Outgoing Fellow at University of Grenoble,
France for the completion. The paper was further revised under the project
``Secant varieties, computational complexity, and toric degenerations''
realised within the Homing Plus programme of Foundation for Polish Science,
cofinanced from European Union, Regional Development Fund.
Our thanks also to Joseph Landsberg and
Colleen Robles at Texas A\&M and John Cannon at The University of Sydney
for their hospitality and support during our collaborative visits.

%% file: toricmapsbackground.tex
\section{Simple extensions of rings}
\label{sect:simple_ring_extensions}

We review some material in the context of
multi-graded rings in \S\ref{sect:auxilliary}, then present
some field theory in \S\ref{sect:field-extensions}, and finally
give the key definition of simple extension of rings
in \S\ref{sect:ring_extensions}.

\subsection{Auxilliary algebra and geometry}
\label{sect:auxilliary}

\subsubsection{Homogeneous ideals}

 We outline standard points about ideals in rings
graded over finitely-generated abelian grading semigroups.
The cases we have in mind include the Cox ring $\cox X$ of a
toric variety $X$, extensions $\cox X [f^{-1}]$ for a
homogeneous polynomial $f$, quotients $\cox X / I$ for some
homogeneous ideal $I$ and combinations of these. Recall that
$\cox X$ has a distinguished ideal, the irrelevant ideal $B_X$.
In our applications, the grading group is $H=\Hom(G_X,\CC^*)$.
We write $H$ additively with identity $0\in H$; we often consider
elements of degree~$0$ in the rings above. 

\begin{defin}
Let $\crl$ be a graded ring. A homogeneous ideal $\gotp \ideal \crl$
is \textbf{homogeneously prime} if and only if whenever a
homogeneous $h\in \gotp$ factorises $h= fg$ with homogeneous
factors $f, g \in \crl$, then either $f\in \gotp $ or $g \in \gotp$.
\end{defin}

This notion is also called $G$-prime; see \cite{perling_toric_varieties_as_spectra},
Remark~3.20, or \cite{ludwig_prime_ideals} in an unrelated context.
A homogeneous ideal which is prime is homogeneously prime, but
the converse is not always true.

\begin{example}\label{example:homog_prime_not_prime}
Let $X$ be the affine $\half(1,1)$ singularity, so $\cox X = \CC[x_1, x_2]$
graded by $\ZZ/2$ as multiplication by $-1$ and $B_X=\cox X$.
The ideal generated by $x_1^2-1$ is homogeneously prime but not prime.
It determines an irreducible line $L\subset X$, but regarded
on the GIT cover $\CC^2$ it determines a disjoint union of two lines,
$x_1 = 1$ and $x_1 = -1$, the preimage $\inv{\pi_X} L$.
\end{example}

This is an important consideration throughout. When the grading
group is $\ZZ^k$ with no torsion, it is easy to see that the two concepts coincide.


\begin{prop}{\cite[Proposition~2.4]{cox_homogeneous}}
For homogeneously prime ideal $I \ideal \cox X$ there exists
a unique irreducible subvariety $V(I) \subset X$,
such that a section $f\in\cox X$ vanish identically on $V(I)$
if and only if $f\in I$.

Conversely, for every irreducible subvariety $V \subset X$,
there exists a homogeneously prime ideal $I(V) \ideal \cox X$
contained in the irrelevant ideal $B$ such that $V(I(V)) = V$.
\end{prop}

\begin{defin}
\label{defin:homog_ideal}
Let $\crl$ be a graded ring and $I \ideal \crl$ an ideal. The
\textbf{homogenisation of $I$} is the biggest homogeneous ideal
$\homog I$ contained in $I$.
\end{defin}
It follows that $\homog I$ is the ideal generated by all the homogeneous
elements in~$I$. The following easy proposition contains the essential
observation that an image of an irreducible variety is irreducible.
We use this later to prove that certain multi-valued maps descend
to honest regular maps between toric varieties, even though on the
Cox rings the pathologies of Example~\ref{example:homog_prime_not_prime}
can occur.

\begin{prop}\label{prop:homog_of_kernel_is_homog_prime}
Let $\crl$ be a graded domain. If $\gotp \ideal \crl $ is a prime ideal,
then $\homog \gotp$ is homogeneously prime. In particular, if $R$
any domain and $\alpha\colon \crl \to R$ is any ring homomorphism,
then $\homog{(\ker \alpha)}$ is homogeneously prime.
\end{prop}

If $R$ is a ring and $I \subset R$ is any subset,
then we use $\langle I \rangle$ or $\langle I \rangle_R$
to denote the ideal generated by set $I$.
We use this notation very often in the case that
$S \subset R$ is a subring and $I \ideal S$ is an ideal.
Then $\langle I \rangle_R \ideal R$ is the extension of the
ideal $I$ in the ring $R$.

\begin{defin}{\cite[I.3]{hartshorne}}
Let $\crl$ be a graded ring and let $\gotp \ideal \crl$ be a homogeneously
prime ideal. Then the set $A$ of all homogeneous elements in $\crl$
which are not in $\gotp$ is multiplicative, and the
\textbf{(homogeneous) localisation \hlocalise{\crl}{\gotp}}
is defined to be the set of degree $0$ elements in $\inv{A} \crl$.
It is a local ring with maximal ideal
$(\gotp \cdot A^{-1}\crl) \cap \hlocalise{\crl}{\gotp}$.

If $f\in \crl$ is homogeneous, define the
\textbf{(homogeneous) localisation \hlocalise{\crl}{f}}
to be the set of degree $0$ elements in $\crl[f^{-1}]$.
If $I\ideal \crl$ is a homogeneous ideal,
then $\hlocalise{I}{f}$ is the set of degree $0$ elements in
$\langle I \rangle_{\crl[f^{-1}]}$; equivalently,
\[
    \hlocalise{I}{f}= \langle I \rangle_{\crl[f^{-1}]} \cap \hlocalise{\crl}{f}.
\]
\end{defin}

When $\crl=\cox X$ is the Cox ring of a toric variety $X$ and $Z \subset X$
an irreducible subvariety defined by a homogeneously prime ideal $I(Z) \ideal \cox X$, the localisation
$\hlocalise{\cox X}{I(Z)}$ is equal to the local ring of point $Z$
in the scheme $X$:
\[
  \hlocalise{\cox X}{I(Z)} = \set{q \in \functionfield{X} \mid Z\cap \Reg q \ne \emptyset}.
\]
This is analogous to the usual statement for $\Proj$ of an $\NN$-graded
ring: see \cite[Prop.~II.2.5(a)]{hartshorne}, for example.
Localisation at an element $f$ is also analogous to the case of usual
$\Proj$. Roughly, $S_{(f)}$ consists of global rational functions that
are regular an open subset $X_f= X \setminus (f)$, but there are caveats.
First, if $X$ has nontrivial $\CC^*$-factors, then we assume that the
zero locus of $f$ contains the resulting divisorial components of the
irrelevant locus $\irrel{X}$.
Second, the open subset $X_f$ is not necessarily affine, so regular
functions on $X_f$ might be scarce (or even all constant).




\begin{defin}
An ideal $I \ideal \cox X$ is \textbf{relevant} if it does not
contain any power of the irrelevant ideal $B_X$.
\end{defin}

Note that if $I$ is relevant, then $\homog I$ is relevant too.

\begin{lemma} \label{lem:R_I_generates_KX}
Let $X$ be a toric variety and $\gotp \ideal \cox X$ a homogeneously
prime ideal. Set $R=\hlocalise{\cox X}{\gotp}$. If $\gotp$ is
relevant, then $R$ and $R^{-1}$ generate $\functionfield{X}$.
\end{lemma}

\begin{prf}
Let $A$ be the set of all homogeneous elements in $\cox X$ which are not
in $\gotp$, so $R=(A^{-1} \cox X)^0$. We consider the subset $\mu\subset A$
of monomials not in $\gotp$; 
we will find enough elements to generate $\functionfield{X}$ from that.
We treat $\mu$ naturally as a subset
$\mu\subset\monoposcone{X}$.
Here and below $\monoposcone{X}$ and $\coxmonolatt{X}$ are as defined in \S\ref{sect:Cox_coordinates}.
In fact, since $\gotp$ is homogeneously prime, $\mu$ forms a lattice cone
in $\coxmonolatt{X}$ which is a face of the positive cone $\monoposcone{X}$.

Set $\mu^*$ to be the face of the positive cone of the ray lattice $R_X$
that is dual to $\mu$ (that is, the span of the basis elements $\rho_i$
for which the corresponding Cox variable $x_i$ is not in $\mu$).
Let $\dualcone{(\mu^*)}\subset\coxmonolatt{X}$ be the cone dual to $\mu^*$,
which is precisely
\[
\dualcone{(\mu^*)}  = \{ z - y \mid z\in\monoposcone{X}, y\in\mu \}
\]
so the localisation $A^{-1} \cox X $ contains all the
monomials in $\dualcone{(\mu^*)}$.
For example, if $\coxmonolatt{X} \simeq \ZZ^2$ and $\mu$ is generated by $(1,0)$,
  then $\mu^*$ is generated by $\rho_2$, and $\dualcone{(\mu^*)} = \langle (1,0),(-1,0), (0,1)\rangle$.
Restricting only to those monomials of degree~$0$ with respect to the
gradings is the same as taking the pullback via the principal divisor map
$M_X \hookrightarrow \coxmonolatt{X}$, so to prove the claim it is enough
to prove that this pullback is a cone of maximal dimension in $M_X$.

The pullback above is simply the dual of the image of $\mu^*$ in $N_X$
under the ray lattice map. Since $\gotp$ is relevant, this image cone
is one of the cones in the fan, so it is strictly convex and therefore
its dual is of maximal dimension, as required.
\end{prf}



\subsubsection{Equations defining subschemes}
\label{sect:equations_of_subschemes}

 Subschemes are defined by ideals in Cox rings. We discuss different
choices here, which then arise later when considering
images and preimages of subschemes. 

\begin{defin}
  Let $X$ be a toric variety with Cox ring $\cox X$.
   If $I \ideal \cox X$ is a homogeneous ideal, then we write
  $R = \cox X / I$ for the graded quotient ring, and for
  $h\in \cox X$ we write $\tilde{h}$ for $h +I \in R$. 

  Suppose $A \subset X$ is a closed subscheme.
  \begin{itemize}
    \item We say \textbf{$I$ defines $A$}
          if for every affine open subset $X_h=X \setminus (h)$
          for some homogeneous $h$ in $\cox X$ we have equality of schemes:
            $A \cap X_h = \Spec \hlocalise{R}{\tilde{h}}$.
    \item We say \textbf{$I$ maximally defines $A$}
          if $I$ defines $A$ and $I' \subset I$ for any other
          $I'\ideal \cox X$ which defines~$A$.
    \item We say \textbf{$I$ freely defines $A$}
          if $I$ defines $A$ and $I$ is generated by $\fromto{f_1}{f_k}$
          for some homogeneous $f_i \in \cox X$, such that $f_i$ defines a Cartier divisor.
  \end{itemize}
\end{defin}


\begin{lemma}\label{lem:defining_ideals_are_additive}
Suppose $A_1, A_2 \subset X$ are two closed subschemes defined by homogeneous ideals $I_{A_1}, I_{A_2}$,
respectively.
Then the scheme theoretic intersection $A_1 \cap A_2$  is defined by $I_{A_1} +I_{A_2}$.
\end{lemma}

\begin{example}
   Let $X=\PP(1,1,2)$ with Cox coordinates $ x_1, x_2, x_3$
   and let $A$ be the coordinate locus $x_2=0$.
   Then the ideal $I_{\mathrm{max}} = \langle x_2 \rangle$ maximally defines $A$
   and $I_{\mathrm{free}} = \langle x_1 x_2, {x_2}^2  \rangle$ freely defines~$A$.
\end{example}

In practice, ideals maximally defining a subscheme are often the simpler ones
and describe global properties of the scheme, while
ideals freely defining a subscheme say more about local properties.
For instance in the example above we immediately see that $A$ is not a local complete intersection.

Following Kajiwara \cite{kajiwara}~1.5, we say that a toric
variety $X$ {\em has enough Cartier divisors} if the
complement of each torus invariant affine patch on $X$
supports an effective $T$-invariant Cartier divisor.
We also say that an ideal $I \ideal \cox{X}$ is \textbf{saturated}
if $(I:B_X) = I$, or equivalently if the scheme in $\CC^m$ defined
by $I$ has no (embedded) components with support
on $\irrel X \subset \CC^m$.

\begin{prop}
   Let $X$ be a toric variety and $A \subset X$ a closed subscheme. Then there
   exists a unique homogeneous ideal $I_{\mathrm{max}} \ideal \cox X$
   maximally defining $A$, and this ideal is saturated.
   If, furthermore, $X$ has enough Cartier divisors (in the sense above),
   then there exists a saturated ideal $I_{\mathrm{free}}$ freely defining~$A$.
\end{prop}
Recall that if $X$ is $\QQ$-factorial or $X$ is quasiprojective, then it has
enough Cartier divisors. Our methods do not require this condition, except
where stated.

The saturatedness property is essential for accurate calculations of image
of a subvariety.
\begin{example}\label{example:saturation}
  Let $X = \PP^1 \times \CC$, $Y =\CC$ and let $\varphi\colon X \to Y$
  be the projection
  described in coordinates as $\Phi(x_1,x_2,x_3) = (x_3)$.
  Then $\cox{X}= \CC[x_1,x_2,x_3]$ with $B_X = \langle x_1, x_2 \rangle$.
  Let $I_{1}= \langle x_1 x_3 , x_2  x_3 \rangle$ and $I_{2} = \langle x_3 \rangle$.
  Then $I_{1}$ is not saturated, its saturation is $I_{2}$
  and the scheme theoretic image $\overline{\varphi}(A)$ of the scheme $A \subset X$
  given by either of these ideals is equal to the scheme given by
  $\langle y_1 \rangle \ideal \cox{Y}$.
  This ideal is obtained as $\inv{(\varphi^*)} (I_{2})$ whereas
  $\inv{(\varphi^*)} (I_{1}) = \langle 0 \rangle$.
\end{example}

\subsubsection{Rational maps}\label{sect:rational_maps}

We assemble standard facts about the image of subschemes under rational maps.
Let $X$ and $Y$ be two (irreducible) algebraic varieties with fields of rational
functions $\functionfield{X}$ and $\functionfield{Y}$.
Suppose $A \subset X$ is a closed subscheme; we denote the
corresponding ideal sheaf by $\ccI_{A} \ideal \ccO_X$.

Given a rational map $ \varphi \colon X \dashrightarrow Y$,
we denote by $\Reg \varphi \subset X$ the maximal open subset on which
$\varphi$ is regular and by $\phireg$ the restricted (regular) map
$\varphi|_{\Reg \varphi}$.
Suppose $U \subset \Reg\varphi$ is an open subset.
By definition, the scheme-theoretic image $\phibar|_U(A) \subset Y$ of
$A$ under $\varphi$ restricted to $U$ is
the minimal closed subscheme of $Y$
such that $\phireg|_{A \cap U}$ factorises through $\phibar|_U(A)$.
Set theoretically, $\phibar|_U(A)$ is supported on $\overline{\phireg(A \cap U)}$.
We write $\phibar(A)$ for $\phibar|_{\Reg\varphi}(A)$.

For a closed irreducible subvariety $Z \subset Y$ let
$\ccO_{Y,Z} \subset \functionfield{Y}$ be the local ring
of $Z$ with maximal ideal $\gotm_{Y,Z}$.
The next proposition is standard; see \cite[\S{}I.4]{hartshorne}
or \cite[\S{}V.1.1]{eisenbud_harris}, for example.

\begin{prop}\label{prop:rational_maps}
   Let $\varphi\colon X \dashrightarrow Y$ be a rational map between
  algebraic varieties. 
  \begin{enumerate}
   \item  If $Z=\phibar(X)$, then
	 $Z$ is reduced and irreducible and pullback
	 determines a ring homomorphism $\varphi^* \colon
	 \ccO_{Y,Z} \to \functionfield{X}$ with kernel $\gotm_{Y,Z}$. 
   \item \label{item:ring_homomorphism_determines_map}
         Conversely, suppose $R \subset \functionfield{Y}$ is a subring such
         that $R$ and $R^{-1}$ generate $\functionfield{Y}$ (as a ring). Then
         every ring homomorphism $\alpha \colon R \to \functionfield{X}$
         uniquely determines a rational map $\psi\colon X \dashrightarrow Y$
         such that $\psi^*|_R = \alpha$ and $R \subset \ccO_{Y,Z}$,
         where $Z = \overline{\psi}(X)$.
   \item \label{item:image_under_rational_map}
         If $A \subset X$ is a closed subscheme and
         $V\subset Y$ is an open affine subset, then
         \[
           \ccI_{\phibar(A)}(V) =
           \inv{(\varphi^*)} \ccI_{A} (\phireginv V) \ideal \ccO_Y(V).
         \]
   \item \label{item:preimage_under_rational_map}
   If $B \subset Y$ is a closed subscheme
   and $U\subset \Reg \varphi$ is an open affine subset, then
        \[
          \ccI_{\phireginv(B)}(U) =
          \left\langle \varphi^* \ccI_B \right\rangle
          \ideal \ccO_{\Reg \varphi}(U)
        \]
        determines the ideal sheaf of the preimage of $B$, also denoted
        $\ccI_B \cdot \ccO_{\Reg \varphi}$ in this context.
    \end{enumerate}
\end{prop}


 Analogous algorithms
compute the image of a point and the preimage of a subscheme under
a map between toric varieties expressed in Cox coordinates; see
\S\ref{sect:image_of_point} and \S\ref{sect:preimage}. 

The next proposition describes the locus where a rational map is
regular; it is used later to prove the existence of `complete' descriptions.

\begin{prop}\label{prop:algebraic_description_of_regularity_locus}
  Let $\varphi \colon X \dashrightarrow Y$ be a rational map of
  irreducible varieties. Let $\set{V_i}$ be an affine cover of $Y$ and $I$
  be the set of those $i$ for which $V_i\cap\varphi(X)$ is nonempty.
  Let $G_i$ be a set of generators of the affine coordinate ring $\ccO_{V_i}$.
  Then the locus where $\varphi$ is regular is
  \[
   \Reg \varphi = \bigcup_{i \in I} \bigcap_{g \in G_i} \Reg \varphi^*g.
  \]
\end{prop}

\begin{prf}
  It is enough to assume that $Y = V_1$ is affine and then, by composing
  it with closed immersion into an affine space, that $Y$ is an affine space
  and $G_1$ is the set of coordinate functions. In that case the statement
  is clear.
\end{prf}

\subsection{Field extensions}
\label{sect:field-extensions}

Throughout this subsection we assume $\FF$ is a field which contains
all the roots of unity.
We denote the algebraic closure of $\FF$ by $\overline{\FF}$.
Our main interest is $\FF=\CC( \fromto{x_1}{x_m})$
or a finite extension of this.

\begin{lemma} \label{lemma:power_defines_extension}
Let $\mvsa \in \overline{\FF}$ be such that $\mvsa^r \in \FF$ for $r>0$
and assume $r$ is minimal with this property.
  Then the polynomial $t^{r} -{\mvsa}^{r} \in \FF[t]$ is the minimal polynomial of $\mvsa$.
  In particular, the extension $\FF \subset \FF(\mvsa)$ is of degree $r$.
\end{lemma}

\begin{prf}
  Let $\epsilon$ be a primitive $r$-th root of unity. Then in $\overline{\FF}[t]$ we have
  \[
    t^{r} -{\mvsa}^{r} = (t - {\mvsa}) (t - \epsilon {\mvsa})\dotsm (t - {\epsilon}^{r-1} {\mvsa}).
  \]
  If $p \in \FF[t]$ is the minimal polynomial of ${\mvsa}$,
  then $p$ divides $t^{r} -\mvsa^{r}$ (see \cite[\S V.1]{lang}).
  Hence (up to a scalar in $\FF$) $p$ must be a product of a subset of $j$
  of the factors of $t^r-\mvsa^r$ above for some $0<j\le r$.
  But then $p(0) = {\epsilon}^{N} {\mvsa}^j$ for some power $N$.
  Hence ${\mvsa}^j \in \FF$, and so by minimality of $r$ we
  must have $j=r$ and $p= t^{r} -{\mvsa}^{r}$ as claimed.
  The degree calculation follows by \cite[Prop.~V.1.4]{lang}.
\end{prf}

\begin{cor}
  Consider a sequence of field extensions
  \[
    \FF=\FF_0 \subset \FF_1 \subset \dotsb \subset \FF_a = \FF(\mvsa_1, \dotsc, \mvsa_a)
  \]
  where $\FF_i = \FF_{i-1}(\mvsa_i)$ and each $\mvsa_i$ to some power is in $\FF$.
  Set $r_i$ to be the minimal positive integer such that ${\mvsa_i}^{r_i} \in \FF_{i-1}$.
  Then the collection
  \[
    \set{{\mvsa_1}^{j_1} \dotsm {\mvsa_a}^{j_a} \mid j_i \in \set{0,\dotsc, r_i -1}, \  i \in  \set{1,\dotsc, a} }
  \]
  forms a basis of $\FF(\mvsa_1, \dotsc, \mvsa_a)$ as a $\FF$-vector space.
\end{cor}

\begin{prf}
  Follows immediately from Lemma~\ref{lemma:power_defines_extension} and \cite[Prop.~V.1.2]{lang}.
\end{prf}

The following lemma is elementary, however we did not find any reference for this fact.
\begin{lemma}\label{lemma:irrational_part_is_the_same1}
  Assume $\mvsa_0, \dotsc, \mvsa_a \in \overline{\FF}$
  are all such that ${\mvsa_i}^{r_i} \in \FF$  for some $r_i > 0$  and
  $\mvsa_0+ \dotsb + \mvsa_a =0 $.
  Then the set $\Xi=\set{\mvsa_0, \dotsc, \mvsa_a}$ divides into a union
  of disjoint subsets
  $
    \Xi = \Xi_1 \sqcup \dotsb \sqcup \Xi_b
  $
  such that for each $j$  all $\mvsa \in \Xi_j$ are proportional over $\FF$ and
  $
    \sum_{\mvsa \in \Xi_j} \mvsa =0.
  $
\end{lemma}

\begin{prf}
  We argue by induction on $a$.
  If $a = 0$, then there is nothing to prove,
  so assume the result holds for all values less than $a\ge 1$ and
 that $\mvsa_i\not=0$ for every~$i$.

  Let $r_i$ be the minimal positive integers for which ${\mvsa_i}^{r_i} \in \FF$
  and let $\epsilon_i$ be a primitive $r_i$-th root of unity.
  Without loss of generality we may assume that $r_0$ is maximal among the $r_i$.
  By Lemma~\ref{lemma:power_defines_extension},
  $t^{r_i}- \mvsa_{i}^{r_i} \in \FF[t]$ is the minimal polynomial of $\mvsa_{i}$.

  Consider $\mvsa_0 = -(\mvsa_1+ \dotsb + \mvsa_a)$.
  The polynomial
  \[
    q(t)=\prod_{
     \substack{
       j_1 \in \set{0, \dotsc r_1-1} \\
       \vdots\\
       j_a \in \set{0, \dotsc r_a-1}}}
    (t + {\epsilon_1}^{j_1} \mvsa_1 + \dotsb + {\epsilon_a}^{j_a} \mvsa_a)
  \]
  is in $\FF[t]$ and it vanishes at $\mvsa_0$.
  Hence the irreducible polynomial $t^{r_0} -{\mvsa_0}^{r_0}$
  must divide $q(t)$.
  In particular
  \[
   {\epsilon_0} (\mvsa_1 +\dotsc+ \mvsa_a) =
   {\epsilon_1}^{j_1} \mvsa_1 + \dotsb + {\epsilon_a}^{j_a} \mvsa_a
  \quad\text{for some }j_1, \dotsc, j_a.
  \]
  Writing $\mvsb_i = (\epsilon_0 - {\epsilon_i}^{j_i}) \mvsa_i$ for $i \in \set{1,\dotsc a}$,
  this equation becomes $\mvsb_1 + \dotsb + \mvsb_a =0$.
  Now, by the inductive assumption, the $\mvsb_i$ divide into groups
  $\Xi_1, \dotsc, \Xi_b$, each of whose elements are proportional over $\FF$,
  and for which
  \begin{equation}\label{eq:sum_Zeta_k_is_0}
    \sum_{\mvsb \in \Xi_k} \mvsb =0
  \quad\text{for each }k.
  \end{equation}

  We consider three cases.
  In the first case, suppose there exist two different numbers $i_1$ and $i_2$,
  such that $\mvsb_{i_1}$ and $\mvsb_{i_2}$ belong to the same set $\Xi_k$ and
  such that $\epsilon_0 - {\epsilon_{i_1}}^{j_{i_1}}  \ne 0$ and $\epsilon_0 - {\epsilon_{i_2}}^{j_{i_2}}  \ne 0$.
  Without loss of generality, assume $i_1 =a- 1$  and $i_2=a$.
  Then
  \begin{align*}
               \mvsa_{a-1} & = \frac{\mvsb_{a-1}}{\epsilon_0 - {\epsilon_{a-1}}^{j_{a-1}}} = g_{a-1} \mvsb \\
  \text{ and }  \ \mvsa_{a} & = \frac{\mvsb_{a}}{\epsilon_0 - {\epsilon_{a}}^{j_{a}}} =  g_{a} \mvsb
  \end{align*}
  for some $g_{a-1}, g_{a} \in \FF$.
  So the $a$-tuple of elements
  $
    \mvsa_0, \mvsa_1, \dotsc, \mvsa_{a-2}, (\mvsa_{a-1} + \mvsa_{a})
  $
  satisfies the conditions of the lemma and we use our inductive assumption to conclude.
  Note that $\mvsa_{a-1}$ and $\mvsa_{a}$ either form a new group on their own
  (if $\mvsa_{a-1} = - \mvsa_{a}$)
  or they both must be proportional to the elements
  of one of the groups existing by the inductive assumption.

  In the second case, suppose that within a given group $\Xi_k$ there is only one such $i$
  that $\epsilon_0 - {\epsilon_{i}}^{j_{i}}  \ne 0$.
  Then from Equation \eqref{eq:sum_Zeta_k_is_0} we deduce that $\mvsa_i =0$,
  contrary to our assumption.

  Finally, as the third case, suppose that $\epsilon_0 = {\epsilon_{i}}^{j_{i}}$ for all $i$.
  In particular, $r_0$ divides $r_i$. But since we assumed $r_0$ to be maximal among $r_i$,
  we have
  $
    r_0 = r_i$
  for all $i$.

  We have shown
  that for every $(a+1)$-tuple satisfying the hypotheses of the lemma
  either all elements of the tuple are divided into groups proportional over $\FF$
  or their minimal powers are all equal.

  To conlude, consider the following $(a+1)$-tuple which also satisfies the hypotheses:
  \[
    1, \frac{\mvsa_1}{\mvsa_0}, \dotsc , \frac{\mvsa_a}{\mvsa_0}.
  \]
  If it divides into the appropriate  groups proportional over $\FF$,
  then so do the $\mvsa_i$.
  On the other hand if the minimal powers are equal, then all they are equal to $1$
  (because $1^1 \in \FF$).
  In that case, $\frac{\mvsa_i}{\mvsa_0} \in \FF$ and so again the $\mvsa_i$
  are proportional over $\FF$ (actually forming just one group $\Xi_1$ in this case).
\end{prf}

\begin{cor}\label{cor:irrational_part_is_the_same}
 Let $k\subset \FF$ be a subfield, and $V$ be a $k$-vector space with a $k$-linear map
 $i\colon V \to \overline{\FF}$
 such that  for every $\mvsb \in i(V)$ there is some $r>0$ for which
 $\mvsb^r \in \FF$.
 Then there exists $\mvsa \in \overline{\FF}$ and a $k$-linear map $j\colon V \to \FF$, such that
 $
 i(v) = j(v) \cdot \mvsa$
 for all $v\in V$.
\end{cor}
We express the conclusion of this corollary by saying that
the elements of $i(V)$  have a common irrational part.

\begin{prf}
 If $\dim i(V) =0$, then there is nothing to prove, so assume $\dim i(V) \ge 1$.
 Fix a non-zero element $\mvsa \in i(V)$ and take any other element $\mvsb \in i(v)$.
 Apply Lemma~\ref{lemma:irrational_part_is_the_same1} for the triple $\mvsa, \mvsb, -(\mvsa+\mvsb)$
 to conclude that
 $
   \mvsb = h(\mvsb) \cdot \mvsa$
 for some $h(\mvsb) \in \FF$.
 The implicit map $h$ is clearly $k$-linear in $\mvsb$, so define a $k$-linear map $j$ by
 $
   j(v)=h(i(v)).
 $
 These $j$ and $\mvsa$ have the required properties.
\end{prf}

\subsection{Simple ring extensions}\label{sect:ring_extensions}

Let $k$ be a field which contains all roots of unity and
$\crl$ an integral $k$-domain with field of fractions $\FF$.
Let $\overline{\FF}$ be the algebraic closure of $\FF$.
\begin{defin}\label{def:simple_extension}
  A ring $\mrl$ is called a \textbf{simple extension of $\crl$}
  if there exist $\mvsa_1, \dotsc, \mvsa_a\in\overline{\FF}$, with each
  ${\mvsa_i}^{r_i} \in \crl$ for some $r_i >0$ (which is assumed to be minimal),
  for which
  \begin{enumerate}
    \item\label{item:define-simple-extension-generators}
    $\mrl= \crl[\mvsa_1, \dotsc, \mvsa_a]$,
    \item\label{item:define-simple-extension-free}
    $\mrl$ is a free $\crl$-module with basis
          $\set{{\mvsa_1}^{l_1}  \dotsm {\mvsa_a}^{l_a}  \mid 0 \le l_i < r_i}$, and
    \item\label{item:define-simple-extension-normal}
    For any $\mvsb$ in the field of fractions $K(\mrl) \subset \overline{\FF}$
    of $\mrl$, if $\mvsb^r \in \crl$ for some integer $r > 0$,
    then $\mvsb \in \mrl$.
  \end{enumerate}
  The elements $\mvsa_1, \dotsc, \mvsa_a$ are
  called the \textbf{distinguished generators of $\mrl$ over $\crl$.}
\end{defin}

We establish some basic properties of simple ring extensions
as a corollary of~\S\ref{sect:field-extensions}.

\begin{cor}\label{cor:about_K(S)_and_zeta}
  Let $\crl$, $\FF$, $\mrl$ and the $\mvsa_i$ be as in
  Definition~\ref{def:simple_extension}\ref{item:define-simple-extension-generators} and
  \ref{item:define-simple-extension-free}.
  Let $K(\mrl) \subset \overline{\FF}$ be the field of fractions of $\mrl$.
  Let $\mvsb \in K(\mrl)$ be such, that $\mvsb^r \in \FF$ for some $r$.
  Then:
  \begin{enumerate}
    \item \label{item:basis_of_K(S)}
          $K(\mrl)$ is a vector space over $\FF$ with basis
          $\set{{\mvsa_1}^{l_1}  \dotsm {\mvsa_a}^{l_a}  \mid 0 \le l_i < r_i}$.
    \item \label{item:zeta_is_simple_to_write}
          $\mvsb = g\cdot {\mvsa_1}^{l_1}  \dotsm {\mvsa_a}^{l_a}$
          for some $g \in \FF$ and $0 \le l_i < r_i$.
    \item \label{item:when_is_zeta_regular}
          $\mvsb \in \mrl$ if and only if $g\in \crl$, where $\mvsb$ is
          expressed in the basis as in \ref{item:zeta_is_simple_to_write}.
          In particular, $\mrl \cap \FF = \crl$.
    \item \label{item:polynomial_for_xi_is_irreducible}
          Fix any $j \in \set{1,\dotsc, a}$.
          Let $\mrl_{j-1}$ be the ring $\crl[\mvsa_1, \dotsc, \mvsa_{j-1}]$
          and let $K(\mrl_{j-1})$ be its field of fractions.
          Then the polynomial $t^{r_j} - {\mvsa_j}^{r_j}$
          is irreducible in $K(\mrl_{j-1})[t]$.
  \end{enumerate}
\end{cor}

\begin{prf}
  To prove that \ref{item:basis_of_K(S)} holds, observe that $K(\mrl)$ is
  $\FF$-generated by the listed elements because $\mrl$ is.
  On the other hand if there were an $\FF$-linear relation between these generators,
  then after clearing the denominator there  would be a relation between these $\crl$-generators of $\mrl$,
  contradicting the assumption that $\mrl$ is the free module.

  \smallskip

  To prove \ref{item:zeta_is_simple_to_write} using \ref{item:basis_of_K(S)},
  write
  $
    \mvsb  =  \mvsb_1 + \dotsb + \mvsb_b
  $
  where each $\mvsb_i$ is of the form $g_i\cdot {\mvsa_1}^{l_{1,i}}  \dotsm {\mvsa_a}^{l_{a,i}}$.
  Setting $\mvsb_0 = - \mvsb$ we can apply Lemma~\ref{lemma:irrational_part_is_the_same1}
  to deduce that actually the $\mvsb_i$ divide into groups of elements
  proportional over $\FF$ such that the sum in each group is 0.
  In particular, $\mvsb_0$ must be either $0$ or proportional over $\FF$ to at
  least one of the $\mvsb_i$,
  which finishes the proof of \ref{item:zeta_is_simple_to_write}.
  Part~\ref{item:when_is_zeta_regular} follows immediately from
  \ref{item:zeta_is_simple_to_write}.

   \smallskip

   In \ref{item:polynomial_for_xi_is_irreducible}, $r_j$ is also the minimal positive integer,
   such that ${\mvsa_j}^{r_j} \in K(\mrl_{j-1})$,
   for otherwise, we would have an $\FF$-linear relation between smaller powers of the $\mvsa_i$
   contrary to \ref{item:basis_of_K(S)}.
   So the conclusion follows from Lemma~\ref{lemma:power_defines_extension}.
\end{prf}

Our main concern is a particular class of simple extensions.

\begin{prop}\label{propostion:simple_extension_of_poly_ring}
  Let $\crl=\CC[\fromto{x_1}{x_m}]$ be a polynomial ring.
  Let $\fromto{g_1}{g_a}$ be square free, pairwise coprime polynomials and
  $\fromto{r_1}{r_a}$ be positive integers. Set $\mvsa_i = \sqrt[r_i]{g_i}$.
  Then $\mrl= \crl[\mvsa_1, \dotsc, \mvsa_a]$ is a simple extension of $\crl$
  with distinguished generators $\mvsa_1, \dotsc, \mvsa_a$.
\end{prop}

\begin{prf}
 Since the polynomials are pairwise coprime,
 there is no polynomial relation between the $\mvsa_j$
 other than those generated by ${\mvsa_j}^r - g_j = 0$.
 Thus $\mrl$ is a free module over $\crl$ with the desired basis and $\mrl$
 satisfies \ref{item:define-simple-extension-generators} and \ref{item:define-simple-extension-free} of
 Definition~\ref{def:simple_extension}.

 Suppose $\mvsb\in K(\mrl)$ and $\mvsb^r \in \crl$ as in
 Definition~\ref{def:simple_extension}\ref{item:define-simple-extension-normal}.
 Then by Corollary~\ref{cor:about_K(S)_and_zeta}\ref{item:zeta_is_simple_to_write}
 we can write $\mvsb = \frac{g}{h_1 \dotsm h_b} \cdot {\mvsa_1}^{l_1}  \dotsm {\mvsa_a}^{l_a}$,
 where $g$  and the $h_j$ are non-constant polynomials in $\crl$, the $h_j$ are irreducible, and
 none of the $h_j$ divides $g$. We claim $b=0$ so that the denominator does not exist,
 as required by  Definition~\ref{def:simple_extension}\ref{item:define-simple-extension-normal}.
 Suppose on contrary, that $b \ge 1$. Since the $g_i$ are pairwise coprime,
 there is at most $1$ of  $\fromto{g_1}{g_a}$
 which is divisible by $h_1$ (say this is $g_i$)
 and since $g_i$ is square free, it can only divide $h_1$ with multiplicity $1$.
 Thus the multiplicity of $h_1$ in $\mvsb^r \in \crl$ is $-r+\frac{r l_i}{r_i}$,
 which is always negative, a contradiction.
\end{prf}

A simple verification of the definition confirms 
that simple extensions behave well under localisation.
\begin{prop}\label{prop:simple_extension_under_localisation}
  Suppose $\crl \subset \mrl$ is a simple ring extension
  and that $f \in \crl$.
  Then $\crl[\inv f] \subset \mrl[\inv f]$
  is a simple ring extension with the same set of distinguished generators.
\end{prop}

For $\mvsb \in K(\mrl)$ with $\mvsb^r \in \FF$ write
$\mvsb = g \cdot {\mvsa_1}^{l_1} \dotsm {\mvsa_a}^{l_a}$ with $0 \le l_i < r_i$ and $g \in \FF$
as in Corollary~\ref{cor:about_K(S)_and_zeta}\ref{item:zeta_is_simple_to_write}.
Define \textbf{the floor  $\lfloor \mvsb \rfloor$} and
 \textbf{the ceiling $\lceil\mvsb \rceil$ of $\mvsb$} to be
\[
  \lfloor \mvsb \rfloor := g
  \qquad\textrm{and}\qquad
  \lceil\mvsb \rceil:=
  	g \cdot {\mvsa_1}^{\epsilon_1 r_1} \dotsm {\mvsa_a}^{\epsilon_a r_a}
\]
respectively, where $\epsilon_i = \lceil l_i/r_i \rceil$ is either $0$
(if $l_i =0 $) or $1$ (if $l_i >0$).
They are both elements of $\FF$, and are related by
$
\left\lfloor\dfrac{1}{\mvsb}\right\rfloor = \dfrac{1}{\lceil\mvsb\rceil}.
$
\begin{prop}\label{prop:floor-and-ceiling}
  For $\mvsb$ as above, the floor and ceiling of $\mvsb$ satisfy both
  \[
   \mvsb \in \mrl \iff \lfloor \mvsb \rfloor \in \crl
   \qquad\textrm{and}\qquad
   \mvsb \in \mrl \Longrightarrow \lceil\mvsb \rceil \in \crl.
   \]
   Moreover, if $\mvsb$ is an invertible element of $\mrl$, then $\lfloor \mvsb \rfloor$
   and $\lceil\mvsb \rceil$  are invertible elements of $\crl$.
\end{prop}
\noprf

So far we did not exploit the property~\ref{item:define-simple-extension-normal} of
Definition~\ref{def:simple_extension}.
It is a normality condition, and it has two important consequences.
First, if $\mvsb \in K(\mrl)$ satisfies $\mvsb^r \in \FF$ for some $r>0$,
 then $\mvsb$ is regular on $\Spec \crl$
 (as a multi-valued function, meaning that it has no poles;
  see Definition~\ref{def:multivalued_section_is_regular_etc})
 if and only if $\mvsb \in \mrl$.
This is made precise in the proof of Lemma~\ref{lemma:floor_is_invertible}.
Meanwhile we illustrate it with an example.
\begin{example}
  Suppose $\crl:= \CC[x_1,x_2]$, and let  $\mrl' = \crl[\mvsa]$
  where $\mvsa:=\sqrt[4]{x_1 {x_2}^2}$.
  Then the extension $\crl \subset \mrl'$ satisfies conditions
  \ref{item:define-simple-extension-generators}--\ref{item:define-simple-extension-free}
  of Definition~\ref{def:simple_extension}, but does not satisfy \ref{item:define-simple-extension-normal}: for example, the
  multi-valued function $\mvsb:= \sqrt{x_1} \in K(\mrl')$ has no poles,
  but $\mvsb = \mvsa^2 / x_2 \notin \mrl'$.
  Instead, we may consider a slightly bigger ring $\mrl = \crl[\sqrt[4]{x_1}, \sqrt[2]{x_2}]$.
  Then $\crl \subset \mrl$ is a simple extension and $\mvsa, \mvsb \in \mrl$.
\end{example}
The second consequence of \ref{def:simple_extension}\ref{item:define-simple-extension-normal}
is the uniqueness of $\lfloor \cdot \rfloor$ and $\lceil \cdot \rceil$ operations.
\begin{prop}\label{prop:uniqueness_of_floor}
  Suppose $\mvsb \in \overline{\FF}$ is such that $\mvsb^r \in \FF$.
  Then up to an invertible element in $\crl$,
  $\lfloor \mvsb \rfloor$ and $\lceil \mvsb \rceil$ are well defined elements of $\FF$,
  independent on the choice of simple ring extension $\crl \subset \mrl$
  such that $K(\mrl)$ contains $\mvsb$.
\end{prop}
\begin{prf}
  It is enough to prove the statement for $\lfloor \mvsb \rfloor$.
  More precisely, suppose
  $\mrl:=\crl[\fromto{\mvsa_1}{\mvsa_a}]$ and $\mrl':=\crl[\fromto{\mvsa'_1}{\mvsa'_b}]$
  are two simple ring extensions of $\crl$ with $\mvsb \in \mrl, \mrl'$.
  Write $\mvsb = g \cdot {\mvsa_1}^{l_1}  \dotsm {\mvsa_a}^{l_a}
               = g' \cdot {\mvsa'_1}^{m_1}  \dotsm {\mvsa'_b}^{m_b}$.
  We have to prove $g'/g \in \crl$ (inverting the roles of $\mrl$ and $\mrl'$
  we also get $g/g' \in \crl$).

  Observe that ${\mvsb}/{g} = {\mvsa_1}^{l_1}  \dotsm {\mvsa_a}^{l_a} \in \mrl$,
  thus  $\left({\mvsb}/{g}\right)^r \in \crl$ for some $r$.
  By Definition~\ref{def:simple_extension}\ref{item:define-simple-extension-normal}
  also ${\mvsb}/{g} \in \mrl'$.
  Since
  $
    {\mvsb}/{g} = ({g'}/{g}) \cdot {\mvsa'_1}^{m_1}  \dotsm {\mvsa'_b}^{m_b}
  $
  by Corollary~\ref{cor:about_K(S)_and_zeta}\ref{item:when_is_zeta_regular}
  we have $g' / g \in S$ as claimed.
\end{prf}

The next corollary shows that the intersection with $\crl$ is readily calculated
for certain ideals in simple extensions $\mrl\supset \crl$.


\begin{cor}\label{cor:generators_of_elimination_ideal_are_simple}
  Let $I \ideal \mrl$ be an ideal generated by $\mvsb_1, \dotsc, \mvsb_{\beta}$
  where each $\mvsb_i$ satisfies $\mvsb_i^{r_i} \in \crl$ for some $r_i>0$.
  Then
  \[
    I \cap \crl = \bigl\langle  \lceil\mvsb_1\rceil, \dotsc, \lceil\mvsb_{\beta}\rceil \bigr\rangle  \ideal \crl.
  \]
  In particular intersecting ideals generated by such $\mvsb_i$ in $\mrl$
  with $\crl$ is additive:
  \[
    (I_1 + I_2 ) \cap \crl  = (I_1 \cap \crl) + (I_2 \cap \crl).
  \]
\end{cor}

\begin{prf}
  This is repeated application of Lemma~\ref{lemma:elimination_of_homogeneous_irrational_element} below,
keeping in mind Corollary~\ref{cor:about_K(S)_and_zeta}\ref{item:polynomial_for_xi_is_irreducible}.
\end{prf}

\begin{lemma}\label{lemma:elimination_of_homogeneous_irrational_element}
 Let $\crl$ be an integral domain. Consider an integral domain
 $\mrl = \crl[\mvsa] / \left\langle \mvsa^r - g\right\rangle $
 for some  $r\in \ZZ$, $r > 0 $ and $g \in \crl$ for which $\mrl$ is
 a free $\crl$-module with basis $1, \fromto{\mvsa}{\mvsa^{r-1}}$
 (in particular, $\mvsa^r - g $ is irreducible over $\crl$).
 Furthermore assume $I$ is an ideal in $\mrl$ generated as
 \[
   I = \left\langle f_1,\dotsc,f_{\alpha},
       f_{\alpha+1}\mvsa^{m_{\alpha+1}},\dotsc, f_{\beta}\mvsa^{m_{\beta}} \right\rangle
 \]
 where $f_i \in \crl$ and $0 < m_{i} < r$.
 Then
 \[
   I \cap \crl = \left\langle f_1,\dotsc,f_\alpha, f_{\alpha+1} g,\dotsc, f_{\beta} g  \right\rangle.
 \]
\end{lemma}

\begin{prf}
 Clearly the listed generators are in $I \cap \crl$.

 So  consider $h \in I$:
 \[
   h = \left( \sum_{i=1}^{\alpha} \sum_{j=0}^{r-1} h_{i, j} f_i\mvsa^{j} \right) +
       \left( \sum_{i=\alpha+1}^{\beta} \sum_{j=0}^{r-1} h_{i, j} f_i\mvsa^{j+m_{i}} \right)
 \]
 for some $h_{i,j}$ in $\crl$.
 Rewrite $h$ as:
 \[
   h= \left( \sum_{i=1}^{\alpha} h_{i, 0} f_i \right) +
       \left( \sum_{i=\alpha+1}^{\beta} h_{i, r-m_i} f_i g \right)  +
        \mvsa \left( \dotsc \right) + \dotsb + \mvsa^{r-1} \left( \dotsc \right).
 \]
 If $h \in \crl$, then the summands with $\mvsa^i$ for $i \in \set{1,\dotsc,r-1}$ are all $0$
 (because $\mrl$ is a free $\crl$ module with basis $1, \fromto{\mvsa}{\mvsa^{r-1}}$).
 Hence:
 \[
   h= \left( \sum_{i=1}^{\alpha} h_{i, 0} f_i \right) +
       \left( \sum_{i=\alpha+1}^{\beta} h_{i, r-m_i} f_i g \right)
 \]
 which is an element of
 $
    \left\langle f_1,\dotsc,f_\alpha, f_{\alpha+1} g,\dotsc, f_{\beta} g  \right\rangle \ideal \crl
 $
 as claimed.
\end{prf}

\begin{lemma}\label{lem:extending_to_simple_extensions}
  Let $\crl$, $\FF$, $\mrl$ and the $\mvsa_i$ be as in
  Definition~\ref{def:simple_extension}.
  Analogously, let $\mrl'$ be a simple extension of an integral $k$-domain $\crl'$
  and let $\FF'$ be the field of fractions of $\crl'$.
  Assume $\Phi^* \colon \crl \to \mrl'$ is a homomorphism.
  Then $\Phi^*$ can be extended (non-uniquely) to a homomorphism
  $\widetilde{\Phi^*}\colon  \mrl \to \overline{\FF}$ as in the diagram:
  \[
    \xymatrix{
    \crl \ar[rrd]^{\Phi^*} \ar@{^{(}->}[d] && \crl' \ar@{^{(}->}[d]\\
    \mrl \ar@{.>}[rrd]^{\widetilde{\Phi^*}} \ar@{^(->}[d] && \mrl' \ar@{^{(}->}[d]\\
    \overline{\FF}               && \overline{\FF'}}
  \]
  (so that, in particular, the diagonal square is commutative).
  The extension can be chosen as follows.
  For every $i$, suppose ${\mvsa_i}^{r_i} \in \crl$ is the
  (minimal) defining property of $\mvsa_i$ and set $g_i:= {\mvsa_i}^{r_i}$.
  Then set
  \[
    \widetilde{\Phi^*}(\mvsa_i) := \sqrt[r_i]{\Phi^*(g_i)} \in \overline{\FF'}
  \]
  for any choice of the $r_i$th root.
  \end{lemma}

\begin{prf}
  Since the only polynomial relations between
  $\fromto{\mvsa_1}{\mvsa_a}$ are $g_i - {\mvsa_i}^{r_i}$,
  $\widetilde{\Phi^*}$ really defines a homomorphism.
\end{prf}

%% file: toricmapsmultivalued.tex

\section{Roots and multi-valued maps}
\label{sect:multival}

In this section we introduce the main technical tool to study descriptions of
maps between toric varieties. We extend the field of rational functions to
include special elements of its algebraic closure, so-called multi-valued
functions.
We use these multi-valued functions to define multi-valued maps in the same way
rational functions are used to define rational maps.

We fix notation for this section, and indeed for the rest of this paper.
We work with two toric varieties  $X$ and $Y$ and their Cox covers:
\[
 \xymatrix{
    \CC^m  = \Spec \cox X \ar@{-->}[d]_{\pi_X} &&  \Spec \cox Y  =  \CC^n \ar@{-->}[d]^{\pi_Y}\\
     X                               && Y
 }
\]
where $\cox X = \CC[x_1,\dotsc, x_m]$ and $\cox Y = \CC[\fromto{y_1}{y_n}]$.
Although in this section we work exclusively on the Cox covers $\CC^m$ and
$\CC^n$, and everything could be described with no reference to $X$ and $Y$,
we maintain the connection between the
Cox covers and their toric varieties in the notation.

\subsection{Multi-valued sections}
\label{sect:multi_valued_sections}

\begin{defin}
\label{def:multivalsec}
 A \textbf{multi-valued section on $X$} is an element $\mvsa$ of the algebraic closure $\overline{\coxfield X }$.
 We say $\mvsa$ is \textbf{homogeneous} if $\mvsa^r = f$ for some homogeneous $f \in \coxfield{X}$
  and for some integer $r \ge 1$.
\end{defin}

\begin{notation}
\label{notn:mvs}
If $\mvsa$ is a homogeneous multi-valued section with $\mvsa^r = f$ as above,
then we write $\mvsa=\sqrt[r]{f}$.
It is implicit in this notation that $r$ is minimal
and that an $r$-th root of $f$ has been chosen once and for all, and
any other use of $\sqrt[r]{f}$ in the same calculation 
refer to the same element $\mvsa$.
\end{notation}

The product and quotient of two homogeneous multi-valued sections is again homogeneous, but
their sum is usually not:
 $\sqrt{x_1}+\sqrt{x_2}$ is not homogeneous
even if $x_1$ and $x_2$ have the same degree.
Furthermore, it is not true that every multi-valued section can be expressed as a sum of homogeneous ones.

In the first place, we treat multi-valued sections on $X$ as mildly generalised
rational functions on the affine Cox cover $\CC^m$. In particular, we simply
define when a homogeneous multi-valued section is regular or invertible
on an open subset of $\CC^m$ following the notions for rational functions.

\begin{defin}\label{def:multivalued_section_is_regular_etc}
Let $\mvsa = \sqrt[r]{f}$ be a homogeneous multi-valued section of $X$ with
$f \in \coxfield{X}$ homogeneous. Then $\mvsa$ is \textbf{regular} if
$f\in\cox X$. More generally, $\mvsa$ is \textbf{regular on $U$}, for a
Zariski open subset $U \subset \CC^m$, if $f$ is regular on $U$. If $\mvsa$
is regular on $U$ and does not vanish anywhere on $U$, we say $\mvsa$ is
\textbf{invertible on $U$}.

The \textbf{domain of $\mvsa$}, also called the \textbf{regular locus of
$\mvsa$} and denoted $\Reg\mvsa$, is defined to be the largest open subset
of $\CC^m$ on which $\mvsa$ is regular.
\end{defin}

If $V\subset X$ is a Zariski open subset of $X$ and $\mvsa$ a homogeneous
multi-valued section of $X$, then we say that $\mvsa$ is regular on $V$
if it is regular on the open subset $\pi_X^{-1}(V)\subset \CC^m$.

A typical homogeneous multi-valued section $\mvsa = \sqrt[r]{f}$ is not a
function in the usual sense. Nevertheless, for $\xi\in\Reg\mvsa$ we write
$\mvsa(\xi)$ to denote the finite set of values $a \in \CC$ for which
$a^r = f(\xi)$.

\begin{defin}
A homogeneous multi-valued section $\mvsa = \sqrt[r]{f}$  is \textbf{single valued} if
$r=1$, in which case $\mvsa = f \in \coxfield{X}$.
\end{defin}

This notion relies on the convention of \ref{notn:mvs} that $r$ is assumed
to be minimal. Thus, for example, $\sqrt[r]{1}$ is single valued, since
the minimal choice is $r=1$.
Since we are in characteristic~0
and our ground field contains all roots of unity, there is an equivalent
set-theoretic condition (Proposition~\ref{prop:single_valued_sections}); we omit the proof.

\begin{prop}\label{prop:single_valued_sections}
  A homogeneous multi-valued section $\mvsa \in \overline{\coxfield{X}}$ is single-valued
  if and only if $\mvsa(\xi)$ has exactly one element for a general $\xi \in \Reg \mvsa$.
\end{prop}

Finally we show that linear subspaces of homogeneous multi-valued sections
all have the same irrational part. This is one of the key points that makes
the theory work: if we imagined a map to projective space as being determined
by a basis of a vector space of sections corresponding to a `multi-valued
linear system', then this property would allow us to divide out by the
common irrational part to recover a map defined without radicals.

\begin{prop}\label{prop:image_has_constant_degree_and_irrational_part}
If $V$ is a $\CC$-vector space and $i \colon V \to \overline{\coxfield{X}}$
is a $\CC$-linear map whose image consists of only homogeneous multi-valued
sections, then there exists a homogeneous multi-valued section $\mvsa \in
\overline{\coxfield{X}}$ and a $\CC$-linear map $j \colon V \to
\coxfield{X}$ whose image  consists of homogeneous elements of a constant
degree, and $ i(v) = j(v) \cdot \mvsa$  for all  $v\in V.$
\end{prop}

\begin{prf}
  If $\dim i(V) =0$, then there is nothing to prove, so assume $\dim i(V) \ge 1$.
  Apply Corollary \ref{cor:irrational_part_is_the_same} for $k= \CC$ and $\FF=\coxfield{X}$
 and let $j'$ and $\mvsa'$ be the resulting map and section.
 Let $v_0 \in V$ be a vector such that $j'(v_0)$ is not zero and set $\mvsa:=\mvsa'\cdot j'(v_0)$ and
 $j(v):=j'(v)/ j'(v_0)$.

 By assumption $(j(v)\cdot \mvsa)^r$ is a homogeneous section in $\coxfield{X}$ for some $r$.
 Hence $\mvsa^r \in \coxfield{X}$ and $\mvsa^r$ is homogeneous (take $v = v_0$).
 But then also $j(v)^r$ is homogeneous,
 being a quotient of two homogeneous sections, and so $j(v)$ is homogeneous too.

 It remains to prove that $j(v_1)$ and $j(v_2)$ have the same degree for any $v_1, v_2 \in V$.
 Consider $j(v_1+v_2) =j(v_1)+ j(v_2)$.
 If $j(v_1)$ and $j(v_2)$ had different degrees,
 then the decomposition of $j(v_1+v_2)$ into homogeneous components would have two components,
 but $j(v_1+v_2)$ is also homogeneous, so there can be only one component.
\end{prf}

\subsection{Multi-valued maps}\label{sect:multivalued_maps}

We define multi-valued ``maps'' between affine spaces allowing roots in
their descriptions. We will not consider the largest possible class of maps
that one might define by multi-valued sections, but only a particular case.

\begin{defin}
 A \textbf{multi-valued map $\Phi$ from $\CC^m$ to $\CC^n$} is a
 $\CC$-algebra homomorphism
 $
 	\Phi^* \colon \CC[\CC^n] \to \overline{\functionfield{\CC^m}}
 $
 such that $\Phi^*y_i$ is a homogeneous multi-valued section for each $i=
 \fromto{1}{n}$. We say $\Phi$ is \textbf{regular} on $U \subset \CC^m$ if
 all $\Phi^*y_i$ are regular on $U$.
\end{defin}

\begin{notation}
  If $\Phi$ is a multi-valued map as above, then we write
 \begin{align*}
   \Phi\colon \CC^m & \multito \CC^n \\
           \xi   & \multimapsto \Bigl(\bigl(\Phi^*y_1\bigr)(\xi),
           \dotsc,\bigl(\Phi^*y_n\bigr)(\xi) \Bigr).
 \end{align*}
 Of course evaluating $\Phi$ at a point $\xi \in \CC^m$ is slightly
 delicate. Each component is the evaluation of a multi-valued function, so
 it is a set. However $\Phi(\xi)$ is not necessarily the product of these
 sets, since we must match the roots appearing in the multi-valued sections
 when they are the same, as in  \S\ref{sect:example_weighted_blow_up}. The
 evaluation will be explained in detail in \S\ref{sect:image_of_point}.
\end{notation}

We extend $\Phi^*$ to a subset of rational functions (for which the pullback makes sense, i.e. we do not divide by $0$) by
\[
  \Phi^*\left(\frac{f}{g}\right) = \frac{\Phi^*f}{\Phi^*g}
  \quad \text{ whenever } \Phi^*g \ne 0.
\]
If $q =f/g$ is a reduced expression and $\Phi^*g =0$, then we say $\Phi^* q$
is not defined.

\begin{example}
The toric map of \S\ref{sect:example_weighted_blow_up},
 an affine patch on the blow up of the affine quotient singularity $\half(1,1)$,
 lifts to a multi-valued map
 \begin{align*}
   \Phi\colon \CC^2  & \multito \CC^2 \\
          (s,t) & \multimapsto (\sqrt{s}, t\sqrt{s}).
 \end{align*}
\end{example}

\begin{defin}
 Let $\varphi$ be a rational map $\CC^m \dashrightarrow \CC^n$.
 We can naturally associate a multi-valued map $\Phi$ to $\varphi$,
  by assigning $\Phi^*:=\varphi^*$.
 If a multi-valued map $\Phi$ arises in this way, then we say  $\Phi$ is \textbf{single-valued}.
\end{defin}

The maximal subset $U \subset \CC^m$ on which $\Phi$ is regular is an open affine subset.

\begin{prop}\label{prop:regularity_locus_is_affine}
 Let
 \begin{align*}
   \Phi\colon \CC^m & \multito \CC^n \\
           x   & \multimapsto \left(\left(\frac{f_1}{g_1}\right)^{\frac{1}{r_1}},
                 \dotsc,\left(\frac{f_n}{g_n}\right)^{\frac{1}{r_n}}\right).
 \end{align*}
 be a multi-valued map.
 Assume that $f_i / g_i$ is reduced for each $i$.
 Then the maximal subset $U \subset \CC^m$ where $\Phi$ is regular is the
 complement of the vanishing locus of $g:=g_1\dotsm g_n$: that is,
  \[
   U= \bigl( \CC^m \bigr)_g = \Spec \cox X [g^{-1}].
  \]
  In particular, since the $g_i$ are homogeneous, $U$ is $G_X$-invariant,
  where $G_X$ is the group acting on $\CC^m$, making $X = \CC^m// G_X$,
     see the notation in \S\ref{sect:Cox_coordinates}.
\end{prop}

\begin{prf}
 Clearly $\Phi$ is regular on $\bigl( \CC^m \bigr)_g$.
 Further let $\xi$ be such that $g_i(\xi) =0$ for some $i$.
 Then $\Phi$ is not regular at $\xi$,
 because $\Phi^*y_i = \sqrt[r_i]{(f_i/g_i)}$ is not regular at $\xi$.
\end{prf}

\begin{defin}
\label{defin:regular_locus}
Let $\Phi\colon\CC^m \multito \CC^n$ be a multi-valued map. The
\textbf{domain of $\Phi$}, also called the \textbf{regular locus of $\Phi$}
and denoted $\Reg \Phi$, is defined to be the affine open subset $U$
of Proposition~\ref{prop:regularity_locus_is_affine}.
\end{defin}

\begin{cor}
If a description $\Phi$ is determined by polynomial radicals
\[
x\multimapsto (\sqrt[r_1]{f_1},\dots,\sqrt[r_n]{f_n}),
\]
for polynomials $f_1,\dots,f_n\in\cox X$,
then $\Reg\Phi = \CC^m$.
\end{cor}

\subsection{Map rings of multi-valued maps}
\label{sect:map_ring}

There is another natural way of thinking of a multi-valued map, and it
is the key to the analysis here.
Let $\Phi \colon \CC^m \multito \CC^n$  be a multi-valued map
with corresponding toric varieties $X$ and $Y$.
Choose $\mapring \Phi$ to be any subring in $\overline{\coxfield{X}}$
which satisfies the following properties:
\begin{enumerate}
  \item \label{item:condition_on_mapring_finite_generation}
        $\mapring \Phi = \CC[\Reg \Phi] [\mvsa_1, \dotsc, \mvsa_a]$
        for some homogeneous multi-valued sections $\mvsa_1,\dots,\mvsa_a$, all
        of which are regular on $\Reg\Phi$.
  \item \label{item:condition_on_mapring_containing_image}
        the image  $\Phi^* (\cox Y) = \Phi^*(\CC[\CC^n])$ is contained in
        $\mapring \Phi$.
  \item \label{item:condition_on_mapring_simple_extension}
        $\cox X [\mvsa_1, \dotsc, \mvsa_a]$ is a simple extension of $\cox X$
        with distinguished generators $\mvsa_1,\dots,\mvsa_a$
        (so by Proposition~\ref{prop:simple_extension_under_localisation}
        also $\mapring \Phi$ is a simple extension of $\CC[\Reg \Phi]$ with
        the same generators).
\end{enumerate}

Although such rings are not uniquely determined, they are important
in our considerations.
\begin{defin}
  Any ring $\mapring \Phi$ satisfying
  \ref{item:condition_on_mapring_finite_generation}, \ref{item:condition_on_mapring_containing_image}
  and \ref{item:condition_on_mapring_simple_extension}
  is called \textbf{a map ring of $\Phi$}.
\end{defin}

\begin{prop}
  Let $\Phi\colon \CC^m \multito \CC^n$ be a multi-valued map.
  Then there exists a map ring $\mapring \Phi$ of~$\Phi$.
\end{prop}

This proof is constructive but does not necessarily give the most
efficient way of choosing a map ring.

\begin{prf}
  Let $\Phi^*y_i = \sqrt[r_i]{f_i}$, where $ f_i \in \coxfield{X}$. Let
  $\set{g_1,\dots,g_a}$ be a finite set of homogeneous, square free, and pairwise
  coprime, polynomials in $\cox X$ so that each $f_i$ has an expression as a
  Laurent monomial in the $g_j$. Set $r$ to be lowest
  common multiple of all $r_i$.
  Then set
  \[
    \mvsa_j=\sqrt[r]{g_j} \quad \text{for all } j \in \setfromto 1 a.
  \]
  We claim that $\mapring \Phi= \CC[\Reg \Phi][\fromto{\mvsa_1}{\mvsa_a}]$
  is a map ring of $\Phi$.

  Property~\ref{item:condition_on_mapring_finite_generation} is satisfied by
  construction. It is also clear that each $\Phi^*y_i$ can be expressed in
  terms of the $\mvsa_j$, so $\mapring \Phi$ contains the image of $\cox{Y}$
  which is Property~\ref{item:condition_on_mapring_containing_image}. Finally
  Property~\ref{item:condition_on_mapring_simple_extension} follows by
  Proposition~\ref{propostion:simple_extension_of_poly_ring} since the $g_j$
  are coprime.
\end{prf}

The fact that the map ring is a simple extension has three advantages.
First, the image of a point can be calculated by a simple evaluation. Second, it
allows us to compose appropriate pairs of multi-valued maps. Finally, with
the distinguished generators (which only need be calculated once for each
map), preimages of subvarieties can be calculated, at least away from
certain loci.
Section~\ref{sect:properties_of_descriptions} explains these.

\subsection{Images and preimages under multi-valued maps}

Let $\Phi\colon \CC^m \multito \CC^n$ be a multi-valued map and $\mapring\Phi$
be a ring satisfying
Conditions~\ref{item:condition_on_mapring_finite_generation} and
\ref{item:condition_on_mapring_containing_image} of \S\ref{sect:map_ring}.
Eventually we need $\mapring \Phi$ be a map ring of $\Phi$, but for the sole
purpose of proving Proposition~\ref{prop:image_and_preimage_does_not_depend}
we consider this slightly more general object.

Setting $V(\Phi) = \Spec \mapring\Phi$, we have two natural morphisms:
\[
	\CC^m \stackrel{p_{\Phi}}{\longleftarrow} V(\Phi)
	\stackrel{q_{\Phi}}{\longrightarrow} \CC^n,
\]
where ${p_{\Phi}}^*$ is the inclusion of $\cox{X}$ in $\mapring\Phi$ and
${q_{\Phi}}^*$ is defined by mapping $y_i$ to $\Phi^* y_i \in \mapring \Phi$.
We treat $V(\Phi)$ informally as a correspondence (even though it is
not constructed in the product, and in any case it is finite over
$\Reg \Phi$ but not necessarily over $\CC^m$).
Using this, we can define (set-theoretic)
image and preimage of subsets in a natural way.

\begin{defin}\label{defin:image_preimage_under_multi_valued_map}
	Let $\Phi\colon \CC^m \multito \CC^n$  be a multi-valued map.
	Let $A \subset \Reg \Phi$ be a subset.
	The \textbf{image of $A$ under $\Phi$} is the subset of $\CC^n$ defined by
        \[
           \Phi(A):=  q_{\Phi}\left(\inv{p_{\Phi}}(A)\right).
        \]
	Let $B \subset \CC^n$ be a subset.
	The \textbf{preimage of $B$ under $\Phi$} is the subset of $\Reg \Phi$ defined by
	\[
	  \Phi^{-1}(B):=  p_{\Phi}\left( \inv{q_{\Phi}} (B) \right).
	\]
\end{defin}

In Section~\ref{sect:image_of_point} below, we explain how to evaluate
a multi-valued function $\Phi$ at a point $\xi$; this agrees with
the notion of image just discussed for $A=\{\xi\}$:
$\Phi(A) = \{\Phi(\xi)\}$. When
$\Phi$ is a single-valued map, these definitions give the usual image and
preimage under a rational map.

Since $q_\Phi$ is continuous and $p_\Phi\colon V(\Phi) \to \Reg \Phi$ is finite
and locally free (and thus closed by \cite[Ex.~II.3.5(b)]{hartshorne}
 and open by \cite[Lemmas~042S and 02KB]{stacks_project}), preimage behaves well
with respect to the Zariski topology.
\begin{prop}
\label{prop:preimage_of_closed_is_closed}
If $B \subset \CC^n$ is open, then $\Phi^{-1}(B) \subset \CC^m$ is open.
If $B \subset \CC^n$ is closed, then $\Phi^{-1}(B) \subset \Reg \Phi$ is closed.
\end{prop}

\begin{prop}\label{prop:image_and_preimage_does_not_depend}
  The definitions of image and preimage above are independent of
  the choice of $\mapring \Phi$ satisfying conditions
  \ref{item:condition_on_mapring_finite_generation} and
  \ref{item:condition_on_mapring_containing_image} of \S\ref{sect:map_ring}.
\end{prop}

\begin{prf}
  All the rings satisfying conditions
  \ref{item:condition_on_mapring_finite_generation} and \ref{item:condition_on_mapring_containing_image}
  must contain
  \[
    \mapring \Phi_{min}:= \CC[\Reg \Phi] [\Phi^* y_1, \dotsc, \Phi^* y_n].
  \]
  On the other hand $\mapring \Phi_{min}$ itself satisfies these two conditions.
  So for any $\mapring \Phi$ we have the commutative diagram:
  \[
    \xymatrix{
                        &    &\Spec \mapring \Phi \ar@{->>}[d] \ar@{->>}[ddll] \ar[ddrr]\\
                        &    &\Spec \mapring \Phi_{min}  \ar@{->>}[dll] \ar[rrd]^{}&\\
              \Reg \Phi & & & \text{\phantom{blablabla}} & \CC^n  \\
             }
  \]
  and since the middle vertical arrow is epimorphic it follows that
  it does not matter which way around one carries the subset between $\Reg \Phi$ and $\CC^n$.
\end{prf}

In \S\ref{sect:image}--\ref{sect:preimage} we explain how to
consider scheme-theoretic image and preimage under certain multi-valued maps.
This is more delicate since the scheme structure of the image or preimage may
depend on the choice of map ring $\mapring{\Phi}$.

\begin{prop}\label{prop:image_of_mvm_is_given_by_kernel}
	The ideal of the Zariski closure $\overline{\Phi(\Reg \Phi)}$ of
	$\Phi(\Reg \Phi)$ is the kernel of $\Phi^*$.
\end{prop}
\begin{prf}
  Since the image of $p_\Phi$ is exactly $\Reg\Phi$,
  \[
    \Phi(\Reg \Phi) = q_{\Phi}\left( \inv{p_{\Phi}}( \Reg \Phi) \right) =
    q_{\Phi} \bigl(V(\Phi)\bigr).
  \]
  Now $f \in \cox Y$ vanishes on $q_{\Phi} \bigl(V(\Phi)\bigr)$ if and only if
  $0 = {q_{\Phi}}^* f = \Phi^* f$.
\end{prf}

\subsubsection{Image of a single point}\label{sect:image_of_point}

We consider the image of a single closed point under a multi-valued map and
prove that it can be computed by evaluation with a little care.

\begin{example}\label{example:image_of_a_point}
 Consider the following multi-valued map:
\begin{align*}
   \Phi\colon \CC^2  & \multito \CC^2 \\
          (s,t) & \mapsto (\sqrt[6]{s}, \sqrt[2]{s^3} (t^2+s)).
\end{align*}
The image of the point $(64,-1)$ consists of the 6 points
\[
 (2 \epsilon_6 , 512 \epsilon_6^3 (1+64) )
\]
as $\epsilon_6$ runs over the $6$-th roots of unity.
On the other hand, the point $(2,-512 \times 65)$ is not in the image of $(64,-1)$,
even though $2  = \sqrt[6]{64}$ and $-512\times 65 = - \sqrt{64^3} ((-1)^2 + 64)$.
\end{example}

The crucial observation in this example is that the irrational parts
$\sqrt[6]{s}$ and $\sqrt[2]{s^3}$ are algebraically dependent: in fact,
\[
 (\sqrt[6]{s})^9 = \sqrt[2]{s^3},
\]
so $\sqrt[2]{s^3}$ is already in the extension ring $\CC[s,t][\sqrt[6]{s}]$.
(Some choice of the sixth root must have been made, and here we enforce
that choice on the whole calculation.)

Choose a point $\xi \in \Reg \Phi$ and let $ev_\xi\colon \CC[\Reg \Phi] \to \CC$
be the evaluation map. Consider the diagram:
 \[
    \xymatrix{
    \CC[\CC^n] \ar[rrd]^{\Phi^*}  && \CC[\Reg \Phi] \ar@{^{(}->}[d]  \ar[rr]^{ev_\xi} && \CC  \ar@{=}[d]\\
      && \mapring{\Phi}   \ar@{.>}[rr] ^{\widetilde{ev}_\xi}                          &&\overline{\CC}
             }
 \]
The extensions $\widetilde{ev}_\xi$ exist and they are precisely determined
by any choice of roots of images of the distinguished generators
(see Lemma~\ref{lem:extending_to_simple_extensions}).

\begin{thm}
   Let $\xi\in\Reg\Phi$. Then
   $\Phi(\xi)$ is precisely the set of all those $\eta \in \CC^n$
   whose maximal ideal is a kernel of $\widetilde{ev}_\xi \circ \Phi^*$
   for some extension $\widetilde{ev}_\xi$.
\end{thm}

\begin{prf}
Let $\gotm_\xi = \ker ev_\xi \ideal \Reg \Phi$ be the maximal ideal of $\xi$.
First assume $\eta\in \Phi(\xi)$.
Then there exist a point $\zeta \in V(\Phi)$, such that $q_{\Phi}(\zeta) = \eta$
and $p_{\Phi}(\zeta) =\xi$.
So if $\gotm_\zeta \ideal \mapring \Phi$ is the maximal ideal of $\zeta$,
then $\gotm_\zeta \supset  \langle \gotm_\xi \rangle \ideal \mapring \Phi$.
Consider $ev_\zeta \colon \mapring \Phi \to \mapring \Phi/ \gotm_\zeta \simeq \CC$.
Now clearly $ev_\zeta|_{\CC[\Reg \Phi]}$ is a (nonzero) ring homomorphism,
whose kernel contains the maximal ideal $\gotm_\xi$.
So
\[
  ev_\zeta|_{\CC[\Reg \Phi]} = ev_\xi,
\]
and so $\widetilde{ev}_\xi := ev_\zeta$ is an extension of $ev_\xi$
such that its kernel of $\widetilde{ev}_\xi \circ \Phi^*$ is $\gotm_\eta$.

Now assume we have an extension $\widetilde{ev}_\xi$.
Let $\gotm_\zeta$ be its kernel.
Clearly $\langle \gotm_\xi \rangle \subset \gotm_\zeta$,
so $p_{\Phi}(\zeta) =\xi$ and therefore $q_{\Phi}(\zeta) \in \Phi(\xi)$.
\end{prf}

%% file: toricmapsdescriptions.tex

\section{Descriptions of maps}
\label{sect:descriptions}

Consider two toric varieties $X$ and $Y$ and their
Cox covers $\CC^m = \Spec \cox X$ and $\CC^n = \Spec \cox Y$,
where $\cox X = \CC[\fromto{x_1}{x_m}]$ and $\cox Y = \CC[\fromto{y_1}{y_n}]$.
In this section, we show how to use multi-valued maps
$\Phi\colon \CC^m \multito \CC^n$ to describe rational maps
$\varphi\colon X \dashrightarrow Y$.
In particular, we address

\renewcommand{\theenumi}{\textnormal{\arabic{enumi}}}
\renewcommand{\labelenumi}{\textnormal{\theenumi.}}
\begin{enumerate}
\addtolength{\itemsep}{-0.5\baselineskip}
\item
what it means for a multi-valued map $\Phi$ to describe a rational
map $\varphi$.
\item
which multi-valued maps describe rational maps at all.
\item
that every rational map can be described by a multi-valued map.
\item
a class of multi-valued maps that describe rational maps particularly well, or completely.
\end{enumerate}
An algorithm for finding such a complete description
$\Phi$ of a given $\varphi$ is implicit in the proofs.

\subsection{The agreement locus}
\label{sect:agreement}

Let $\Phi\colon \CC^m \multito \CC^n$ be a multi-valued map.
It fits into a diagram
\begin{equation}
\label{eq:Phi_setup}
   \xymatrix{
     \CC^m \xymultito{r}^{\Phi}
       \ar@{-->}[d]^{\pi_X} & \CC^n\ar@{-->}[d]^{\pi_Y}\\
	    X&  Y }
\end{equation}
The regular locus $\Reg\Phi \subset \CC^m$ of $\Phi$,
where its denominators do not vanish as in
Definition~\ref{defin:regular_locus}, contains a finer subset,
the locus where $\pi_Y \circ \Phi$ is a well-defined map of sets:
\[
  \Reg_Y\Phi:=
    \set{ \xi \in \Reg \Phi \mid
    \Phi(\xi) \cap \Reg \pi_Y \ne \emptyset \text{ and }
 	\# \pi_Y(\Phi(\xi)) = 1}.
\]
This locus $\Reg_Y\Phi$ may be empty. On the other hand, if $\Reg_Y\Phi$
contains a nonempty open subset, then we regard $\Phi$ as being adapted to
$Y$; under this assumption, it makes sense to ask where $\Phi$ agrees
with a rational map $X\dashrightarrow Y$.

\begin{defin}
 Given a multi-valued map $\Phi\colon \CC^m \multito \CC^n$
 and a rational map $\varphi\colon X \dashrightarrow Y$, in the notation above,
 the \textbf{agreement locus of $\Phi$ and $\varphi$} is
 \[
   \agree{\Phi}{\varphi} = \set{
      \xi\in \Reg_Y\Phi \cap {\pi_X}^{-1}(\Reg\varphi) \mid
      \pi_Y\circ\Phi(\xi) = \varphi\circ\pi_X(\xi) }.
 \]
\end{defin}
In other words, the agreement locus is the set of points where both
compositions $\pi_Y \circ \Phi$ and $\varphi \circ \pi_X$ are well-defined
as maps of sets and they have the same values. 

\begin{rem}
At this point, even if $\Reg_Y\Phi$ contains an open dense subset this agreement locus 
  could be contained in a proper closed subset, 
  equal to a finite number of points, 
  or even empty.
In this paper, we are interested in the case, when the agreement locus contains an open dense subset
  (see Definition~\ref{defin:description}), but it is easy to imagine, it is only dense in some subvariety $Z$ of $X$,
  (for instance, $Z$ could be a Mori dream space with $m$ generators of the Cox ring of $Z$,
   and $X$ could be a toric variety containing $Z$).
Then we could study descriptions of maps between $Z$ and $Y$ (or yet another subvariety of $Y$).
We would not comment further on this possibility.
\end{rem}

Perhaps the next definition does not seem surprising, but it really is the key one in this paper.
As written it is purely set-theoretic---what is surprising is that with our general restrictions on the multivalued maps, 
  the set-theoretic properties suffice to prove many algebraic conditions, 
  such as Propositions~\ref{prop:equivalence_of_homog_conds} and \ref{prop:Y_Q_factorial_then_max_agr_not_necessary}.

\begin{defin}\label{defin:description}
	We say \textbf{$\Phi$ is a description of $\varphi$},
	or that \textbf{$\Phi$ describes $\varphi$},
	if $\agree \Phi \varphi$ contains an open dense subset of $\CC^m$.
\end{defin}

When we have a multi-valued map $\Phi$ that describes a rational map $\varphi$,
we say that $\varphi$ is given in Cox coordinates by
\begin{align*}
\varphi\colon X & \dashrightarrow Y \\
  x & \mapsto \left[\left(\Phi^* y_1\right)(x), \ldots, \left(\Phi^* y_n\right)(x)\right]
\end{align*}
leaving implicit that $\Phi^*y_i$ is really only evaluated on some
$\xi\in\agree{\Phi}{\varphi}$ for which $x=[\xi]$.

Section~\ref{sect:Motivating_example} has several examples of descriptions
of maps, and here is another.

\begin{example}
The diagonal embedding of $\PP^1 \hookrightarrow \PP^1 \times \PP^1$
has the description
 \begin{align*}
 \Phi\colon
    [x_1,x_2] & \mapsto [x_1,x_2, x_1,x_2].
 \end{align*}
In this case,
$\ker \Phi^* = \left\langle y_1 - y_3, y_2 - y_4\right\rangle$
is not a homogeneous ideal with respect to the gradings
\[
\begin{pmatrix} 1 & 1 & 0 & 0 \\ 0 & 0 & 1 & 1 \end{pmatrix}
\text{ on the Cox coordinates }\fromto{y_1}{y_4},
\]
in contrast to the case of projective spaces.
It is easy to see in this case that the homogeneous part of the
kernel is $\left< y_1y_4 - y_2y_3 \right>$, and that this defines
the image of the embedding.
\end{example}

\subsection{Homogeneity and relevance conditions}

 We determine when a multi-valued map $\Phi$ between
the Cox covers of two toric varieties $X$ and $Y$
describes a rational map $X\dashrightarrow Y$. 
First we show the equivalence of four conditions analogous to
the usual homogeneity conditions
for maps between projective spaces.  Together, they are referred to as
\emph{the homogeneity condition}. 

\renewcommand{\labelenumi}{\textnormal{(\theenumi)}}

\begin{prop}
\label{prop:equivalence_of_homog_conds}
Let $\Phi$ be a multi-valued map as in \eqref{eq:Phi_setup} above and
consider the set
\[
T=\set{ y_i \mid i \in \setfromto{1}{n} \text{ and }\ \Phi^* y_i \ne 0}
\]
of Cox generators of $\cox Y$ that pull back nontrivially
under $\Phi$.
The following conditions are equivalent:
\begin{enumerate}
\renewcommand{\theenumi}{\textnormal{\ref{item:homogeneity_cond}\arabic{enumi}}}
\item
\label{item:homogeneity_cond_all}
If $q\in \coxfield{Y}$ is homogeneous and $\Phi^*q$ is defined, then
$\Phi^*q$ is a homogeneous multi-valued section on $X$.
\item
\label{item:homogeneity_cond_on degree_0}
If $q\in \functionfield{Y}$ and $\Phi^*q$ is defined, then
$\Phi^*q\in\functionfield{X}$.
\item
\label{item:homogeneity_cond_monomial}
There exist rational monomials $t_1,\dotsc,t_k$ generating
$\functionfield{Y} \cap \CC(T) = \CC(T)^0$ as a field extension of $\CC$
such that $\Phi^* t_i$ are homogeneous single-valued sections of degree $0$.
\item
\label{item:homogeneity_cond_points}
For all $\xi,\xi' \in \Reg\Phi$ with $\xi' \in G_X \cdot \xi$,
if $\eta \in \Phi(\xi)$ and $\eta' \in \Phi(\xi')$ then $\eta' \in G_Y \cdot \eta$.
\renewcommand{\theenumi}{\textnormal{\ref{item:homogeneity_cond_points}${}^\prime$}}
\item
\label{item:homogeneity_cond_points_on_open}
There exists an open dense subset $U \subset \Reg \Phi$, such that
for all $\xi,\xi' \in U$ with $\xi' \in G_X \cdot \xi$,
if $\eta \in \Phi(\xi)$ and $\eta' \in \Phi(\xi')$ then $\eta' \in G_Y \cdot \eta$.
\end{enumerate}
\end{prop}

\ref{item:homogeneity_cond_on degree_0} is the usual treatment
of rational maps $X \dashrightarrow Y$ as a map of function fields,
taking care with the domain in case the map is not dominant.
We use this to recover a rational map from a description
(see Theorem~\ref{thm:Phi_homog_and_relevant_is_a_description}),
and it is also convenient in calculations, as in the introduction.
\ref{item:homogeneity_cond_monomial} is the same condition
expressed for a finite number of generators, which is useful when
deciding whether an expression determines a rational map; we
also use it to construct a description of a rational map
(see Theorem~\ref{thm:existence_of_description}).

\ref{item:homogeneity_cond_all} is used
to prove Proposition~\ref{thm:disagreement_is_of_codimension_1},
calculating the dimension of the complement of the agreement locus.
It is not much help for deciding whether a
given expression determines a rational map, as the example in
\S\ref{sect:mvmls} illustrates.
\ref{item:homogeneity_cond_points} is the geometric condition
that  $\Phi$ maps $G_X$-orbits into $G_Y$-orbits. This is a closed
condition, which is expressed as \ref{item:homogeneity_cond_points_on_open}.
\ref{item:homogeneity_cond_points} and \ref{item:homogeneity_cond_points_on_open}
are used to give conditions for a multi-valued map
to be a description of some rational map (see
Proposition~\ref{prop:description_satisfy_homog_and_relev} and
Theorem~\ref{thm:Phi_homog_and_relevant_is_a_description})
and in the calculations of agreement locus in \S\ref{sect:agreement_revisited}.

\begin{prf}[ of Proposition~\ref{prop:equivalence_of_homog_conds}]
Suppose \ref{item:homogeneity_cond_all} holds for $\Phi$.
Let $V \subset \coxfield{Y}$ be the subspace of
homogeneous sections of degree $0$ for which the pullback by $\Phi$ is
defined. Denote the restriction of $\Phi^*$ to $V$ by
$i\colon V \to \overline{\coxfield{X}}$. Since $i(1)=1$ is rational
and has degree $0$,
Proposition~\ref{prop:image_has_constant_degree_and_irrational_part} implies
that all elements of $i(V)$ are rational and of degree $0$.
Therefore \ref{item:homogeneity_cond_on degree_0} holds for $\Phi$.

Suppose \ref{item:homogeneity_cond_on degree_0} holds. Since
$\CC(T)^0 \subset \coxfield{Y}^0$, any monomial generating set
$\fromto{t_1}{t_k}$ of $\CC(T)^0$ satisfies
\ref{item:homogeneity_cond_monomial} for $\Phi$.

Suppose \ref{item:homogeneity_cond_monomial} holds for $\Phi$; we show
that \ref{item:homogeneity_cond_all} holds.
Let $q\in \coxfield{Y}$ be any homogeneous function. Write
\[
   q = \frac{\mu_1 + \dotsb + \mu_{\alpha}}{\nu_1 + \dotsb + \nu_{\beta}},
\]
where the $\mu_i$ and $\nu_j$ are monomial terms in $\cox Y$ with
$\deg \mu_i =d_1$ and  $\deg \nu_j =d_2$ for all $i$ and~$j$. Assume that
$\Phi^*(\nu_1 + \dotsb + \nu_{\beta}) \ne 0$, so $\Phi^*q$ is defined.

Certainly each $\Phi^*\mu_i$ is a homogeneous multi-valued section.
Therefore the Laurent monomial $\mu_{i_1}/\mu_{i_2}$ is homogeneous of degree $0$
and either $\Phi^*(\mu_{i_1})=0$ or $\Phi^*(\mu_{i_2})=0$
or $\Phi^*(\mu_{i_1}/\mu_{i_2})$ is a nonzero homogeneous degree~$0$
rational section in $\functionfield{X}$. In particular, for every $i$,
\[
    \Phi^*(\mu_i) = f_i \cdot \mvsa
\]
where $\mvsa \in \overline{\coxfield{X}}$ is a fixed homogeneous multi-valued
section (independent of $i$) and $f_i \in \functionfield{X}$. So
\[
    \Phi^*(\mu_1  + \dotsb + \mu_{\alpha}) = (f_1 + \dotsb + f_{\alpha}) \mvsa.
\]
Similarly,
$\Phi^*(\nu_1  + \dotsb + \nu_{\beta}) = (g_1 + \dotsb + g_{\beta}) \mvsb \not= 0$,
for some $\mvsb \in \overline{\coxfield{X}}$ and $g_j\in\functionfield{X}$.
So
\[
   \Phi^*(q) = h \cdot {\varepsilon}
\]
where $\varepsilon = \mvsa/\mvsb \in \overline{\coxfield{X}}$ is homogeneous
and $h = (\sum f_i) / (\sum g_j) \in \functionfield{X}$.
So $\Phi^*(q)$ is homogeneous and \ref{item:homogeneity_cond_all} holds.

It remains to prove the equivalence of \ref{item:homogeneity_cond_on degree_0},
  \ref{item:homogeneity_cond_points} and \ref{item:homogeneity_cond_points_on_open}.

    Suppose \ref{item:homogeneity_cond_on degree_0} holds.
    Let $\xi \in \Reg \Phi$  and consider $G_X \cdot \xi$.
    The claim of \ref{item:homogeneity_cond_points} is that $\Phi(G_X \cdot\xi)$ is contained in one $G_Y$ orbit.
    Let $A \subset \setfromto{1}{n}$ be the set of those $i$, that $\Phi^* y_i$ vanishes at $\xi$.
    Since $\Phi^* y_i$ is homogeneous, if $i\in A$,
        then $\Phi^* y_i$ vanishes identically on the orbit $G_X \cdot \xi$
        and if $i\notin A$, then $\Phi^* y_i$ is nowhere zero on $G_X \cdot \xi$.
    Thus $\Phi(G_X \cdot\xi)$ is contained in a torus $T \subset \CC^n$
        given by $y_i =0 $ for $i \in A$ and $y_j \ne 0$ for $i \notin A$.
    By definitions, the group $G_Y$ preserves $T$.
    We consider the quotient torus $T/G_Y = \Spec({\ccO_T}^{G_Y})$ and obtain the following diagram:
    \begin{equation*}
       \xymatrix{
                                     &       &\inv{p_{\Phi}}(G_X \cdot \xi) \ar@{^(->}[d] \ar@{->>}[ddll] \ar[ddrr]\\
                                     &       &V(\Phi) \ar@{->>}[dl]_{p_{\Phi}} \ar[rd]&\\
         G_X \cdot \xi \ar@{^(->}[r] & \Reg \Phi  & &   \CC^n & T  \ar@{_(->}[l] \ar@{->>}[r] & T/G_Y \\
                }
    \end{equation*}
    The claim of \ref{item:homogeneity_cond_points} is that image of $\inv{p_{\Phi}}(G_X \cdot \xi)$
       under the composed map $\inv{p_{\Phi}}(G_X \cdot \xi \to  T/G_Y$ is a single point.
    Equivalently, for any regular function on $T/G_Y$, the pullback is a constant function on
       $\inv{p_{\Phi}}(G_X \cdot \xi)$.
    A regular function on $T/G_Y$ is a $G_Y$-invariant regular function on $T$. Any $G_Y$-invariant function on $T$ is a restriction of a degree zero rational function $q$ on $\CC^n$, whose pullback $\Phi^*q$ is defined,
       and is regular at $\xi$.
    By \ref{item:homogeneity_cond_on degree_0} we have $\Phi^*q \in \functionfield{X}$,
       in particular the pullback of $q$ to $\inv{p_{\Phi}}(G_X \cdot \xi)$
       is equal to ${p_{\Phi}}^*(\Phi^*q|_{G_X \cdot \xi})$.
    Since $\Phi^*q|_{G_X \cdot \xi}$ is $G_X$-invariant, $\Phi^*q|_{G_X \cdot \xi} \equiv \Phi^*q(\xi)$, i.e.,
       it is a constant function. Its pullback by ${p_{\Phi}}^*$ is therefore also a constant function
       on   $\inv{p_{\Phi}}(G_X \cdot \xi)$ and the claim of \ref{item:homogeneity_cond_points} is proved.

    If \ref{item:homogeneity_cond_points} holds, then clearly \ref{item:homogeneity_cond_points_on_open} holds too.

    Finally, suppose \ref{item:homogeneity_cond_points_on_open} holds.
    Let $q\in \functionfield{Y}$ be such that $\Phi^*q$ is defined.
    Suppose $\xi \in U$ is general.
    The possible values taken by $\Phi^*q$ at $\xi$ are simply those values taken by $q$
       at the points of the image set $\Phi(\xi)$.
    Setting $\xi'=\xi$ in \ref{item:homogeneity_cond_points_on_open}
       shows that $\Phi(\xi)$ is contained in a single $G_Y$-orbit.
    and so
    $
       \Phi^*q(\xi) = \set{q(\eta) \mid \eta \in \Phi(\xi)}
    $
    is a single point.
    Therefore $\Phi^* q \in \coxfield{X}$ by Proposition~\ref{prop:single_valued_sections}.
    In any case, for any $\xi$, $\xi'$ as in \ref{item:homogeneity_cond_points_on_open},
    \[
       (\Phi^*q)(\xi) = q(\Phi(\xi)) = q(\Phi(\xi')) = (\Phi^*q)(\xi')
    \]
    since $q$ is constant on $G_Y$-orbits. That is, $\Phi^* q$ is contant on a general $G_X$ orbit,
       so $\Phi^* q$ is $G_X$-invariant.
\end{prf}

Next we note the equivalence of another three conditions,
jointly referred to as \emph{the relevance condition}.

\begin{prop}
\label{prop:equivalence_of_rel_conds}
Let $\Phi$ be a multi-valued map as in \eqref{eq:Phi_setup} above and
consider the set
\[
R_0=\set{ \rho_i \in \Sigma^{(1)}_Y \mid i \in \setfromto{1}{n}
\text{ and }\ \Phi^* y_i = 0}
\]
of rays of the fan $\Sigma_Y$ of $Y$ which correspond to Cox
generators of $\cox Y$ that pull back trivially under~$\Phi$.
The following conditions are equivalent:
\begin{enumerate}
\renewcommand{\theenumi}{\textnormal{\ref{item:rel_cond}\arabic{enumi}}}
\item
\label{item:rel_cond_points}
The image of $\Phi$ is not contained in the irrelevant locus of $Y$.
\item
\label{item:rel_cond_kernel}
$\ker \Phi^*$ does not contain the irrelevant ideal $B_Y$ of $Y$,
that is, $\ker \Phi^*$ is a relevant ideal.
\item
\label{item:rel_cond_cones}
The rays of $R_0$ are all contained in a single cone of $\Sigma_Y$.
\end{enumerate}
\end{prop}

\begin{prf}
The equivalence of the first two conditions is immediate (even acknowledging
the multi values of $\Phi$).
If $\sigma$ is a maximal cone of $\Sigma_Y$ containing all the rays
of $R_0$, then the standard generator $m_\sigma\in B_Y$ determined by
$\sigma$ satisfies $\Phi^*m_\sigma\not=0$. Thus $m_\sigma$ is not contained
in $\ker\Phi^*$, and so neither is $B_Y$. Conversely, if there is
no maximal cone containing all the rays of $R_0$, then every standard
generator of $B_Y$ contains at least one such ray. Therefore
$B_Y\subset \ker\Phi^*$.
\end{prf}

\renewcommand{\theenumi}{\textnormal{\Alph{enumi}}}

\begin{defin}
\label{defin:homogrel}
Let $\Phi$ be a multi-valued map as in \eqref{eq:Phi_setup} above.
\begin{enumerate}
\item
\label{item:homogeneity_cond}
  We say that $\Phi$ satisfies the \textbf{homogeneity condition}
  if any of the equivalent conditions
  \ref{item:homogeneity_cond_all},
  \ref{item:homogeneity_cond_on degree_0},
  \ref{item:homogeneity_cond_monomial},
  \ref{item:homogeneity_cond_points},
  \ref{item:homogeneity_cond_points_on_open}
  of Proposition~\ref{prop:equivalence_of_homog_conds} hold for $\Phi$.
\item
\label{item:rel_cond}
  We say that $\Phi$ satisfies the \textbf{relevance condition}
  if any of the equivalent conditions
  \ref{item:rel_cond_points},
  \ref{item:rel_cond_kernel},
  \ref{item:rel_cond_cones}
  of Proposition~\ref{prop:equivalence_of_rel_conds} hold for $\Phi$.
\end{enumerate}
\end{defin}

\renewcommand{\theenumi}{\textnormal{(\roman{enumi})}}
\renewcommand{\labelenumi}{\theenumi}

\begin{prop}\label{prop:description_satisfy_homog_and_relev}
If $\Phi$ is  a description of a rational map
$\varphi \colon X \dashrightarrow Y$, then $\Phi$ satisfies the
homogeneity and relevance conditions of Definition~\ref{defin:homogrel}.
\end{prop}

\begin{prf}
By Definition~\ref{defin:description} of description, $\pi_Y \circ \Phi$
is defined on an open subset of $\CC^m$, so $\Phi(x)$ cannot be contained
in the irrelevant locus for those points. Therefore $\Phi$ satisfies the
relevance condition \ref{item:rel_cond_points}.

Since $\Phi$ is a description the agreement locus $\agree{\Phi}{\varphi}$
contains an open dense subset of $\Reg\Phi$.
The homogeneity condition~\ref{item:homogeneity_cond_points_on_open} is
satisfied on this set.
\end{prf}

The converse is the main point: the homogeneity and relevance
conditions guarantee that a multi-valued map is a description of
a uniquely-determined rational map.

\begin{thm}\label{thm:Phi_homog_and_relevant_is_a_description}
Let $\Phi$ be a multi-valued map as in \eqref{eq:Phi_setup} above
that satisfies the homogeneity and relevance conditions of
Definition~\ref{defin:homogrel}.
\begin{enumerate}
\item
\label{item:Phi_determines_rational_map}
By its action on rational functions, $\Phi^*$ naturally determines
a rational map $\varphi \colon X \dashrightarrow Y$.
\item
\label{item:Phi_describes_phi_iff_phi_eq_psi}
$\Phi$ is a description of some map $\psi\colon X\dashrightarrow Y$
if and only if $\psi = \varphi$.
\end{enumerate}
\end{thm}

\begin{prf}
 To prove \ref{item:Phi_determines_rational_map} first note
 that $\gotp:=\homog{(\ker \Phi^*)} \ideal \cox Y$ is homogeneously prime by
 Proposition~\ref{prop:homog_of_kernel_is_homog_prime},
 so that the following localisation makes sense:
 \[
   R:=\hlocalise{\cox Y}{\gotp}.
 \]
 We claim $\Phi^*$ naturally determines a ring homomorphism:
    \[
      \functionfield{Y} \supset R \stackrel{\Phi^*}{\longrightarrow}\functionfield{X}.
    \]
 This is because by definition
 \[
   R = \set{ \frac{f}{g} \mid f,g \in \cox Y , \ g \notin \gotp
   \text{ and $f$, $g$ are homog.~of the same degree}   }.
 \]
 Since $\gotp$ is generated by all homogeneous sections in $\ker \Phi^*$,
 we can also replace the condition $g \notin \gotp$ with $g \notin \ker \Phi^*$:
 \[
   R= \set{  \frac{f}{g} \mid f,g \in \cox Y , \ g \notin \ker \Phi^*
\text{ and $f$, $g$ are homog.~of the same degree}   }.
 \]
 In particular, if $\frac{f}{g} \in R$, then the pull-back by $\Phi$ is defined
 (because $\Phi^* g$ is not zero).

 By the homogeneity condition \ref{item:homogeneity_cond_on degree_0},
 \[
   \Phi^* \left( \frac{f}{g} \right) \in \cox{X}_0 \cong \functionfield{X}.
 \]
 So we have a ring homomorphism $R \to \functionfield{X}$ as claimed.

 Note, that by Lemma~\ref{lem:R_I_generates_KX} together with the relevance condition \ref{item:rel_cond_kernel},
 $R$ and $R^{-1}$ together generate $\functionfield{Y}$.
 Hence by Proposition~\ref{prop:rational_maps}\ref{item:ring_homomorphism_determines_map}
 the ring homomorphism
 \[
   \Phi^*\colon R  \to \functionfield{X}
 \]
  determines a rational map $\varphi\colon X \dashrightarrow Y$ which is characterised
  by its action on rational functions $q\in\functionfield{Y}$ being $\varphi^*(q) = \Phi^*(q)$.

 Next we have to prove that $\Phi$ describes $\varphi$.
 Consider the open subset
 \[
   U=\set{\xi \in \Reg \Phi \mid \xi \notin \irrel X, \Phi(\xi) \nsubseteq \irrel Y}
 \]
 of $\Reg \Phi$; note that
 it contains a non-empty open subset of $\CC^m$ by the relevance condition.
 Choose any $\xi \in U$.
 By the homogeneity condition \ref{item:homogeneity_cond_points},
 $\pi_Y(\Phi(\xi))$ is a single point $y$.
 We claim $y = \varphi(\pi_X(\xi))$,
 so that $\xi\in \agree{\Phi}{\varphi}$.

 To prove the claim, we set $x=[\xi] = \pi_X(\xi)$ and evaluate rational functions
 $q \in \functionfield{Y}$ at $\varphi(x)$ and~$y$:
 \[
   q(\varphi(x)) = (\varphi^* q) (x) = (\Phi^*q) ([\xi]) = q([\Phi(\xi)]) = q(y).
 \]
 So no rational function on $Y$ can distinguish between $\varphi(x)$ and $y$ and therefore $y= \varphi(x)$.
 Hence $U \subset \agree{\Phi}{\varphi}$ and $\Phi$ describes $\varphi$.

 Finally we note that if $\psi\colon X \dashrightarrow Y$ is another rational map
 which is also described by $\Phi$, then for $\xi \in \agree{\Phi}{\psi}$  with $x = [\xi]$
 and for a rational function $q\in K(Y)$ we have
 \[
  (\psi^* q) (x) =   q(\psi(x)) =    q([\Phi(\xi)]) = (\Phi^*q) (\xi) = (\varphi^*q) (x).
 \]
 Hence $\psi^* = \varphi^*$ and therefore $\psi =\varphi$.
\end{prf}

\begin{cor}\label{cor:torus_orbit_containment}
Let $\Phi$ be a description of a rational map $\varphi\colon X \dashrightarrow Y$.
\begin{enumerate}
\item
Let $\sigma \in \Sigma_Y$ be the smallest cone which contains all rays
whose corresponding coordinate $y_i$ is pulled back to $0$ by $\Phi$.
Then the closed toric stratum corresponding to $\sigma$ is
the smallest closed stratum of $Y$ that contains $\varphi(X)$.
\item
The assignment
\[
    \Psi^*y_i :=
    \begin{cases}
       0         & \text{if the $i$-th ray of $\Sigma_Y$ is in $\sigma$,}\\
       \Phi^*y_i & \text{otherwise}
    \end{cases}
\]
defines a multi-valued map $\Psi$, and $\Psi$ also describes $\varphi$.
\item
If, furthermore, $Y$ is $\QQ$-factorial, then $\Phi^*y_i = 0$ if and only if
$\varphi(X)$ is contained in the locus $y_i = 0$.
\end{enumerate}
\end{cor}

\begin{prf}
When $\pi_Y\colon \CC^n\dashrightarrow Y$ is a geometric quotient,
$\eta\in\CC^n$ is a semistable point and $y_i$ is a Cox coordinate, then
\[
y_i(\eta) = 0 \iff \pi_Y(\eta) \in \Supp (y_i),
\]
where $(y_i)$ is the divisor on $Y$ corresponding to $y_i$.
So if $Y$ is $\QQ$-factorial and $\xi\in \agree \Phi \varphi$, then
  \begin{align*}
    \varphi(\pi_X(\xi)) \in \Supp(y_i) & \iff \pi_Y(\Phi(\xi)) \in \Supp(y_i) \\
                                         & \iff y_i(\Phi(\xi)) =0 \\
                                         & \iff (\Phi^* y_i)(\xi) = 0,
  \end{align*}
which proves the final statement.

If the quotient $\pi_Y$ is not geometric, then
we have only
\[
y_i(\eta) = 0 \ \Longrightarrow\ \pi_Y(\eta) \in \Supp (y_i),
\]
so that $\varphi(X)$ is contained in the intersection of the
supports of divisors $(y_i)$ for those $y_i$ with $\Phi^*y_i=0$.
On $Y$, this intersection is the toric stratum corresponding to
the cone $\sigma$. In particular, $\varphi(X)\subset \Supp(z)$
for every Cox coordinate $z$ corresponding to a ray of $\sigma$,
whether or not $\Phi^*z$ is zero.
So for any $\xi\in\agree\Phi\varphi$ we have
 $\pi_Y(\Phi(\xi)) = \pi_Y(\Psi(\xi))$.
Since $\Psi$ and $\Phi$ therefore have the same agreement locus,
$\Psi$ is also a description of $\varphi$.
\end{prf}

\subsection{Existence of descriptions}\label{sect:existence_of_description}

The previous section shows that descriptions of rational maps are
characterised by the homogeneity and relevance conditions.
Now we show that every rational map does have a description.

\begin{thm}\label{thm:existence_of_description}
   Let $\varphi\colon X \dashrightarrow Y$ be a rational map of toric varieties.
   Then there exists a description $\Phi\colon\CC^m \multito \CC^n$ of $\varphi$.
\end{thm}

\begin{prf}
   We construct $\fromto{\Phi^* y_1}{\Phi^* y_n}$ inductively.
   Set  $\Phi^*y_i =0$ if and only if $\varphi(X) \subset \Supp(y_i)$.
   So assume without loss of generality that $\varphi(X)$ is contained in
   $y_1 = \ldots = y_s =0 $ only, for some $s \in \setfromto{0}{n}$.
   Fix $\Phi^*y_i=0$ for $i \in \setfromto{1}{s}$.

   Assume $\Phi^*y_i$ is fixed for all $i \in \setfromto{1}{k-1}$
   for some $k \in \setfromto{s+1}{n}$.
   Let $\FF \subset \CC(\fromto{y_{s+1}}{y_n})$
   be the subfield generated by degree $0$ functions in $\CC(\fromto{y_{s+1}}{y_n})$
   and by $\fromto{y_{s+1}}{y_{k-1}}$.

   If $y_k \in \FF$, then there is a unique way to express $\Phi^* y_k$:
   write $y_k= \mu \cdot \nu$, where $\mu \in \CC(\fromto{y_{s+1}}{y_n})$ is a monomial of degree $0$
   and $\nu$ is a monomial in $\fromto{y_{s+1}}{y_{k-1}}$.
   Then $\Phi^* \mu = \varphi^* \mu$ and $\Phi^* \nu$ is already fixed.
   So set
   \[
     \Phi^* y_k=\varphi^* \mu \cdot \Phi^* \nu.
   \]

   Similarly, if ${y_k}^r \in \FF$ for some $r > 0$, then assume $r$ is minimal such $r$
   and again write ${y_k}^r= \mu \cdot \nu$,
   where $\mu \in \CC(\fromto{y_{s+1}}{y_n})$ is a monomial of degree $0$
   and $\nu$ is a monomial in $\fromto{y_{s+1}}{y_{k-1}}$.
   Then set
   \[
     \Phi^* y_i= \sqrt[r]{\varphi^* \mu \cdot \Phi^* \nu}.
   \]

   Otherwise, if ${y_k}^r \notin \FF$ for any $r >0$, then we have complete freedom to
   choose $\Phi^*y_k$ to be any homogeneous multi-valued section we like.
   For instance, we may fix $\Phi^* y_k =1$.

   Proceeding by induction, we eventually fix all
 $\fromto{\Phi^*y_1}{\Phi^*y_n}$
   and hence define the multi-valued map $\Phi\colon \CC^m \multito \CC^n$.

   We must now show that $\Phi$ defined above indeed describes $\varphi$.
   Firstly, we observe $\Phi$ satisfies homogeneity condition \ref{item:homogeneity_cond_monomial}:
   Let $\mu \in \CC(\fromto{y_{s+1}}{y_n})$ be a monomial of degree $0$.
   Assume there is a nontrivial contribution of $y_k$ in $\mu$
   and there is no contribution of $y_i$ for $i>k$.
   Then
   \[
     {y_k}^r = \mu \cdot \nu
   \]
   where $\nu$ is a monomial in $\fromto{y_{s+1}}{y_{k-1}}$.
   By our construction:
   \[
     \left(\Phi^* y_k\right)^r = \varphi^* \mu \cdot \Phi^*\nu.
   \]
   Therefore
   \begin{equation}\label{eq:existence_of_description}
     \Phi^* \mu = \frac{\left(\Phi^* y_k\right)^r}{\Phi^*{\nu}} = \varphi^* \mu.
   \end{equation}
   In particular $\Phi^*\mu$ is homogeneous of degree $0$,
   so homogeneity condition \ref{item:homogeneity_cond_monomial} holds.

   Also the locus $y_1 = \ldots = y_s=0$ is the (non-empty) toric stratum containing $\varphi(X)$,
   so $\Phi$ satisfies the relevance condition of Definition~\ref{defin:homogrel}.

   Finally, by Equation \eqref{eq:existence_of_description}
   the two ring homomorphisms $\Phi^*$ and $\varphi^*$ agree
   so by Theorem~\ref{thm:Phi_homog_and_relevant_is_a_description}
   indeed $\Phi$ describes $\varphi$.
\end{prf}

The descriptions obtained by following the algorithm of this proof are
not the favoured ones we discussed in the introduction.
For instance for $\varphi= \id_{\PP^1}$ we get
\begin{align*}
  \Phi\colon  \CC^2  & \to \CC^2\\
         [x_1,x_2] & \mapsto [1, \frac{x_2}{x_1}],
\end{align*}
and for the embedding $\varphi\colon \PP^1 \hookrightarrow \PP(1,1,2)$
of \S\ref{example:coordinate_axis_on_P112},
we get
\begin{align*}
  \Phi\colon  \CC^2  & \to \CC^3\\
         [x_1,x_2] & \mapsto [1,0, \frac{x_2}{{x_1}^2}].
\end{align*}
We explain how to modify the descriptions obtained here in
\S\ref{sect:existence_of_complete_description}.

\subsection{The agreement locus revisited}
\label{sect:agreement_revisited}

In this section we calculate the agreement locus for any description.

\begin{prop}\label{thm:agreement_locus}
 Let $\Phi$ be a description of $\varphi$.
 Then
 \[
   \agree{\Phi}{\varphi} = \Reg{\Phi} \setminus \Bigl(\irrel X  \ \cup \ \Phi^{-1}\bigl(\irrel Y\bigr)\Bigr).
 \]
\end{prop}

\begin{prf}
 By the definition of the agreement locus,
 if $\xi \in \agree{\Phi}{\varphi}$, then
 \[
    \xi \in \Reg \Phi  \setminus \irrel X.
 \]
 The homogeneity condition holds for $\Phi$, so, for such $\xi$,
 $\Phi(\xi)$ is contained in a single orbit
 by condition~\ref{item:homogeneity_cond_points} of
 Proposition~\ref{prop:equivalence_of_homog_conds}.
 Since $\pi_Y(\Phi(\xi))$ is defined it follows that no point in
 $\Phi(\xi)$ is in $\irrel Y$,
 which proves the first inclusion:
 \[
   \agree{\Phi}{\varphi} \subset \Reg{\Phi} \setminus
   \Bigl(\irrel X  \ \cup \ \Phi^{-1}\bigl(\irrel Y\bigr)\Bigr).
 \]

 To prove the other inclusion, take
 $\xi \in \Reg{\Phi} \setminus \Bigl(\irrel X  \ \cup
 \ \Phi^{-1}\bigl(\irrel Y\bigr)\Bigr)$
 and set $y = \pi_Y(\Phi(\xi)) \in Y$.
 We must prove, that $x=\pi_X(\xi) \in \Reg \varphi$ and that $\varphi(x) = y$,
 in other words that $\varphi^*$ maps the local ring $\ccO_{Y, y} \subset \CC(Y)$
 into the local ring $\ccO_{X, x} \subset \CC(X)$.
 So take any $q \in \ccO_{Y, y}$.
 By the proof of Theorem~\ref{thm:Phi_homog_and_relevant_is_a_description},
 \[
   \varphi^*q = \Phi^*q \quad\text{as elements of $\CC(X)$}.
 \]
 Since a lift of $y$ to $\CC^m$ is in the image of $\Phi$,
 it follows that $\Phi^*q$ is defined and hence $\varphi^*q$ is defined.
 Hence we can calculate:
 \[
   (\varphi^*q)(x) = (\Phi^*q)(\xi) = q(\Phi(\xi)) = q(y),
 \]
 where the outer equalities hold because rational functions can
 be evaluated on any representative of a point in the Cox cover.
 Since $q$ is regular at $y$, also $\varphi^*q \in \ccO_{X,x}$ as claimed.
 So $\varphi(x) = y$ and thus $\xi \in \agree{\Phi}{\varphi}$.
\end{prf}

\begin{cor}\label{cor:pi_arg_is_open}
  The agreement locus $\agree{\Phi}{\varphi}$ is an open $G_X$-invariant
  subset of $\CC^m$ (and of $\Reg \Phi$).
  In addition, if $X$ is $\QQ$-factorial, then $\pi_X(\agree{\Phi}{\varphi})$ is open.
  In general, $\pi_X(\agree{\Phi}{\varphi})$ contains an open dense subset of $X$.
\end{cor}

\begin{prf}
  $\Reg \Phi$ is an open $G_X$-invariant subset by
  Proposition~\ref{prop:regularity_locus_is_affine}. $\irrel X$ is clearly
  closed and $G_X$-invariant. Finally, $\irrel Y$ is a $G_Y$-invariant subset
  of $\CC^n$, so by homogeneity condition~\ref{item:homogeneity_cond_points}
  also $\Phi^{-1}\bigl(\irrel Y\bigr)$ is $G_X$-invariant, and it is closed in
  $\Reg \Phi$ by Proposition~\ref{prop:preimage_of_closed_is_closed}. Thus
  $\agree{\Phi}{\varphi}$ is open and $G_X$-invariant by
  Proposition~\ref{thm:agreement_locus}.
\end{prf}

The definition of the agreement locus gives
$\agree{\Phi}{\varphi} \subset {\pi_X}^{-1}(\Reg \varphi)$.
In \S\ref{sec:complete_descriptions}, we distinguish those
descriptions for which the equality holds.
In the meantime, we  call the difference between the two sets
the \textbf{disagreement locus}.

\begin{prop}\label{thm:disagreement_is_of_codimension_1}
	Let $\varphi\colon X \dashrightarrow Y$ be a rational map between
	two toric varieties $X$ and $Y$ with a description
	$\Phi\colon \CC^m \multito \CC^n$.
	Consider two open subsets $U_2 \subset U_1$ of $\CC^m$:
	\[
		U_1 = {\pi_X}^{-1}(\Reg \varphi) \quad\text{and}\quad
		U_2 = \agree{\Phi}{\varphi}.
	\]
	The disagreement locus $D=U_1 \setminus U_2$ is 
	a closed subset in $U_1$ purely of codimension~$1$ in $U_1$ or is empty.
\end{prop}

\begin{prf}
   Since $U_2$ is a non-empty open subset of $U_1$
   by Proposition~\ref{thm:agreement_locus} (it is an intersection of three open subsets),
   clearly $D$ is a proper closed subset in $U_1$.
   By Proposition~\ref{thm:agreement_locus} we have an equality
   $
     U_2 = \Reg{\Phi} \setminus \Bigl(\irrel X  \ \cup \ \Phi^{-1}\bigl(\irrel Y\bigr)\Bigr).
   $
   Note that $\irrel X$ is disjoint from $U_1$ (because $\pi_X$ is not regular on $\irrel X$).
   Therefore
   \begin{align*}
     D & =  \underbrace{\Bigl(U_1 \setminus  \Reg{\Phi} \Bigl)}_{=:\Dind}
            \cup
            \underbrace{\Bigl(U_1 \cap \Phi^{-1}\bigl(\irrel Y \bigr)\Bigl)}_{=:\Dirrel}.
   \end{align*}
    By Proposition~\ref{prop:regularity_locus_is_affine} the locus $\Dind$ is indeed purely of codimension~$1$ (or empty).
    It therefore remains to prove that also $\Dirrel$ is purely of codimension~$1$ or empty.

    Assume $\Dirrel$ is not empty and choose arbitrary $\xi \in \Dirrel$.
    We have to prove the codimension of $\Dirrel$ at $\xi$ is $1$.
    Since $\xi \in U_1$ the rational map $\varphi$ is regular at $x = [\xi]$.
    Consider $y = \varphi(x)$ and its toric open affine neighbourhood $V \subset Y$,
    such that $V$ is given by non-vanishing of certain coordinates, say
    \[
      V = \setfromto{y_1 \ne 0}{y_k \ne 0} = \set{y_1 \dotsm y_k \ne 0}.
    \]
    Set $\mvsa=\Phi^* (y_1 \dotsm y_k)$.
    By homogeneity condition~\ref{item:homogeneity_cond_all}, there exists
    $f \in \CC[\Reg {\Phi}]$ such that $\mvsa^r = f$ for some $r \ge 1$.
    We claim that $f(\xi)=0$ and that for all $\xi'$ in the locus $\set{f=0}$
    and in some sufficiently small open neighbourhood of $\xi$
    we have $\xi' \in \Dirrel$.

    First we prove $f(\xi) =0$.
    Since $\xi \in \Phi^{-1}\bigl(\irrel Y \bigr)$
    it follows $\Phi(\xi)$ and $\irrel Y$ have non-empty intersection.
    As usual, since $\Phi(\xi)$ is contained in a single torus orbit
    by the homogeneity condition~\ref{item:homogeneity_cond_points},
    we have $\Phi(\xi) \subset \irrel Y$.
    In particular, $\Phi(\xi)$ is disjoint from ${\pi_Y}^{-1} (V)$,
    in other words the product $y_1 \dotsm y_k$ vanishes on $\Phi(\xi)$.
    So $\mvsa$ vanishes at $\xi$ and therefore $f$ vanishes on $\xi$.

    We prove further that $\xi' \in \set{f=0}$ implies $\xi' \in \Dirrel$,
    at least on some neighbourhood of $\xi$.
    More precisely, we take this neighbourhood to be
    $
       (\varphi \circ {\pi_X})^{-1}(V) \cap \Reg \Phi.
    $
    Since $\Phi$ is regular at such $\xi'$:
    \[
      0= f(\xi') = \mvsa^r(\xi') = \Phi^*(y_1 \dotsm y_k)^r (\xi') = (y_1 \dotsm y_k)^r (\Phi(\xi')),
    \]
    so $\Phi(\xi')$ is contained in the locus $y_1 \dotsm y_k =0$.
    Therefore $\Phi(\xi')$ is disjoint from ${\pi_Y}^{-1}(V)$
    and hence the set $\pi_Y(\Phi(\xi'))$ (if non-empty) is not in $V$.
    On the other hand $\varphi(x')$ is contained in $V$ by our choice of open
    neighbourhood of $\xi$.
    We conclude, that $\xi'$ cannot be in the agreement locus $U_2$.
    But $\xi' \in \Reg{\Phi}$ and $\xi' \notin \irrel X$ (again by our choice
    of open neighbourhood of $\xi$).
    Therefore by Proposition~\ref{thm:agreement_locus} there is no other
    possibility than $\xi' \in \Phi^{-1}\bigl(\irrel Y\bigr)$ so that
    $\xi' \in \Dirrel$ as claimed.

    Hence $\Dirrel$  locally near $\xi$ contains a subset $\set{f=0}$ purely of
    codimension~$1$.
    Since the same holds true for every $\xi\in \Dirrel$ and $\Dirrel \ne U_1$,
    we conclude that $\Dirrel$ is purely of codimension~$1$.
\end{prf}

\subsection{Existence of complete descriptions}
\label{sec:complete_descriptions}
\label{sect:existence_of_complete_description}

The map $\Phi(x_1,x_2) = (x_1^3,x_1^2x_2)$ is a description of
the identity map on $\PP^1$. As written, it does not evaluate
automatically at the point $(0,1)\in\PP^1$: that point is not
in the agreement locus. We can modify the description to
increase the agreement locus following the usual argument that
rational maps are defined in (regular) codimension~1.
The divisor $(x_1)$ contains the bad locus, and the components
of $\Phi$ have multiplicities $\nu_0=(3,2)$ along this divisor.
Using the exponent vector $\nu'=(1,0)$ of $x_1$ itself to push
$\nu_0$ down into the span of the
gradings on the Cox ring, $\nu=\nu_0-\nu' = (3,2)-(1,0)=(2,2)$
computes a scaling factor with which to modify $\Phi$: define
\[
\Phi_{\mathrm{new}} = x_1^{-\nu}\cdot\Phi = [x_1,x_2].
\]
The agreement locus of $\Phi_{\mathrm{new}}$ is larger than
that of $\Phi$, and this new description is better
behaved at $(0,1)$.

We use this notion of `complete agreement' to define complete
descriptions, and then apply the argument above to show that
complete descriptions exists.
In Section~\ref{sect:properties_of_descriptions} we prove
a series of additional properties of complete descriptions.

\renewcommand{\theenumi}{\textnormal{\Alph{enumi}}}
\renewcommand{\labelenumi}{\theenumi.}

\begin{defin}
\label{defin:complete_agreement}
A description $\Phi$ of $\varphi\colon X \dashrightarrow Y$ is \textbf{complete}
if it satisfies
 \begin{enumerate}
   \setcounter{enumi}{2}
   \item \label{item:complete_agreement}
           $\agree{\Phi}{\varphi} = \pi_X^{-1} (\Reg \varphi)$.
 \end{enumerate}
\end{defin}

Proposition~\ref{thm:agreement_locus}, together with this definition,
has an immediate corollary.

\begin{cor}\label{cor:regularity_locus_from_complete_description}
  If $\Phi$ is a complete description of $\varphi$, then
  \[
    \Reg \varphi = \pi_X( \Reg \Phi  \setminus \Phi^{-1}(\irrel Y )).
  \]
  In particular $\varphi$ is regular on $X$ if and only if
  $\Phi$ is regular on $\CC^m \setminus \irrel X$
  and $\Phi^{-1}(\irrel Y)$ is contained in $\irrel X$.
\end{cor}
If $X$ is not a product with $\CC^*$ as one of factors, then saying
$\Phi$ is regular on $\CC^m\setminus\irrel X$ is equivalent to saying
$\Phi$ is regular on $\CC^m$ (because $\irrel X$ is of codimension at least 2), in which case the regularity criterion
for $\varphi$ is the natural statement one would expect,
analogous to standard statement for maps between projective spaces.

The main claim of this article is that complete descriptions always exist
and that they have the properties listed in \S\ref{sect:Motivating_example}.
We establish the properties later in Section~\ref{sect:properties_of_descriptions}.
First we prove the existence.

Let $\Phi$ be a description of a rational map of toric varieties
$\varphi\colon X \dashrightarrow Y$. If $Y$ is a projective space and
$\Phi$ is single-valued, then the procedure for computing a complete
description of $\Phi$ is well known: first clear the denominators
in the sequence $\fromto{\Phi^* y_1}{\Phi^*y_n}$ and then divide
through by the GCD of the resulting polynomials. The proof of our
existence theorem imitates this.

\begin{thm}\label{thm:existence_of_complete_descriptions}
  Let $\varphi\colon X \dashrightarrow Y$ be a rational map of toric varieties.
  Then there exists a complete description $\Phi\colon\CC^m \multito \CC^n$ of $\varphi$.
\end{thm}

Before we start the proof, we discuss the freedom that we have in
choosing a description of a rational map. Let $\Phi$ be a description
of a rational map of toric varieties $\varphi\colon X \dashrightarrow Y$.
If $f \in \cox X$ and $w=(\fromto{w_1}{w_n})$ is a rational linear
combination of $\CC^*$-weights of $Y$, then we can define a multi-valued map
\begin{align*}
    f^w \cdot \Phi\colon \CC^m &\multito \CC^n\\
    x  &  \mapsto \left(\fromto{f^{w_1}\Phi^*y_1}{f^{w_n}\Phi^*y_n}\right)
\end{align*}
which describes the same map $\varphi$ (this follows easily from the proof of
Theorem~\ref{thm:Phi_homog_and_relevant_is_a_description}). Of course,
if $\Phi^* y_i=0$ for some $i$, then there is no harm in replacing the
$i$-th coordinate of $w$ with an arbitrary rational number.

More precisely, we consider the $n$-tuple $w$ as an element of
$\raylattice{Y} \otimes \QQ \simeq \QQ^n$. We define a map of vector spaces $L$,
whose kernel describes the freedom of taking $w$.
Since, by Proposition~\ref{prop:description_satisfy_homog_and_relev},
$\Phi$ satisfies the relevance condition \ref{item:rel_cond_cones} of Definition~\ref{defin:homogrel},
there is a smallest cone $\sigma \in \Sigma_Y$ which contains all the rays
whose corresponding Cox generators $y_i$ lie in $\ker\Phi^*$.
By Corollary~\ref{cor:torus_orbit_containment}, we may assume that
\[
  \Phi^*y_i = 0  \iff \rho_i \in \sigma,
\]
modifying $\Phi$ if necessary.

Let $\starfan{\sigma}$ be the star of $\sigma$, namely the subfan of $\Sigma_Y$
comprising those cones that contain $\sigma$ and their faces.
This fan corresponds to the smallest invariant open neighbourhood
 of the toric stratum containing $\varphi(X)$.
Let $\quofan$ be the quotient fan of $\starfan{\sigma}$ by $\sigma$;
this is the fan of the toric stratum containing $\varphi(X)$ regarded as
a toric variety in its own right.
(If $\varphi(X)$ is not contained in any toric stratum
of $Y$, then both $\starfan{\sigma}$ and $\quofan$ are equal to $\Sigma_Y$.)
Let $L$ be the natural map from $R_Y\otimes\QQ$ to the ambient
rational vector space of $\quofan$ (the composition of the ray lattice
map $\raylattice{Y} \to N_Y$ and the quotient map).
This fits into a diagram of lattices as follows.
\[
\xymatrix{
R_Y \ar[r] \ar[rd]_{L}&  N_Y \ar[d] & \text{dual to}  & M_Y \\
 & \overline{N_{Y(\sigma)}} &\text{dual to}& \overline{M_{Y(\sigma)}} \ar@{^(->}[u]}
\]
where $\overline{N_{Y(\sigma)}}$ is the lattice containing the
quotient fan $\quofan$.

\begin{lemma}\label{lemma:agreement_loci_equal_away_of_f=0}
For any $w\in \ker L$ and nonzero $f\in\cox X$,
both $f^w \cdot \Phi$ and $\Phi$ describe the same map
$\varphi\colon X\dashrightarrow Y$.
Moreover, the agreement locus of the two descriptions
is equal away from the locus $\set{f=0}$: that is,
\[
    \agree{\Phi}{\varphi} \setminus \set{f=0} =
        \agree{f^w \cdot \Phi}{\varphi} \setminus \set{f=0}.
\]
\end{lemma}

\begin{prf}
That the two multi-valued maps describe the same map $\varphi$ follows from
the above considerations: the kernel of the ray lattice map gives the
freedom to choose a linear combination of $\CC^*$-weights, whereas the
pullback of the kernel of the quotient map reflects the freedom to
multiply $0$ coordinates in the description $\Phi$ by anything.

By Proposition~\ref{thm:agreement_locus},
\begin{align*}
    \agree{\Phi}{\varphi} &= \Reg{\Phi}
    \setminus \Bigl(\irrel X  \ \cup \ \Phi^{-1}\bigl(\irrel Y\bigr)\Bigr)
	\quad\text{and} \\
    \agree{f^w \cdot \Phi}{\varphi} &= \Reg{(f^w \cdot \Phi)}
    \setminus
    \Bigl(\irrel X  \ \cup \ (f^w \cdot \Phi)^{-1}\bigl(\irrel Y\bigr)\Bigr).
\end{align*}
Clearly $\Reg{\Phi}$ and $\Reg{(f^w \cdot \Phi)}$ are equal away from $\set{f=0}$,
and also $\irrel X$ does not depend on $\Phi$. Therefore it remains to
compare $\Phi^{-1}\bigl(\irrel Y\bigr)$ with
$(f^w \cdot \Phi)^{-1}\bigl(\irrel Y\bigr)$.

Let $A$ be an irreducible component of $\irrel Y$ defined by the vanishing
of some coordinates, without loss of generality the coordinates
$\fromto{y_1}{y_s}$. Now, for $\xi\in\Reg\Phi$,
\[
    \xi \in \Phi^{-1}( A ) \ \text{ if and only if }\
    \Phi^*y_1 (\xi) = \cdots = \Phi^* y_s (\xi) = 0
\]
whereas, for $\xi\in\Reg(f^w\cdot\Phi)$,
\[
    \xi \in (f^w \cdot \Phi)^{-1}( A ) \ \text{ if and only if }\
    (f^{w_1}\Phi^*y_1) (\xi) = \cdots = (f^{w_s}\Phi^* y_s) (\xi) = 0.
\]
Therefore $\Phi^{-1}(A)$ and $(f^w \cdot \Phi)^{-1}(A)$ are equal
away from $\set{f=0}$, as claimed.
\end{prf}

Now we are ready to prove the theorem.

\begin{prf}[ of Theorem~\ref{thm:existence_of_complete_descriptions}]
By Theorem~\ref{thm:existence_of_description} there is a description
$
    \Phi\colon\CC^m \multito \CC^n
$
of $\varphi$. By Proposition~\ref{thm:disagreement_is_of_codimension_1},
the disagreement locus
$
    D = \pi_X^{-1} (\Reg \varphi) \setminus \agree{\Phi}{\varphi}
$
is a union of codimension $1$ components. If $D$ is empty, then the
theorem is proved, so suppose it is not empty; we must modify $\Phi$
so that the new description is defined on those components which
cover the locus where $\varphi$ is defined.

Choose any homogeneously prime component of $D$ and pick a homogeneously
irreducible polynomial $f\in \cox X$ that vanishes along it. We aim
to replace $\Phi$ by $f^w \cdot \Phi$ for some vector $w$ so that
$\agree{f^w \cdot \Phi}{\varphi}$ contains a general point of $\set{f=0}$.

\paragraph{Step 1: interpret disagreement in terms of a fan.}
Let $v_i\in\QQ$ be the multiplicity of $f$ in $\Phi^*y_i$ and
consider $v = (\fromto{v_1}{v_n})$ as a point in $R_Y\otimes\QQ$,
where $R_Y$ is the ray lattice of $Y$.

\begin{lemma}\label{lem:vanishing_order_along_divisor}
Let $m$ be an integral linear form on the lattice containing $\quofan$, and $\chi^m$ be
the corresponding rational function on $Y$. Then the order of
vanishing of $\varphi^*\chi^m$ along the divisor $(f)$ is equal to
$\left\langle L(v), m\right\rangle$. In particular, $L(v)$ is an
integral point in the lattice of $\quofan$.
\end{lemma}

  \begin{prf}
    $L^*m$ is the monomial expressed in terms of Cox coordinates of $Y$.
    So
    $
       \varphi^* \chi^m = \Phi^* \chi^{L^*m}.
    $
    Now the order of $\Phi^* y_i= \Phi^* \chi^{e_i}$ along $(f)$ is by definition
    $v_i = \left\langle v, e_i\right\rangle$,
    so  the order of $\Phi^*\chi^{L^*m}$ along $(f)$ is
    $
     \left\langle v, L^*m\right\rangle =  \left\langle L(v), m\right\rangle.
    $
  \end{prf}

\begin{cor}\label{cor:Lv_in_support_of_fan}
If $L(v)$ is not in the support of $\quofan$, then $\varphi$ is
not regular on $(f)$.
\end{cor}

  \begin{prf}
    Let $\tau$ be any cone in $\quofan$.
    Since $L(v) \notin \tau$, there exists $m_{\tau} \in \tau^{\vee}$ such that
    $\left\langle L(v), m_{\tau}\right\rangle < 0$.
    Then by Lemma~\ref{lem:vanishing_order_along_divisor}
    the rational function $\varphi^* \chi^{m_{\tau}}$ has a pole along $(f)$.
    Let $U_{\tau}$ be the affine open subset corresponding to a cone in $\starfan{\sigma}$,
    which maps to $\tau$.
    Note that the collection of such $U_{\tau}$ for all $\tau \in \quofan$
    will cover the image of $\varphi$.
    By Proposition~\ref{prop:algebraic_description_of_regularity_locus},
    this implies that $\varphi$ is not regular on $(f)$.
  \end{prf}

Thus if $L(v)$ does not lie in the support of $\quofan$, then $(f)$
is not part of the disagreement locus, contradicting our initial setup.
In short, we may assume that $L(v)$ lies the support of $\quofan$.

\paragraph{Step 2: modify $\Phi$.}
Let $\tau_{quo}$ be the cone in $\quofan$ of minimal dimension that contains $L(v)$,
and $\tau_{star}$ be a cone in $\starfan{\sigma}$ that maps exactly onto $\tau_{quo}$ and is maximal with
this property.

By definition of $\tau_{star}$, there is a vector $u\in \tau_{star}$ that maps
to $L(v)$, and so by choosing a vector $v'$ of $R_Y\otimes\QQ$
in the hyperplane quadrant above $\tau_{star}$ which maps to $u$,
we have $v-v' \in \ker L$.
We may assume that the coordinates of this vector $v'=(\fromto{v_1'}{v_n'})$
satisfy $v_i'=0$ if the $i^{\text{th}}$
ray of $\Sigma_Y$ is not in $\tau_{star}$ and otherwise $v_i'\ge 0$.
We define
\[
    \Phi_{\mathrm{new}}:= f^{v' - v} \cdot \Phi.
\]
By Lemma~\ref{lemma:agreement_loci_equal_away_of_f=0} the two
descriptions of $\varphi$ have the same (dis)agreement locus away
from $\set{f=0}$.

\paragraph{Step 3: $\agree{\Phi_{\mathrm{new}}}{\varphi}$ contains a
general point of $\set{f=0}$.}

By Proposition~\ref{thm:agreement_locus}, it is enough to prove the
following two statements:
\begin{itemize}
    \item $\Phi_{\mathrm{new}}$ is regular on a general point of $(f)$.
    \item $\Phi_{\mathrm{new}}$ does not map general point of $(f)$
	into the irrelevant locus of $Y$.
\end{itemize}
The first is immediate: $f^{-v}\cdot\Phi$ is regular along $(f)$, since
$f$ does not appear in any component $\Phi^*y_i$, and
as each component $v_i'$ of $v'$ is non-negative, $\Phi_{\mathrm{new}}$
is also regular there. Moreover, this shows that if $x \in \set{f=0}$ is a
general point, then $\Phi_{\mathrm{new}}(x)$ has zero $y_i$-coordinate
if and only if either $\Phi^*y_i=0 $ or $v'_i > 0$. In particular,
if the $i$-th ray of $\Sigma_Y$ is not in $\tau_{star}$, then
$\Phi_{\mathrm{new}}(x)$ has non-zero $i$-th coordinate. This means
that the standard generator of $B_Y$ determined by $\tau_{star}$ is nonzero
at $\Phi_{\mathrm{new}}(x)$, and so $\Phi_{\mathrm{new}}(x)$ is not in
the irrelevant locus of $Y$. Therefore $\agree{\Phi_{\mathrm{new}}}{\varphi}$
contains a general point of $\set{f=0}$ as claimed.

Thus we have obtained a description $\Phi_{\mathrm{new}}$ of
$\varphi$ whose disagreement locus contains one component less
than that of $\Phi$. Continuing inductively, we obtain a description
with an empty disagreement locus, namely a complete description.
\end{prf}

\begin{example}\label{ex:complete_descriptions_are_not_unique}
Complete descriptions are not unique. For example, take $X$ to be $\CC$ with
coordinate $x$ and $Y$ to be the non-$\QQ$-factorial base of the standard flop
from \S\ref{sect:using_descriptions}: namely $\cox Y = \CC[y_1,y_2,y_3,y_4]$
graded by $\ZZ$ in degrees $(1,1,-1,-1)$.
Then the map $X \longrightarrow Y$ given by
$[x] \longmapsto [x^t,x^t,x^{1-t},x^{1-t}]$
is a complete description for any rational $t$ in the interval $[0,1]$.

If the target is $\QQ$-factorial, and the map is regular in
codimension~1, then a complete description is unique up to
multiplication by scalars using the whole group action, but we
do not use this fact. (In this example, the values $t=0$ and $1$ define
maps to the two $\QQ$-factorialisations of the cone; other values of
$t$ do not satisfy the relevance condition.)
\end{example}

\section{Geometry of descriptions}
\label{sect:properties_of_descriptions}

In this section, we prove that images and preimages of subschemes behave as
well as the first examples could allow, and we compute descriptions of
compositions of maps, where composition makes sense. We work throughout with a
rational map $\varphi$ together with a description $\Phi$ (not necessarily a
complete description, unless explicitly mentioned) as in the diagram
\[
   \xymatrix{
       \CC^m \xymultito{r}^{\Phi}
       \ar@{-->}[d]^{\pi_X} & \CC^n\ar@{-->}[d]^{\pi_Y}\\
	    X \ar@{-->}[r]^{\varphi} &  Y }
\]

From the start we insisted that descriptions should behave
well when pulling back Cartier divisors. We prove this `local Cartier
pullback' property now, and then present a few additional conditions below
that are closely related to the complete agreement property that characterises
complete descriptions.

\subsection{Properties D--F of complete descriptions}
\label{sect:DEF}

Recall from \S\ref{sect:ring_extensions}:
if $\mvsb$  is a homogeneous multi-valued section in the field of fractions of $\mapring{\Phi}$,
then $\lfloor \mvsb \rfloor$ and $\lceil \mvsb \rceil$ are both homogeneous (single-valued) sections in
$\coxfield X$.

\begin{prop}
\label{prop:very_complete_description_Cartier}
Let $D = (f)$ be a Weil divisor on $Y$, for some $f \in \coxfield{Y}$,
whose support does not contain $\varphi(\Reg \varphi)$.
Consider an open subset $V\subset Y$ for which $D|_{V}$ is Cartier.
Denote the interior of $\pi_X(\agree{\Phi}{\varphi})$ by $\agr \subset X$ and
let $U=\varphi^{-1}(V) \cap \agr$.
Write $\Phi^*f = \lceil \Phi^*f \rceil \cdot \mvsa$
for some homogeneous multi-valued section $\mvsa$ on $X$.

Then $\mvsa$ is invertible on $\inv{\pi_X} (U)$ and
the Cartier divisor ${\varphi|_{U}}^*(D|_{V})$ on $U$ is equal to the
restriction $E|_U$, where $E= (\lceil \Phi^*f \rceil)$
denotes the divisor on $X$ defined by $\lceil \Phi^*f \rceil$.
\end{prop}

Note that if $\Phi$ is a complete description, then $\agr = \Reg \varphi$
and so $U = \varphi^{-1}(V)$.
Also if $D$ is a Cartier divisor on $Y$, then we may take $V=Y$.
Thus, if both of these hold, the statement of the proposition has a much easier form;
see condition~\ref{item:complete_description_Cartier} below.
We observe the following lemma before we prove the proposition.

\begin{lemma}\label{lemma:floor_is_invertible}
  Let $\mvsb$ be a homogeneous multi-valued section in the field of fractions of $\mapring{\Phi}$.
  If $W \subset \CC^m$ is an open subset on which $\mvsb$ is invertible,
  then $\lfloor \mvsb \rfloor, \lceil \mvsb \rceil \in \cox{X}$ are also invertible on~$W$.
\end{lemma}

\begin{prf}
  By definition, $\mvsb = \sqrt[r]{g}$ is invertible on $W$
  if and only if $g\in \coxfield{X}$ is invertible on~$W$.
  For some (reduced) $f \in \cox X$ the locus $Z =\set{f =0}$
  is the codimension $1$ locus of $\CC^m \setminus W$,
  so that $\ccO_{\CC^m}(W) = \cox{X}[\inv f]$.
  Now $g$ is invertible on $W$ if and only if $g, \inv g \in \cox{X}[\inv f]$.
  By Proposition~\ref{prop:simple_extension_under_localisation},
  $\cox{X}[\inv f]  \subset \mapring{\Phi}[\inv f]$
  is a simple ring extension, so $\mvsb,\inv\mvsb \in \mapring{\Phi}[\inv f]$
  by Definition~\ref{def:simple_extension}\ref{item:define-simple-extension-normal}.
  Thus by Proposition~\ref{prop:floor-and-ceiling} both $\lfloor \mvsb \rfloor$
  and $\lceil \mvsb \rceil$ are invertible elements in $\cox{X}[\inv f]= \ccO_{\CC^m}(W)$,
  and so they are both invertible on $W$ as claimed.
\end{prf}

\begin{prf}[ of Proposition~\ref{prop:very_complete_description_Cartier}]
We first work locally on an open subset $V'\subset V$ where
$D|_{V'}$ is principal and defined by $h \in \functionfield{Y}$.
Set $k=h/f$. By construction, $k \in \coxfield{Y}$ is invertible on $\inv{\pi_Y}(V')$.
Suppose $U'= \varphi^{-1}(V') \cap \agr$.
We claim that $\Phi^*k$ is invertible on $W':= {\pi_X}^{-1}(U')$.
To show this, we simply check that $(\Phi^*k)(\xi)$ is nonzero for any $\xi\in W'$.
But $(\Phi^*k)(\xi) = k(\eta)$ for any $\eta\in\Phi(\xi)$, and for such $\eta$
we have $\pi_Y(\eta) =\varphi\circ \pi_X (\xi) \in V'$ so $k(\eta)\not=0$.
Thus $\Phi^* k$ is invertible.
It follows from Lemma~\ref{lemma:floor_is_invertible} that
$\left\lfloor {\Phi^*k}\right\rfloor$ is also invertible on $W'$.

Since $f = h/k$ and $\Phi^* h = \varphi^* h$, hence
$
   \Phi^*f = \dfrac{\varphi^* h}{\Phi^*k}
$
and
$
   \lceil \Phi^*f \rceil = \dfrac{\varphi^* h}{\left\lfloor {\Phi^*k}\right\rfloor}.
$
It then follows from $\Phi^* f = \lceil \Phi^* f \rceil \cdot \mvsa$ that
$
\mvsa = \dfrac{\left\lfloor {\Phi^*k}\right\rfloor}{\Phi^*k},
$
and so $\mvsa$ is invertible on $W'$ and ${\varphi|_{U'}}^* (D|_{U'}) = E|_{U'}$.

The same conclusion is true for any $V' \subset V$ on which $D|_{V'}$ is principal.
Since such $V'$ cover $V$ and the corresponding $U'$ cover $U$,
it follows that $\mvsa$ is invertible on $U$ and ${\varphi|_{U}}^* (D|_{U}) = E|_{U}$
as claimed.
\end{prf}

\begin{defin}
\label{defin:properties_DEFG}
Let $\Phi$ be a description of $\varphi\colon X \dashrightarrow Y$. We recall
the complete agreement property \ref{item:complete_agreement} of
Definition~\ref{defin:complete_agreement} and define some other properties of
$\Phi$:
 \begin{enumerate}
\setcounter{enumi}{2}
\item
\textbf{Complete agreement}:
	$\agree{\Phi}{\varphi} = \pi_X^{-1} (\Reg \varphi)$.
\item
\textbf{Global Cartier pullback}:
\label{item:complete_description_Cartier}
Let $D =(f)$ be a Cartier divisor on $Y$ for some $f \in \coxfield{Y}$ whose
support does not contain $\varphi(\Reg \varphi)$.
Write $\Phi^*f = \lceil \Phi^*f \rceil \cdot \mvsa$
for a homogeneous multi-valued section $\mvsa$ on $X$.
Let $E= (\lceil \Phi^*f \rceil)$ be the divisor on $X$ defined by $\lceil \Phi^*f \rceil$.
Then $\mvsa$ is invertible on $\inv{\pi_X} (\Reg \varphi)$ and
the Cartier divisor ${\varphi}^*D$ on $\Reg \varphi$ is equal to the
restriction $E|_{\Reg \varphi}$.
\item
\textbf{Weil preimage:}
\label{item:complete_description_Weil}
	 If $D=(f)$ is an effective Weil divisor on $Y$ for some function $f \in \cox{Y}$
	 and $\varphi(\Reg \varphi)$ is not contained in the support of $D$, then
	 $\Phi^*f$ is regular on $\Reg \varphi$ and its set-theoretic zero locus
	 agrees with the set $\varphi^{-1}(D)$, the preimage of the support of $D$.
\item
\textbf{Coordinate divisors preimage:}
\label{item:complete_description_coordinates}
	The same as Condition E with $D=(y_i)$ for all $i \in \setfromto 1 n$.
\end{enumerate}
\end{defin}
Note that in Condition~\ref{item:complete_description_Cartier},
if $X$ has no torus factors and $\varphi$ is regular (or at least regular in codimension $1$),
then $\mvsa$ is necessarily a constant in $\CC$.

\begin{prop}\label{prop:Y_Q_factorial_then_max_agr_not_necessary}
Let $\Phi$ be a description of $\varphi\colon X \dashrightarrow Y$.
We have the following implications between the properties
of Definition~\ref{defin:properties_DEFG}:
\begin{center}
\ref{item:complete_description_Weil}
	  $\Longrightarrow$
	  \ref{item:complete_description_coordinates}
	  $\Longrightarrow$
	  \ref{item:complete_agreement}
	  $\Longrightarrow$
	  \ref{item:complete_description_Cartier}.
\end{center}
If, furthermore, $Y$ is $\QQ$-factorial, then
\ref{item:complete_description_Cartier} $\Longrightarrow$
\ref{item:complete_description_Weil},
so that all conditions \ref{item:complete_agreement},
		 \ref{item:complete_description_Cartier},
		 \ref{item:complete_description_Weil},
		 \ref{item:complete_description_coordinates} are equivalent.
\end{prop}

\begin{prf}
  The implication \ref{item:complete_description_Weil}
  $\Rightarrow$ \ref{item:complete_description_coordinates} is clear.

  Assume \ref{item:complete_description_coordinates} holds so that $\Phi^* y_i$
  is regular on $\pi_X^{-1}(\Reg \varphi)$. So, in particular,
  $
    \Reg \Phi \supset \pi_X^{-1}(\Reg \varphi).
  $
  Now  assume (by changing the order of coordinates, if necessary) that
  $y_1, \dotsc, y_s$ define a component of the irrelevant locus of $Y$.
  Then the intersection
  $
   (y_1) \cap \dotsb \cap (y_s)
  $
  of divisors on $Y$ is empty, and so the intersection
  $
   \varphi^{-1}((y_1)) \cap \dotsb \cap \varphi^{-1}((y_s))
  $
 is also empty as a subset of $\Reg \varphi$.
 So the zero locus of $\Phi^*y_1, \dotsc, \Phi^* y_s$
 does not intersect $\pi_X^{-1}(\Reg \varphi)$.
 Proposition~\ref{thm:agreement_locus} now implies that
 $\agree{\Phi}{\varphi} = \pi_X^{-1} (\Reg \varphi)$,
 and so property~\ref{item:complete_agreement} holds.
 The implication \ref{item:complete_agreement}
 $\Rightarrow$ \ref{item:complete_description_Cartier}
 follows from Proposition~\ref{prop:very_complete_description_Cartier},
 with $V = Y$ and $\agr = \Reg \varphi$.

 Finally assume $Y$ is $\QQ$-factorial and
 property~\ref{item:complete_description_Cartier} holds.
 Then, since every Weil divisor is $\QQ$-Cartier,
 property~\ref{item:complete_description_Weil} follows automatically.
\end{prf}

\renewcommand{\theenumi}{\textnormal{\normalfont{(\roman{enumi})}}}
\renewcommand{\labelenumi}{\theenumi}

\subsection{Image of a subscheme}\label{sect:image}

Suppose $A\subset X$ is a closed subscheme defined by a homogeneous ideal
$I_A \ideal \cox X$. We seek the ideal in
$\cox Y$ of the scheme-theoretic image of $A$ under
$\varphi\colon X\dashrightarrow Y$.
Recall the notation
$\phibar|_U(A)$ for the closure of the image of $\varphi|_U(A\cap U)$
where $U\subset \Reg\varphi$ is open.

We define an ideal
$J_A \ideal \cox Y$ by
\[
  J_A:= \homog{\Bigl(
  	(\Phi^*)^{-1} \left( \langle I_A \rangle_{\mapring \Phi}
  \right)\Bigr)},
\]
the homogeneous preimage of the ideal that $I_A$ generates in the map ring.

\begin{thm}\label{thm:image_under_description_agrees}
	Let $\varphi\colon X\dashrightarrow Y$ be a rational map of toric
	varieties with a description $\Phi\colon \CC^m \multito \CC^n$,
	and let $\agr$ be the interior of $\pi_X(\agree{\Phi}{\varphi})$;
	in particular, $\agr\subset\Reg\varphi$.

	Suppose $A\subset X$ is a closed subscheme defined by a homogeneous,
	saturated ideal $I_A \ideal \cox X$, and define $J_A$ as above.
	Then the scheme-theoretic image $\phibar|_{\agr}(A) \subset Y$ and the
	subscheme $B \subset Y$ defined by $J_A$ are equal.
	In particular
	\begin{enumerate}
		\item $B$ is independent of the choice of map ring $\mapring{\Phi}$
			(and of the choice of the saturated ideal $I_A$).
		\item If $\Phi$ is a complete description of $\varphi$, then $\phibar(A)$
			and $B$ are equal.
	\end{enumerate}
\end{thm}

\begin{prf}
Let $V$ be a standard open affine toric
subset of $Y$ given by nonvanishing of some coordinates, say
\[
  V=\set{y \in Y \mid y_i \ne 0 \text{ for every } i \in E}
\]
where $E \subset \setfromto{1}{n}$ is some appropriate subset.
Denoting
$
  {\upsilon} = \prod_{i \in E} {y_i}^{r_i},
$
where $r_i$ are the minimal positive integers such that
$\Phi^* {y_i}^{r_i} \in \CC[\Reg \Phi]$, we set $\ccO_Y(V)$ to be the
homogeneous localisation of $\cox Y$ at $\upsilon$ so that
$V = \Spec \ccO_Y(V)$.
Of course, such open subsets form an open cover of $Y$.
It is enough to prove that $B \cap V=\phibar|_{\agr}(A) \cap V$,
which we do below by comparing their ideals in $\ccO_Y(V)$.

In the first place, suppose $\Phi^*\upsilon =0$.
Then by Corollary~\ref{cor:torus_orbit_containment}
the locus $\phibar(X)$ is disjoint from $V$, so
in this case  we need to prove $B \cap V = \emptyset$.
But $\upsilon \in (\Phi^*)^{-1} \left( \langle I_A \rangle_{\mapring \Phi} \right)$,
so $\upsilon \in J_A$, and indeed $B \cap V = \emptyset$.

So assume that $\Phi^* \upsilon \ne 0$ and consider $\inv \varphi(V) \cap \agr$.
It is an open subset of $X$,
 and thus it has a covering by open affine subsets of $X$.
Any open affine set is the complement of a closed set in codimension~1, so
there exists a finite subset $G \subset \cox X$,
such that $\inv \varphi(V) \cap \agr = \bigcup_{g \in G} X_g$ and each $X_g = X \setminus \Supp(g)$ is affine.
Then $X_g = \Spec \ccO_X(X_g)$ where $\ccO_X(X_g)$ is the homogeneous localisation
of $\cox X$ at $g$.

For any $g \in G$, the following diagram shows the natural
relationships between subrings of a common field $\overline{\coxfield{X}}$
on the left and subrings of $\coxfield{Y}$ on the right.
\begin{equation} \label{eq:diagram_for_scheme_calculations}
   \DiagramForImage{}{}{}{}{}{}{}
\end{equation}

Since $\Phi^*(\upsilon) \ne 0$, it is natural to extend the domain of
$\Phi^*$ to $\cox Y [ {\upsilon}^{-1}]$ (we don't need to specify the precise
subset of $\overline{\coxfield{X}}$ that is the image of these elements).
With that, by Theorem~\ref{thm:Phi_homog_and_relevant_is_a_description},
we have
$
  \Phi^*f = \varphi^*f \quad  \text{for all } f \in \ccO_Y(V).
$
In this sense Diagram \eqref{eq:diagram_for_scheme_calculations} is commutative.
Moreover, since $\varphi(X_g) \subset V$ and
$\inv{\pi_X} (X_g)\subset \agree{\Phi}{\varphi}$,
it follows that $\Phi^*(\upsilon)$ is invertible in $\cox X[g^{-1}]$.

It is enough to prove that
the following two ideals in $\ccO_Y(V)$ are equal:
  \begin{align*}
     I(\phibar|_{\agr}(A) \cap V)&= \bigcap_{g \in G}(\varphi^*)^{-1}\left(\hlocalise{(I_A)}{g}\right) \\
     \text{and}\quad I(B)&=\hlocalise{({J_A})}{{\upsilon}}.
  \end{align*}
  The intersection consists of precisely those functions on $V$
  whose preimage in any $X_g$ is in the ideal of $A$ there---which is
  why it defines the image.

  We redraw  Diagram~\eqref{eq:diagram_for_scheme_calculations}
  marking where each ideal lives:
  \[
    \DiagramForImage{
        \phantom{{}_{\cox X^{-1}}}\hlocalise{(I_A)}{g} \ideal
      }{
        \langle I_A \rangle_{\cox X [g^{-1}]} \ideal
      }{
        \phantom{\langle \rangle_{\cox X [g^{-1}]}}
        I_A \ideal
      }{
        \phantom{{}_{ [g^{-1}]}}
        \langle I_A\rangle_{\mapring \Phi} \ideal
      }{
        \vartriangleright I(\phibar|_{\agr}(A) \cap V),  I(B)
      }{
        \vartriangleright \langle J_A  \rangle_{ \cox Y [{\upsilon}^{-1}]} \phantom{aaaaaaaa}
      }{
        \vartriangleright J_A \phantom{aaaaaaaaaaaaaaaa}
      }
  \]
  The idea of the proof is now straightforward:
  grab an element $q$ in one of the ideals $I(\phibar(A))$ or $I(B)$
  and drag it around Diagram \eqref{eq:diagram_for_scheme_calculations}
  to see that in fact $q$ is also in the other ideal.
  We exploit the ``commutativity'' of the diagram and our choice that
  $\Phi^*(\upsilon)$ is a homogeneous single-valued section which is
  invertible on $X_g$. Here are the details.

  Take $q\in \ccO_Y(V)$.
  Then $q \in  I(B)$ if and only if $q = \tilde q / {\upsilon}^l$
  for some $\tilde q \in J_A$ and $l\in \ZZ$, so:
  \begin{alignat*}{2}
        q \in I(B) & \Iff
                \Phi^*(q \cdot {\upsilon}^l) \in \langle I_A \rangle_{\mapring \Phi} \\
                     & \Iff
                \varphi^*(q) \cdot \Phi^*({\upsilon}^l) \in \langle I_A \rangle_{\mapring \Phi}.\\
 \intertext{%
  Since $\Phi^*({\upsilon}^l) \in \CC[\Reg \Phi]$
  by Corollary \ref{cor:about_K(S)_and_zeta}, we have $\varphi^*(q)\cdot \Phi^*(\upsilon^l) \in \CC[\Reg \Phi]$.
  At this point, our insistence that $\mapring\Phi$ is a simple extension is key.
  By Corollary \ref{cor:generators_of_elimination_ideal_are_simple}
  we can continue the chain of equivalences:
  }
        \dots & \Iff
                \varphi^*(q) \cdot \Phi^*({\upsilon}^l) \in I_A.\\
  \intertext{%
   But $\Phi^*({\upsilon}^l)$  is invertible on each $X_g$, so we continue:}
       \dots   &  \Iff
               \varphi^*(q)  \in \langle I_A \rangle_{\cox X[g^{-1}]}
               			&&\text{for every $g \in G$}.\\
  \intertext{%
   The implication $\Longleftarrow$ above needs a careful explanation,
   as this implication does not hold if $I_A$ is not saturated (as in Example~\ref{example:saturation}, say).
   We postpone the proof of this implication until later,
   meanwhile we continue the series of implications:}
       \dots   &  \Iff
               \varphi^*(q)  \in \hlocalise{(I_A)}{g}
               			&&\text{for every $g \in G$}\\
               &  \Iff
                 q  \in (\varphi^*)^{-1}\left(\hlocalise{(I_A)}{g}\right)
                 		&&\text{for every $g \in G$}\\
               & \Iff
               		q \in I(\phibar|_{\agr}(A \cap V)).
  \end{alignat*}

   It remains to prove the missing implication for $q=\tilde{q}/\upsilon^l$
   as above:
   \[
     \varphi^*(q)  \in \langle I_A \rangle_{\cox X[g^{-1}]}
               		\text{ for every $g \in G$} \Longrightarrow
     \varphi^*(q) \cdot \Phi^*({\upsilon}^l) \in I_A.
   \]
   Let $\hat A \subset \Reg \Phi$
   be the subscheme defined by $\langle I_A \rangle_{\CC[\Reg \Phi]}$.
   Suppose $U_g =\set{g \ne 0} \subset \CC^m$.
   The claim of the implication is that if $\varphi^*(q)$ vanishes on $\hat A \cap U_g$ for all $g \in G$,
   then it vanishes on $\hat A \cap (\Reg \Phi \cap \set{\Phi^* \upsilon \ne 0})$.
   Since $I_A$ is saturated, $\hat A = \overline{\hat A \setminus \irrel{X}}$,
   where the closure is taken in $\Reg\Phi$,
   so it is enough to prove the following inclusion of open subsets:
   \[
     (\Reg \Phi \setminus \irrel X) \cap \set{\Phi^* \upsilon \ne 0} \subset \bigcup_{g\in G} U_g.
   \]
   Suppose $\xi \in (\Reg \Phi \setminus \irrel X) \cap \set{\Phi^* \upsilon \ne 0}$.
   Then $\Phi^* \upsilon (\xi) \ne 0$ and $\upsilon(\Phi(\xi)) \ne 0$, so
   $\pi_Y \circ \Phi(\xi) \in V$.
   In particular $\Phi(\xi)$ is not contained in $\irrel Y$ and by Proposition~\ref{thm:agreement_locus},
   $\xi \in \agree{\Phi}{\varphi}$ and
   $\pi_Y \circ \Phi(\xi) = \phireg \circ \pi_X (\xi)$.
   Thus $\pi_X(\xi) \in \phireginv(V)$ and there exists $g \in G$
   such that $\pi_X(\xi) \in X_g$, so in particular $g(\xi) \ne 0$
   and thus $\xi \in U_g$, as claimed.
\end{prf}

\subsection{Preimage of a subscheme} \label{sect:preimage}

Consider as usual a rational map of toric varieties $\varphi\colon
X\dashrightarrow Y$ with a description $\Phi\colon \CC^m \multito \CC^n$ and
fixed choice of map ring $\Phi^*\colon \cox Y \longrightarrow \mrl(\Phi)$. We
study the problem of finding the preimage of a closed subscheme $B \subset Y$
under $\varphi$. Our main interest is to compute $\phireginv(B)$, the
scheme-theoretic preimage under $\phireg\colon \Reg\varphi\to Y$, but
inevitably the subschemes of $X$ we define are concerned with the closure of
this.

\subsubsection{The regular preimage ideal $J_B$}
Suppose that $B$ is defined by the ideal $I_B \ideal \cox Y$.
We consider a related ideal $J_B \ideal \CC[\Reg\Phi]$ which is
the intersection of the  ideal in $\mapring \Phi$ generated by
$\Phi^* (I_B)$ with $\CC[\Reg \Phi]$:
\[
  J_B= \CC[\Reg \Phi] \cap \langle \Phi^*(I_B) \rangle_{\mapring \Phi}
  	\ideal \CC[\Reg \Phi].
\]
We refer to $J_B$ as \textbf{the regular preimage ideal}.

We check first that the calculation of $J_B$ depends only on $\Phi$
being homogeneous.

\begin{prop}
\label{prop:computations_of_preimage}
 Let $I_B = \left\langle f_1, \dotsc, f_{\beta} \right\rangle \ideal \cox Y$
 be a homogeneous ideal generated by homogeneous sections $f_i$.
 Suppose that $\Phi\colon \CC^m \multito \CC^n$ is a multi-valued map
 satisfying a homogeneity condition~\ref{item:homogeneity_cond}
 (this holds if $\Phi$ is the description of
 some rational map $\varphi\colon X \dashrightarrow Y$).
 Then
 \[
	J_B =
	\bigl\langle
   		\lceil\Phi^*f_1\rceil, \dotsc, \lceil\Phi^*f_{\beta}\rceil
   	\bigr\rangle
 \quad \text{as an ideal of $\CC[\Reg \Phi]$.}
 \]
\end{prop}

\begin{prf}
Follows immediately from the definitions and
Corollary~\ref{cor:generators_of_elimination_ideal_are_simple}.
\end{prf}
Note also that $J_B$ does not depend on the choice of map ring $\mapring{\Phi}$
(see Proposition~\ref{prop:uniqueness_of_floor}).

\subsubsection{Computable preimages}
The relationship between the preimage $\overline{\phireginv(B)}$ and the
regular preimage ideal $J_B$ is a little delicate---we have already seen a
counter-example to an over-optimistic statement in
\S\ref{sect:ideals_of_subvarieties_of_toric_varieties}---and so we identify a
general property which will permit computation of preimages under certain
conditions.

\begin{defin}
Let $\varphi\colon X\dashrightarrow Y$ be a rational map of toric
varieties with a description $\Phi\colon \CC^m \multito \CC^n$.
Fix a closed subscheme $B \subset Y$ with homogeneous defining ideal
$I_B \ideal \cox Y$.
Let $J_B$ be the regular preimage ideal as defined above.

For any open subset $W \subset Y$, we say {\bf $B$ has a
computable preimage on $W$ with respect to $I_B$} (and with respect to $\Phi$
and $\mrl(\Phi)$) if and only if
the subscheme of $X$ defined by $J_B$ equals $\phireginv(B)$
on  $\phireginv(W)$.
\end{defin}

The particular description $\Phi$ we are working with at any time is fixed, so
we do not usually mention $\Phi$. As stated, this property also depends on the
choice of map ring $\mrl(\Phi)$, but this is illusory and we also do not
mention it; see Corollary~\ref{cor:preimage_indep} below.

\begin{thm}\label{thm:preimage_under_description_agrees}
Let $\varphi\colon X\dashrightarrow Y$ be a rational map of toric varieties
with a description $\Phi\colon \CC^m \multito \CC^n$.
Let $B \subset Y$ be a closed subscheme with homogeneous defining ideal
$I_B = \langle \fromto{f_1}{f_{\beta}} \rangle \subset \cox Y$.
If $W \subset Y$ is an open subset on which each divisor
$(f_i)|_{W}$ is Cartier and $\pi_X^{-1}\phireginv(W) \subset \agree{\Phi}{\varphi}$,
 then $B$ has a computable preimage on $W$.

Moreover, on the interior of $\pi_X(\agree{\Phi}{\varphi})$
the scheme defined by $J_B$ is a subscheme of $\phireginv(B)$.
\end{thm}

The main content of this result is that our ability to compute a preimage
for $B$ depends in part on the equations we use to define $B$.

\begin{cor}
\label{cor:preimage_indep}
Let $X,Y,\varphi,\Phi$ and $B$ be as in the theorem.
\begin{enumerate}
\item
	If $\Phi$ is a complete description,
        then the subscheme $B$ has computable preimage on the smooth locus
	$Y_0$ of $Y$.
\item\label{item:free_preimage}
	If $\Phi$ is a complete description and $I_B$ freely defines $B$,
        then $B$ has a computable preimage on $Y$.
\item\label{item:disjoint_B_from_sings}
	If $\Phi$ is a complete description and $\phireginv(B) = \overline{\phireginv(B \cap Y_0)}$
        (which happens for instance, when $B$ is disjoint from the singularities of $Y$),
        then $B$ has a computable preimage on $Y$
\end{enumerate}
\end{cor}
The conditions do not impose requirements on the existence of many Cartier
divisors on $Y$.

The proof of the main part of the theorem is in two steps which we state as
separate lemmas. The first step reduces the theorem to the case where $I_B$
is a principal ideal. In the second we observe that the computable preimage
property holds on the Cartier locus of principal ideals.
The proof of the ``moreover'' part of
the theorem very similar to the proof of Theorem~\ref{thm:image_under_description_agrees},
so we omit it: one would choose suitable open affine covers of $X$ and $Y$,
and prove the appropriate inclusion
of ideals; the calculations may be simplified by observing additivity of ideals
and reducing to the case where $I_B$ is a principal ideal.

\begin{lemma}\label{lem:preimage_is_additive}
  Having a computable preimage is additive in the following sense.
  Let $B_1$ and $B_2$ be two subschemes in $Y$.
  Suppose $W\subset Y$ is an open subset on which both $B_1$ and $B_2$
  have a computable preimage with respect to their defining ideals
  $I_{B_1}, I_{B_2} \subset \cox Y$ respectively.
  Then the closed subscheme $B_1\cap B_2$
  has a computable preimage on $W$ with respect to $I_{B_1}+I_{B_2}$.
\end{lemma}

\begin{prf}
  Let $B = B_1 \cap B_2 \subset Y$. It is enough to prove that
  $J_B = J_{B_1} + J_{B_2}$, since this sum defines the intersection
  of the preimages of $B_1$ and $B_2$ on the open subset $\phireginv(W)$
  (see Lemma~\ref{lem:defining_ideals_are_additive}).
  The equality $J_B = J_{B_1} + J_{B_2}$ follows from Proposition~\ref{prop:computations_of_preimage}:
 for homogeneous ideals $I_{B_1}$, $I_{B_2} \subset \cox Y$,
 \[
   \CC[\Reg \Phi] \cap \langle  \Phi^* (I_{B_1} + I_{B_2}) \rangle_{\mapring \Phi}
   =  J_{B_1} + J_{B_2}.
 \]
\end{prf}

\begin{lemma}\label{lem:preimage_of_principal_agrees}
  If $f \in \cox Y$ is a polynomial and $W\subset Y$ an open subset on which
  the restriction $(f)|_{W}$ of the Weil divisor $(f)$ on $Y$ is Cartier
  and $\pi_X^{-1}\phireginv(W) \subset \agree{\Phi}{\varphi}$,
  then $B$ has computable preimage on $W$ with respect to its defining
  ideal $I_B=\langle f\rangle$.
\end{lemma}

\begin{prf}
  By Proposition~\ref{prop:computations_of_preimage}, $J_B$ is principal and
  generated by $\lceil \Phi^* f \rceil$.
  By Proposition~\ref{prop:very_complete_description_Cartier} we have
  $\Phi^* f = \lceil \Phi^* f \rceil \cdot \mvsa$, where $\mvsa$ is a
  homogeneous multi-valued section invertible on $\pi_X^{-1}\phireginv(W)$.
  Moreover $\lceil \Phi^* f \rceil$ defines the divisor $\varphi^* D$ on
  $\varphi^{-1}(W)$. Since the definition of pullback of a Cartier divisor
  agrees with the definition of preimage of the underlying scheme,
  it follows that $\varphi^{-1}(B)$ is given by the ideal $J_B$ on $W$.
\end{prf}